\documentclass[journal,onecolumn,compsoc]{IEEETran}
\PassOptionsToPackage{hyphens}{url}\usepackage[hidelinks]{hyperref}
\usepackage{cite}
\usepackage{amsmath,amssymb,amsfonts}
\usepackage{algorithmic}
\usepackage{graphicx}
\usepackage{textcomp}
\usepackage{dsfont}
\usepackage{color}
\usepackage{epstopdf}        
\usepackage{enumerate}
\usepackage{lscape}
\usepackage{multicol}
\usepackage{multirow}
\usepackage{calligra}
\usepackage{booktabs}
\usepackage{mathtools}

\usepackage{framed} 
\usepackage{picins}
\usepackage{empheq}
\usepackage{graphics} 
\usepackage{epsfig} 
\usepackage{times} 
\usepackage{amsmath} 
\usepackage{amssymb}  
\usepackage{amsfonts}  
\usepackage{epstopdf}
\usepackage{enumerate}
\usepackage{graphicx} 
\usepackage{color}
\usepackage{comment}

\usepackage{soul}%
\usepackage{cancel}

\usepackage{amsthm}
\newtheorem{thm}{Theorem}
\newtheorem{prop}{Proposition}
\newtheorem{lemma}{Lemma}
\newtheorem{cor}{Corollary}
\newtheorem{definition}{Definition}
\newtheorem{assumption}{Assumption}
\newtheorem{remark}{Remark}
\newtheorem{example}{Example}

\newcommand*{\QEDB}{\hfill\ensuremath{\square}}%
\def\BibTeX{{\rm B\kern-.05em{\sc i\kern-.025em b}\kern-.08em
    T\kern-.1667em\lower.7ex\hbox{E}\kern-.125emX}}

\definecolor{Nicolas}{rgb}{0.017,0.59,0.55}
\newenvironment{SUBENVNicolas}{\color{Nicolas}}{\color{black}}

\definecolor{Jorge}{rgb}{0.217,0.51,0.255}
\newenvironment{SUBENVJorge}{\color{Jorge}}{\color{black}}

\definecolor{Daniel}{rgb}{0.74,0.2,0.64}
\newenvironment{SUBENVDaniel}{\color{Daniel}}{\color{black}}

\allowdisplaybreaks
\begin{document}

\title{Prescribed-Time Control in Switching Systems with Resets: A Hybrid Dynamical Systems Approach}

\author{Daniel E. Ochoa, Nicolas Espitia, Jorge I. Poveda \vspace{-20pt}
\thanks{D. E. Ochoa and J. I. Poveda are with the Department of Electrical and Computer Engineering, University of California, San Diego, La Jolla, CA, 92093, USA. N. Espitia is with the Centre de Recherche en Informatique Signal et Automatique de Lille, CNRS
Centrale Lille, University of Lille, France. Corresponding Author: Jorge I. Poveda ({\tt poveda@ucsd.edu}).}}

\maketitle
\thispagestyle{empty} 

\begin{abstract}
We consider the problem of achieving prescribed-time stability (PT-S) in a class of hybrid dynamical systems that incorporate switching nonlinear dynamics, exogenous inputs, and resets. By ``prescribed-time stability'', we refer to the property of having the states converge to a particular compact set of interest before a given time defined \emph{a priori} by the user. We focus on dynamical systems that achieve this property via time-varying gains. For continuous-time systems, this approach has received significant attention in recent years, with various applications in control, optimization, and estimation problems. However, its extensions beyond continuous-time systems have been limited. This gap motivates this paper, which introduces a novel class of switching conditions for switching systems with resets that incorporate time-varying gains, ensuring the PT-S property even in the presence of unstable modes. The analysis leverages tools from hybrid dynamical system's theory, and a contraction-dilation property that is established for the hybrid time domains of the solutions of the system. We present the model and main results in a general framework, and subsequently apply them to three novel applications: (a) PT regulation of switching plants with no common Lyapunov functions; (b) PT control of dynamic plants with uncertainty and intermittent feedback; and (c) PT decision-making in non-cooperative switching games via hybrid Nash seeking dynamics.
\end{abstract}
%
%

\section{INTRODUCTION}
\IEEEPARstart{R}ecent advances in nonlinear control analysis and design \cite{song2017time,Orlov2022,KRISHNAMURTHY2020108860,KrishnamurthyAdaptive,TRAN2020104605} have reinvigorated the concept of Prescribed-Time  Stability (PT-S), finding successful applications in different domains such as nonlinear regulation \cite{song2017time,Orlov2022}, adaptive control \cite{KRISHNAMURTHY2020108860,KrishnamurthyAdaptive}, systems with delays \cite{Espitia_TAC_2020}, partial differential equations \cite{Espitia2019FTSAutomatica}, and stochastic systems \cite{li2022prescribed}. 

In contrast to asymptotic or exponential stability, the PT-S property guarantees convergence of the trajectories of the system to the desired target (e.g., a particular compact set of interest) before a given time defined \emph{a priori} by the user, and completely independent on the initial conditions. As such, achieving this property requires either time-varying or non-Lipschitz vector fields in the dynamics of the system. Non-Lipschitz autonomous systems that achieve convergence to the point (or set) of interest before a fixed time have been studied in \cite{Fixed_timeTAC,fixedtimeISS,poveda2022fixed}. The state of the art of this technique, usually called ``fixed-time'' (FxT) stability, was recently reviewed in \cite{Song2023}. However, in contrast to this line of research, in this paper we are interested in systems that achieve convergence to the target before a fixed time by employing time-varying gains. Such types of gains have a long history in the context of optimal control and tactical missile guidance systems \cite{slater1973optimal}, and they have received renewed attention during the last few years due to recent breakthroughs in the design and analysis of nonlinear and adaptive controllers in continuous-time systems with convergence properties defined over a \emph{finite time} horizon; see \cite{Song2023} for a recent survey, and the recent works on adaptive systems \cite{song2017time,Orlov2022,KRISHNAMURTHY2020108860,KrishnamurthyAdaptive,TRAN2020104605,todorovski2023practical}, PDEs \cite{Espitia2019FTSAutomatica,STEEVESEJC20203,Espitia_TAC_2020,Irscheid2022}, and systems with delays \cite{ZekraouiEspitiaPerruquetti2021,PrT_delayObserver2021}.
Since this control approach employs a family of ``blow-up'' gains with finite escape times, the solutions of these systems are also defined only on finite-time intervals. For comprehensive discussions on practical applications and implementation strategies to extend the domain of the solutions, as well as the advantages and limitations of PT control, we refer the reader to the recent works \cite{Orlov2022,todorovski2023practical,song2017time,Song2023,Hollowaylinearsystems2019}.

While the study of Prescribed-Time stability properties in continuous-time systems modeled as ordinary differential equations (ODEs) has seen significant progress, PT-S tools for \emph{hybrid systems} that combine continuous-time and discrete-time dynamics have remained mostly unexplored. For example, switching systems with increasing gains were studied in  \cite{gao2019global} using a common Lyapunov function, which, in general, might not exist in switching systems, even if each mode is exponentially stable. Similarly, switched PT-Stable controllers that turn off or ``clip'' the high-gains before the prescribed time is reached were also recently discussed in \cite{orlov2021prescribed}. However, these results consider only one vector field during the convergence phase, and the switching rules can lead to hybrid dynamical systems that are not well-posed in the sense of \cite{bookHDS}. Indeed, to the best of our knowledge, general results on PT-S for switching and hybrid systems, similar to those existing for asymptotic or exponential stabilization \cite{Liu_Tanwani_Liberzon2022}, are still absent in the literature. 
Since switching and hybrid controllers have been shown to provide powerful solutions to complex control \cite{sanfelice2021hybrid,Liberzon_Book}, optimization \cite{OchoaPoveda20LCSS,teel2019first}, and learning problems \cite{poveda2017framework}, there is a clear need for the development of PT-S tools that enable the analysis and design of new algorithms that leverage the advantages of both PT-S and hybrid control. 
In this paper, we show that the PT stability property can be naturally extended to a class of hybrid dynamical systems (HDS) that model nonlinear switching systems with resets, provided the switching signals also incorporate the dynamic effects of the time-varying gains. Specifically, the following are the main contributions of this paper:

\vspace{0.1cm} 
(a) First, we introduce a class of switching signals that preserve the Prescribed-Time Stability property in systems switching between a finite number of PT-S vector fields with exogenous inputs and resets on the states. To derive these conditions, we reformulate the overall dynamical system as a HDS with dynamic gains that induce a suitable time dilation and contraction in the hybrid time domains of its solutions. By leveraging Lyapunov-based constructions for HDS, we show that the original system is PT-Stable provided the switching signal satisfies a ``blow-up'' average dwell-time (BU-ADT) condition. This condition allows a non-linear increase in the number of jumps and switches as the total amount of flow time in the system reaches the prescribed convergence time. To study the effect of exogenous inputs and/or disturbances in the system, we establish PT convergence results via \emph{ISS-like bounds with the Convergence Property}, that parallel those in the literature of Prescribed-Time stability for ODEs \cite[Def. 2]{song2017time}. However, unlike the existing results for ODEs, our convergence bounds, presented in Theorem 1, are written in ``hybrid time'' and they highlight the potentially (asymptotic) stabilizing effect of the resets, as well as the order of the dynamics of the ``blow-up'' gains. To our knowledge, this is the first result that connects the existing tools on Prescribed-Time stability for ODEs \cite{song2017time}, with the setting of HDS \cite{bookHDS}.

\vspace{0.1cm}
(b) Next, we incorporate unstable modes into the switching systems, and we characterize a ``blow-up" average-activation-time (BU-AAT) condition on the amount of time that the system can spend on the unstable modes while preserving the PT-S property. The unstable modes are allowed to have ``blow-up'' time-varying gains with finite-escape times, as well as exogenous inputs and/or disturbances. To study this setting, we construct a hybrid dynamical system with timers, similar in spirit to those considered in \cite{poveda2017framework,Liu_Tanwani_Liberzon2022,yang2014input}, but having increasing evolution rates parameterized by the blow-up gains, which enable faster switching between the stable and unstable modes as the total amount of flow time in the system approaches the prescribed time. A Lyapunov-based construction on a dilated-time scale, and a contraction argument on the hybrid time domains, are used to establish in Theorem 2 a PT-ISS-like result for switched systems with unstable modes.  

\vspace{0.1cm}
(c) To illustrate the applicability of our model and results, we synthesize three different PT-Stable algorithms for the solution of different control and decision-making problems with prescribed-time convergence requirements. First, we consider the problem of PT regulation of input-affine systems, originally studied in \cite{song2017time} for ODEs, and we show in Proposition \ref{Prop:Regulation_Switching_Plants} that existing feedback-based time-varying control laws can be extended to switching systems, provided the switching signals satisfy the proposed BU-ADT constraint. Next, we study the problem of regulation of switching systems using intermittent PT feedback. We consider in Proposition \ref{Prop:Regulation_intermittenUpsilonlants} a class of dynamics with uncertainty, and we propose a new feedback law that extends the results of \cite{song2017time} to the setting of switching systems with on-and-off feedback. Finally,  we consider the problem of prescribed-time Nash equilibrium seeking in switching non-cooperative games using dynamics with resets. We show in Proposition \ref{Prop:Nash_equilibrium}  that such dynamics fit into our model and can be studied using the analytical tools of Theorem \ref{theorem1}, achieving Prescribed-Time stability in switching strongly monotone games with a common Nash Equilibrium. 

The rest of this paper is organized as follows: Section \ref{sec_notation} introduces the notation and some preliminaries on dynamical systems. Sections \ref{analyticalresults} and \ref{sectionISSTheorems} present the main analytical results and the proofs. Section \ref{sectionapplications} presents three different applications, and Section \ref{secconclusions} ends with the conclusions.
\section{PRELIMINARIES}
\label{sec_notation}
Given a closed set $\mathcal{A}\subset\mathbb{R}^n$ and  $z\in\mathbb{R}^n$, we use $|z|_{\mathcal{A}}:=\inf_{s\in\mathcal{A}}\|z-s\|_2$.  A set-valued mapping $M:\mathbb{R}^p\rightrightarrows\mathbb{R}^n$ is outer semicontinuous (OSC) at $z$ if for each sequence $\{z_i,s_i\}\to(z,s)\in\mathbb{R}^p\times\mathbb{R}^n$ satisfying $s_i\in M(z_i)$ for all $i\in\mathbb{Z}_{\geq0}$, we have $s\in M(z)$. A mapping $M$ is locally bounded (LB) at $z$ if there exists an open neighborhood $N_z\subset\mathbb{R}^p$ of $z$ such that  $M(N_z)$ is bounded. The mapping $M$ is OSC and LB relative to a set $K\subset\mathbb{R}^p$ if $M$ is OSC for all $z\in K$ and $M(K):=\cup_{z\in K}M(x)$ is bounded. A function $\gamma:\mathbb{R}_{\ge 0} \to\mathbb{R}_{\ge0}$ is of class $\mathcal{K}$ if it is continuous, strictly increasing, and satisfies $\gamma(0) =0$. A function $\beta:\mathbb{R}_{\geq0}\times\mathbb{R}_{\geq0}\to\mathbb{R}_{\geq0}$ is of class $\mathcal{K}\mathcal{L}$ if it is nondecreasing in its first argument, nonincreasing in its second argument, $\lim_{r\to0^+}\beta(r,s)=0$ for each $s\in\mathbb{R}_{\geq0}$, and  $\lim_{s\to\infty}\beta(r,s)=0$ for each $r\in\mathbb{R}_{\geq0}$. A function $\tilde{\beta}:\mathbb{R}_{\ge 0}\times \mathbb{R}_{\ge 0}\times \mathbb{R}_{\ge0}\to\mathbb{R}_{\ge 0}$ belongs to class $\mathcal{KLL}$ if for every $s\ge 0$, $\tilde{\beta}(\cdot, s, \cdot)$ and $\tilde{\beta}(\cdot, \cdot, s)$ belong to class $\mathcal{KL}$. Throughout the paper, for two (or more) vectors $u,v \in \mathbb{R}^{n}$, we write $(u,v)=[u^{\top},v^{\top}]^{\top}$ to denote their concatenation. We use $\text{diag}\left(\{B_j\}_{j=1}^J\right)$ to denote the block-diagonal matrix obtained from the set of matrices $\{B_j\}_{j=1}^J$. Given a set $\mathcal{O}\subset\mathbb{R}^n$, we use $\mathbb{I}_{\mathcal{O}}(\cdot)$ to denote the indicator function that satisfies $\mathbb{I}_{\mathcal{O}}(x)=1$ if $x\in \mathcal{O}$, and $\mathbb{I}_{\mathcal{O}}(x)=0$ if $x\notin \mathcal{O}$.
\subsection{Switching Systems}
In this paper, we consider switching systems with inputs of the form $\dot{x}=\tilde{f}_{\sigma(t)}(x,u,t)$, where $x(t_0)=x_0$, $x\in\mathbb{R}^n$, $u\in\mathbb{R}^m$ is an exogenous input, and  $\sigma:\mathbb{R}_{\geq0}\to\mathcal{Q}$ is a right-continuous, piecewise constant signal that maps the current time $t$ to a finite set of indexes (or modes) $\mathcal{Q}=\{1,2,\ldots,\overline{q}\}$, $\overline{q}\in \mathbb{Z}_{\ge1}$.  For each $q\in \mathcal{Q}$, $\tilde{f}_q:\mathbb{R}^n\times\mathbb{R}^m\times\mathbb{R}_{\geq0}\to\mathbb{R}^n$ is assumed to be continuous with respect to all arguments. Following the notation of \cite{Liu_Tanwani_Liberzon2022}, we use $\mathcal{S}$ to denote the set of all right-continuous, piecewise constant signals from $\mathbb{R}_{\geq0}$ to $\mathcal{Q}$, with a locally finite number of discontinuities. Such functions are referred to as switching signals. For each signal $\sigma\in \mathcal{S}$, we also define the collection of switching instants $\mathcal{W}(\sigma):=\{t\geq0:\sigma(t)\neq \sigma(t^-)\}$. In this way, the switching system evolves according to
\begin{equation}\label{vectorfield1}
\dot{x}=\tilde{f}_{\sigma(t)}(x,u,t),~~~~\forall~t\notin \mathcal{W}(\sigma),
\end{equation}
where the solutions $x$ to \eqref{vectorfield1} are understood in the Caratheodory sense over any interval $[t_a,t_b)$ where $\sigma$ is constant. During switching times $t\in \mathcal{W}(\sigma)$, we allow ``jumps'' in the state $x$ via mode-dependent \emph{reset maps} of the form
\begin{equation}\label{resetmap1}
x(t)=R_{\sigma(t^-)}(x(t^-)),~~~~\forall~t\in \mathcal{W}(\sigma),
\end{equation}
where the function $R_q:\mathbb{R}^n\to\mathbb{R}^n$ is assumed to be continuous for each $q\in\mathcal{Q}$. Throughout the paper, we will refer to switching systems of the form \eqref{vectorfield1}-\eqref{resetmap1} as R-Switching systems.

In this paper, we consider switching systems with a mix of stable and unstable modes (defined in Section III.b). We denote the set of stable modes as $\mathcal{Q}_s$, the set of unstable modes as $\mathcal{Q}_u$, such that $\mathcal{Q}_s\cup\mathcal{Q}_u=\mathcal{Q}$ and $\mathcal{Q}_u\cap\mathcal{Q}_s=\emptyset$.

\vspace{0.1cm}
\begin{remark}
By taking $R_q$ to be equal to the identity map, system \eqref{vectorfield1}-\eqref{resetmap1} recovers a standard switching system \cite{Liberzon_Book}. However, the use of general maps $R_q$ opens the door to study PT-S results in certain reset control systems \cite{prieur2018analysis} (by taking $\mathcal{Q}=\{1\}$) as well as more general switched reset controllers (when $|\mathcal{Q}|>1$), see \cite{Liu_Tanwani_Liberzon2022}. It is also possible to consider discontinuous functions $f_q,R_q$ by working with their corresponding Krasovskii regularizations \cite[Def. 4.13]{bookHDS}. \QEDB
\end{remark}
\subsection{Hybrid Dynamical Systems with Inputs}
Since switching systems with resets incorporate continuous-time and discrete-time dynamics, they can be modeled as hybrid dynamical systems (HDS). In this paper, we work with HDS aligned with the framework of \cite{bookHDS}, where the state $z\in\mathbb{R}^n$ evolves according to the dynamics:
\begin{subequations}\label{HDS0}
\begin{align}
    &(z,u)\in \tilde{C}:= C\times\mathbb{R}^m,~~~~~~\dot{z}\in F(z,u),\label{HDS0:flow}\\
    &(z,u)\in \tilde{D}:= D\times\mathbb{R}^m,~~~~z^+\in G(z), \label{HDS0:jump}
\end{align}
\end{subequations}
where $u\in\mathbb{R}^m$ is an exogenous input, $F:\mathbb{R}^n\times\mathbb{R}^m\rightrightarrows\mathbb{R}^n$ is called the flow map, $G:\mathbb{R}^n\rightrightarrows\mathbb{R}^n$ is called the jump map, $\tilde{C}\subset\mathbb{R}^n\times\mathbb{R}^m$ is called the flow set, and $\tilde{D}\subset\mathbb{R}^n\times\mathbb{R}^m$ is called the jump set. We use $\mathcal{H}=(\tilde{C},F,\tilde{D},G,u)$ to denote the \emph{data} of the HDS $\mathcal{H}$. HDS of the form \eqref{HDS0} are a generalization of purely continuous-time systems ($\tilde{D}=\emptyset$) and purely discrete-time systems ($\tilde{C}=\emptyset$).
Solutions to system \eqref{HDS0} are parameterized by a continuous-time index $t\in\mathbb{R}_{\geq0}$, which increases continuously during flows, and a discrete-time index $j\in\mathbb{Z}_{\geq0}$, which increases by one during jumps. Therefore, solutions to \eqref{HDS0} are defined on \emph{hybrid time domains} (HTDs). A set $E\subset\mathbb{R}_{\geq0}\times\mathbb{Z}_{\geq0}$ is called a \textsl{compact} HTD if $E=\cup_{j=0}^{J-1}([t_j,t_{j+1}],j)$ for some finite sequence of times $0=t_0\leq t_1\ldots\leq t_{J}$. The set $E$ is a HTD if for all $(T,J)\in E$, $E\cap([0,T]\times\{0,\ldots,J\})$ is a compact HTD. Given a HTD $E$, let
\begin{align*}
\text{sup}_tE&:=\text{sup}\left\{t\in\mathbb{R}_{\geq0}:\exists~j\in\mathbb{Z}_{\geq0},~\text{such that}~(t,j)\in E\right\}.\\
\text{sup}_jE&:=\text{sup}\left\{j\in\mathbb{Z}_{\geq0}:\exists~t\in\mathbb{R}_{\geq0},~\text{such that}~(t,j)\in E\right\}.
\end{align*}
Also, let $\text{sup}~E:=(\text{sup}_tE,\text{sup}_jE)$, and $\text{length}(E):=\text{sup}_tE+\text{sup}_jE$. The following definition is borrowed from \cite{Cai2009}.

\vspace{0.1cm}
\begin{definition}\label{definitionsolutions1}
A hybrid signal is a function defined on a HTD. A hybrid signal $u:\text{dom}(u)\to \mathbb{R}^m$ is called a hybrid input if $u(\cdot, j)$ is Lebesgue measurable and locally essentially bounded for each $j$. A hybrid signal $z:\text{dom}(z)\to \mathbb{R}^n$ is called a hybrid arc if $z(\cdot, j)$ is locally absolutely continuous for each $j$ such that the interval $I_j:=\{t:(t,j)\in \text{dom}(z)\}$ has nonempty interior. A hybrid arc $z:\text{dom}(z)\to \mathbb{R}^n$ and
a hybrid input $u:\text{dom}(u)\to\mathbb{R}^m$ form a solution pair $(z,u)$ to \eqref{HDS0} if $\text{dom}(z)=\text{dom}(u)$, $(z(0, 0), u(0, 0))\in \tilde{C}\cup \tilde{D}$, and:
\begin{enumerate}
\item For all $j\in\mathbb{Z}_{\geq0}$ such that $I_j$ has nonempty interior, and for almost all $t\in I_j$, $(z(t,j),u(t,j))\in \tilde{C}$ and $\dot{z}(t,j)\in F(z(t,j),u(t,j))$.
\item For all $(t,j)\in\text{dom}(z)$ such that $(t,j+1)\in \text{dom}(z)$, $(z(t,j),u(t,j))\in \tilde{D}$ and $z(t,j+1)\in G(z(t,j))$. \QEDB
\end{enumerate}
\end{definition}
\vspace{0.1cm}
\begin{remark}
According to Definition \ref{definitionsolutions1}, solutions to \eqref{HDS0} are required to satisfy $\text{dom}(z)=\text{dom}(u)$. To establish this correspondence, in this paper the hybrid input $u$ in \eqref{HDS0} is obtained from $u$ in \eqref{vectorfield1} by using (with some abuse of notation) $u(t,j)=u(t)$  during the flows \eqref{HDS0:flow} for each fixed $j$, and by keeping $u$ constant during the jumps \eqref{HDS0:jump}. \QEDB 
\end{remark}

\vspace{0.1cm}
A hybrid solution pair $(z,u)$ is said to be maximal if it cannot be further extended. A hybrid solution pair $(z,u)$ is said to be complete if $\text{length}~\text{dom}(z)=\infty$. This does not necessarily imply that $\text{sup}_t\text{dom}(z)=\infty$, or that $\text{sup}_j\text{dom}(z)=\infty$, although at least one of these two conditions should hold when $z$ is complete.  To simplify notation, in this paper we use
\begin{equation}\label{supu}
|u|_{(t,j)}=\sup_{\substack{(0,0)\leq(\tilde{t},\tilde{j})\leq (t,j)\\(t,j)\in\text{dom}(z)}}\left|u(\tilde{t},\tilde{j})\right|,
\end{equation}
and $|u|_{(t,j)}$ is denote by $|u|_{\infty}$ when $t+j\to\infty$ in \eqref{supu}.
\section{PT-ISS IN A CLASS OF HYBRID DYNAMICAL SYSTEMS}
\label{analyticalresults}
Motivated by the PT-S property studied for ODEs \cite{song2017time,Orlov2022,KRISHNAMURTHY2020108860,KrishnamurthyAdaptive,TRAN2020104605,Espitia_TAC_2020}, we consider a sub-class of HDS of the form \eqref{HDS0}, with $z=(\psi,\mu_k)\in\mathbb{R}^n\times \mathbb{R}_{\ge1}$, the $C$ is given by
\begin{subequations}\label{mainHDSmodel}
\begin{equation}\label{flowsetPT}
C:=\Psi_C\times\mathbb{R}_{\geq1},
\end{equation}
the flow map is given by:
\begin{equation}\label{flowshybrid1a}
\dot{z}=\left(\begin{array}{c}
\vspace{0.1cm}
\dot{\psi}\\
\vspace{0.1cm}
\dot{\mu}_k
\end{array}\right)\in F_T(z,u):=\left(\begin{array}{c}
\vspace{0.1cm}
\mu_k\cdot F_{\Psi}(\psi,u,\mu_k),\\
\dfrac{k}{T}\mu_k^{1+\frac{1}{k}},
\end{array}\right),
\end{equation}
where $T>0$ and $k\geq1$ are tunable parameters, and $F_{\psi}:\mathbb{R}^{n}\times\mathbb{R}^m\times \mathbb{R}_{\ge 1}\rightrightarrows
\mathbb{R}^n$ is a set-valued mapping to be specified below. The set $D$ is given by
\begin{equation}\label{jumpSet}
D=\Psi_D\times\mathbb{R}_{\geq1},
\end{equation}
and the jump map is given by:
\begin{equation}\label{jumpmapPT}
z^+=\left(\begin{array}{c}
\psi^+\\
\mu_k^+
\end{array}\right)\in G(z):=\left(
\begin{array}{c}
G_{\Psi}(\psi)\\
\mu_k
\end{array}\right),
\end{equation}
where $G_{\Psi}:\mathbb{R}^n\rightrightarrows\mathbb{R}^n$ is also to be specified. The following regularity assumption will be considered throughout this paper.
\end{subequations}

\vspace{0.1cm}
\begin{assumption}\label{wellposedassumption}
The sets $\Psi_C,\Psi_D\subset\mathbb{R}^{n}$ are closed. The set-valued maps $F_{\psi}$ and $G_{\psi}$ are OSC and LB with respect to $\Psi_C$, and $\Psi_D$, respectively; and $F_{\psi}$ is convex for all $(\psi,u,\mu_k)\in \Psi_C\times\mathbb{R}^m\times\mathbb{R}_{\geq1}$. \QEDB
\end{assumption}

\vspace{0.1cm}
The HDS \eqref{mainHDSmodel} has a cascade structure, where the dynamics of $\mu_k$ are independent of the state $\psi$. However, by the construction of the flow and jump sets, the dynamics of $\psi$ will determine the structure of the HTDs of the solutions of the system, e.g., purely continuous, purely discrete, eventually continuous, eventually discrete, etc. In this paper, we are interested in signals $\mu_k$ that have finite escape times ``controlled'' by the parameters $(T,k)$ and the initial condition $\mu_k(0)$. This property, formalized in the next lemma, can be readily established by direct integration. For completeness, the proof is presented in the Supplemental Material (see Section \ref{sec:appendix}).

\vspace{0.1cm}
\begin{lemma}\label{lemmaodek}
Let $k\geq1$, and consider the ``\emph{blow-up}'' (BU) ODE $\dot{\mu}_k=\frac{k}{T}\mu_k(t)^{1+\frac{1}{k}}$ with $\mu_k(0)=\mu_0\in\mathbb{R}_{\geq1}$. Then, its unique solution satisfies:
\begin{align}\label{tpmu}
\mu_k(t)
=\frac{T^k}{\left(\Upsilon_{T,k}-t\right)^k}\geq 1,~~~~\forall~t\in[0,\Upsilon_{T,k}),
\end{align}
where $\Upsilon_{T,k}:=T\mu_0^{-\frac{1}{k}}$. \QEDB
\end{lemma}

\vspace{0.1cm}
Based on \eqref{tpmu}, for all $k\geq1$ the mapping $\mu_k(\cdot)$ is continuous in its domain, strictly increasing, and satisfies $\lim_{t\to \Upsilon_{T,k}}\mu_k(t)\!=\!\infty$. Hence, the next lemma follows directly by the definition of solutions to HDSs of the form \eqref{HDS0}.

\vspace{0.05cm}
\begin{lemma}[Bounded Flow-Time]\label{lemma101}
Let $z$ be a maximal solution to the HDS \eqref{mainHDSmodel} with $\psi(0,0)\in\Psi\coloneqq \Psi_{C}\cup \Psi_{D}$ and $\mu_k(0,0)\in\mathbb{R}_{\geq1}$. Then, the HTD of $z$ satisfies $\sup_t(\text{dom}(z))\leq \Upsilon_{T,k}$, where $\Upsilon_{T,k}$ is given by \eqref{tpmu}. \QEDB
\end{lemma}

\vspace{0.1cm}
In words, Lemma \ref{lemma101} states that the total amount of flow-time of \emph{every} solution of the HDS \eqref{mainHDSmodel} will be bounded by $\Upsilon_{T,k}$. We will refer to this quantity as the \emph{prescribed time} (PT), and we emphasize its dependency on the initial value of $\mu_0 = \mu_{k}(0)$ as well as the constants $(T,k)$. In the literature on PT-S in continuous-time systems, $\mu_0$ is usually taken to be equal to one. However, for the sake of generality, in this paper all our results are expressed in terms of $\mu_0\in \mathbb{R}_{\ge 1}$.

A useful property of the BU-ODE is that when normalized by $\mu_k$ its unique solutions (from initial conditions in $\mathbb{R}_{\geq1}$) are always complete and uniformly lower bounded, a property that can also be established directly by integration. 

\vspace{0.1cm}
\begin{lemma}\label{completemulemma}
Let $k\geq1$, and consider the normalized BU-ODE $\frac{d\hat{\mu}_k}{ds}=\frac{k}{T}\hat{\mu}_k(s)^{\frac{1}{k}}$ with $\hat{\mu}_k(0)=\mu_0\in\mathbb{R}_{\geq1}$, evolving in the $s$-time scale. Then, its unique solution satisfies: (a) For $k=1$: $\hat{\mu}_k(s)=\mu_0 e^{\frac{s}{T}}\geq 1$ for all $s\geq0$; (b) For $k>1$: $\hat{\mu}_k(s)=\left(\frac{(k-1)}{T}s+\mu_0^{\frac{k-1}{k}}\right)^{\frac{k}{k-1}}\geq 1$,
for all $s\geq0$. \QEDB
\end{lemma}

\vspace{-0.2cm}
\subsection{Time-Scaling of Hybrid Time Domains}
The signals $\mu_k$ generated by the dynamics \eqref{flowshybrid1a} can be used to define a suitable dilation and contraction on the HTD of the solutions to the HDS \eqref{HDS0}. Specifically, for each $(T,k)\in\mathbb{R}_{>0}\times\mathbb{R}_{\geq1}$, and $1\leq a\leq b$, consider the function $\omega_k:\mathbb{R}_{\geq1}\times\mathbb{R}_{\geq1}\to\mathbb{R}_{\geq0}$ defined as
\begin{equation}\label{omegaformula}
\omega_{k}(b,a):=\frac{T}{k}\left(\frac{b^{\rho(k)}-1}{\rho(k)}-\frac{a^{\rho(k)}-1}{\rho(k)}\right),\quad \forall~k>1,
\end{equation}
where $\rho(k)\coloneqq\frac{k-1}{k}$, and  $\omega_{1}(b,a)\coloneqq\lim_{k\to1^+}\omega_{k}(b,a)$ when $k=1$. The following proposition states some useful properties of $\omega_{k}(\cdot,\cdot)$ when evaluated along the trajectories of $\mu_k$. The proof employs standard calculus theorems, and it is presented in the Appendix. 
\begin{prop}\label{transformationk}
Let $(T,k)\in\mathbb{R}_{>0}\times\mathbb{R}_{\geq1}$, $\mu_k$ be given by \eqref{tpmu}, and let $\mathcal{T}_{k}:[0,\Upsilon_{T,k})\to\mathbb{R}_{\geq0}$ be the function
\begin{equation}\label{inverseformula2}
\mathcal{T}_{k}(t):=\omega_{k}(\mu_k(t),\mu_k(0)),~~\forall~t\in[0,\Upsilon_{T,k}).
\end{equation}
Then, $\mathcal{T}_{k}(\cdot)$ satisfies the following properties:
\vspace{0.1cm}
\begin{enumerate}[(P1)]
\item $\lim_{t\to \Upsilon_{T,k}}\mathcal{T}_{k}(t)=\infty.$
\item For any pair $t_2,t_1\in[0,\Upsilon_{T,k})$ such that $t_2\geq t_1$:
\begin{equation*}
\mathcal{T}_{k}(t_2)-\mathcal{T}_{k}(t_1)=\omega_{k}(\mu_k(t_2),\mu_k(t_1)).
\end{equation*}
\item For all $t\in[0,\Upsilon_{T,k})$, $\mathcal{T}_{k}$ satisfies  $$\frac{d\mathcal{T}_{k}(t)}{dt}=\mu_k(t),~~~\mathcal{T}_{k}(0)=0.$$  
\item For all $t\in[0,\Upsilon_{T,k})$, $\mathcal{T}_{k}$ has a well-defined inverse $\mathcal{T}^{-1}_{k}:\mathbb{R}_{\geq0}\to\mathbb{R}_{\geq0}$, which is given by
\begin{align}\label{inversetpk}
\hspace{-\leftmargin}
\mathcal{T}^{-1}_{k}(s)&{=}
    \Upsilon_{T,k}\left(1{-}\!\left(1{+}\frac{(k-1)s}{\Upsilon_{T,k}\mu_0}\right)^{\frac{1}{1-k}}\right),~~k>1,
\end{align}
and by $\mathcal{T}_1^{-1}(s)=\lim_{k\to1^+}\mathcal{T}^{-1}_k(s)$ for $k=1$.
\item For all $s\in\mathbb{R}_{\geq0}$, the inverse of $\mathcal{T}_{k}$ satisfies:
\begin{equation}\label{ODEtau01}
 \frac{d}{ds}\mathcal{T}^{-1}_{k}(s)=\frac{1}{\mu_k\left(\mathcal{T}^{-1}_{k}(s)\right)},~~\mathcal{T}^{-1}_{k}(0)=0.
\end{equation}
\item $\lim_{T\to\infty}\mathcal{T}_{k}(t)=\mu_0^{\frac{k-1}{k^2}}t$ for $k>1$, and $\lim_{T\to\infty}\mathcal{T}_{1}(t)=\mu_0t$ for all $t\geq0$. \QEDB
\end{enumerate}

\end{prop}
%
\begin{figure*}[t!]
    \centering
    \includegraphics[width=0.99\linewidth]{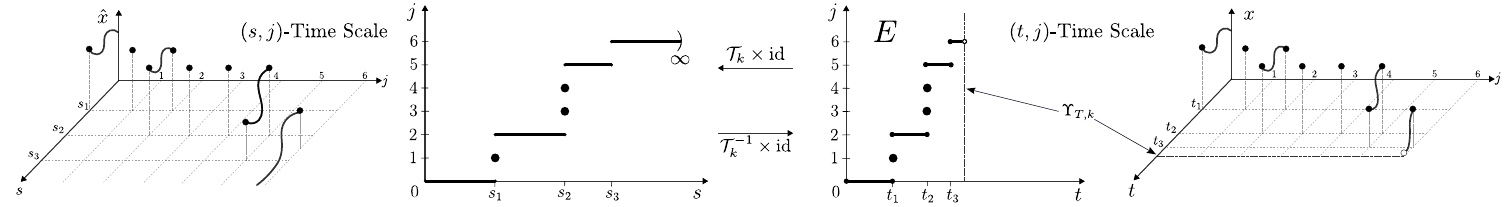}
    \caption{Dilation and contraction of hybrid time domains and hybrid arcs. The structure of the hybrid time domain $E$ in the $(t,j)$-time scale is preserved under the diffeomorphism  $\mathcal{T}_k\times \text{id}$ in the $(s,j)$-time scale.}
        \label{fig:dilation}\vspace{-0.2cm}
\end{figure*}

\vspace{0.25cm}
\begin{remark}\label{casek1}
To contextualize Proposition \ref{transformationk}, consider the special case $k=1$, which is commonly used in the literature on PT-control of ODEs \cite{song2017time,Orlov2022}. In this case, Proposition \ref{transformationk} yields the following ``standard'' mappings:
\begin{subequations}\label{mappingscase1}
\begin{align}
\mathcal{T}_1^{-1}(s)&=\Upsilon_{T,1}\left(1-e^{-\frac{1}{T}s}\right),~~\forall~s\in\mathbb{R}_{\geq0},\label{lnkinv1}\\
\mathcal{T}_1(t)&=T\left(\ln\left(\frac{\mu_1(t)}{\mu_1(0)}\right)\right),~~\forall~t\in[0,\Upsilon_{T,1}).\label{lnk1}
\end{align}
\end{subequations}
This follows because \eqref{inversetpk} can be written as:
$\mathcal{T}^{-1}_k(s)=\Upsilon_{T,1}\left(1-\left(1+\frac{s}{n(k)T}\right)^{-n(k)}\right),$
where $n(k)=\frac{1}{k-1}$. Using  $e^{\frac{s}{T}}=\lim_{n\to\infty}\left(1+\frac{s}{nT}\right)^n$ and the fact that $n\to\infty$ when $k\to1^+$, we obtain \eqref{lnkinv1}. Similarly, using $\rho(k)=\frac{k-1}{k}$, 
$\lim_{\rho\to0}\frac{\mu_1^\rho-1}{\rho}=\ln(\mu_1)$, and $\rho\to0$ if and only if $k\to1$, equation \eqref{lnk1} follows directly from \eqref{omegaformula} by applying the product law for limits.\QEDB
\end{remark}

\vspace{0.1cm}
The properties established in Proposition \ref{transformationk} are used to derive the following result, which indicates a suitable dilation/contraction of the HTDs of system \eqref{mainHDSmodel} when analyzed in a different time scale induced by the diffeomorphisms $\mathcal{T}_{k}$, $k\in\mathbb{Z}_{\geq1}$, see Figure \ref{fig:dilation}.

\vspace{0.1cm}
\begin{prop}[Dilation and Contraction of HTD]\label{dilationProposition}
Let $(T,\mu_0,k)\in\mathbb{R}_{>0}\times\mathbb{R}_{\geq1}\times\mathbb{R}_{\geq1}$, and $\mathcal{T}_{k}$ be given by \eqref{inverseformula2}.
Consider the HDS with state $\hat{z}=(\hat{\psi},\hat{\mu}_k)$ and input $\hat{u}$:
\begin{subequations}\label{rescaledHDS}
\begin{align}
(\hat{z},\hat{u})\in \tilde{C}=C\times\mathbb{R}^m&,~~\dot{\hat{z}}_s\in \frac{1}{\hat{\mu}_k}F_T(\hat{z},\hat{u}).\label{flowrescaled}\\
(\hat{z},\hat{u})\in \tilde{D}=D\times\mathbb{R}^m&,~~\hat{z}^+\in G(\hat{z}).\label{jumprescaled}
\end{align}
\end{subequations}
where the data $\mathcal{H}_r=\{\tilde{C},F_T,\tilde{D},G\}$ in \eqref{rescaledHDS} is the same of \eqref{mainHDSmodel}, and which evolves on the $(s,j)$ hybrid time scale, where $s= \mathcal{T}_{k}(t)$. Then, the following holds:
\begin{enumerate}[(a)]
\item If $(\hat{z},\hat{u})$ is a maximal solution pair of \eqref{rescaledHDS} from the initial condition $z_0$, then the pair of signals defined as $\left(z(t,j),u(t,j)\right):=\left(\hat{z}(s,j),\hat{u}(s,j)\right)$, for all $(s,j)\in\text{dom}(\hat{z})$, is also a maximal solution pair of \eqref{mainHDSmodel} from the initial condition $z_0$ via the time dilation $s=\mathcal{T}_{k}(t)$.
\item  If $(z,u)$ is a maximal solution pair of \eqref{mainHDSmodel} from the initial condition $z_0$, then the pair of signals defined as $\left(\hat{z}(s,j),\hat{u}(s,j)\right):=\left(z(t,j),u(t,j)\right)$ for all $(t,j)\in\text{dom}(z)$, is also a maximal solution pair of \eqref{rescaledHDS} from the initial condition $z_0$ via the time contraction $t=\mathcal{T}_{k}^{-1}(s)$.
\end{enumerate}
\end{prop}
\vspace{0.1cm}
\textbf{Proof:} Let $(T,\mu_0,k)\in\mathbb{R}_{>0}\times\mathbb{R}_{\geq1}\times\mathbb{R}_{\geq1}$, be given. 
We proceed to prove both items of the lemma:

\vspace{0.1cm}\noindent 
(a) Let $(\hat{z},\hat{u})$ be a maximal solution pair of \eqref{rescaledHDS} from $z_0$. Then, for each $j\in\mathbb{Z}_{\geq0}$ such that the interior of $\hat{I}_j:=\{s\geq0:(s,j)\in \text{dom}(\hat{z})\}$ is nonempty, $\hat{z}$ satisfies:
\begin{equation}\label{equaproof123}
\frac{\text{d}}{\text{d}s}\hat{z}(s,j)\in \frac{1}{\hat{\mu}_k(s,j)}F_T(\hat{z}(s,j),\hat{u}(s,j)),
\end{equation}
for almost all $s\in I_j$. Using the chain rule, $z$ satisfies:
\begin{equation*}
\frac{\text{d}}{\text{d}t}z(t,j)=\frac{\text{d}}{\text{d}t}\hat{z}(\mathcal{T}_k(t),j)=\frac{\text{d}}{\text{d}s}\hat{z}(s,j)\cdot \dot{\mathcal{T}}_k(t),
\end{equation*}
and since
$\dot{\mathcal{T}}_{k}(t)
=\mu_k(t)$ for all $t\in[0,\Upsilon_{T,k})$, and $\mu_k$ does not change during the jumps \eqref{jumpmapPT}, using \eqref{equaproof123} we obtain:
\begin{equation*}
\frac{\text{d}}{\text{d}t}z(t,j)=\mu_k(t,j)\frac{\text{d}}{\text{d}s}\hat{z}(s,j)\in \frac{\mu_k(t,j)}{\hat{\mu}_k(s,j)}F_T(\hat{z}(s,j),\hat{u}(s,j)).
\end{equation*}
By construction, $\mu_k(t,j)=\hat{\mu}_k(s,j)$, $u(t,j)=\hat{u}(s,j)$ and $z(t,j)=\hat{z}(s,j)$ via the time dilation $s=\mathcal{T}_{k}(t)$. Therefore, substituting in the above inclusion, and using the continuity of $\mathcal{T}_{k}$, we obtain that $\dot{z}(t,j)$ satisfies \eqref{flowshybrid1a} for almost all $t\in I_j:=\{t\geq0:(t,j)\in\text{dom}(z)\}$. Moreover, note that $\mathcal{T}_{k}(\underline{t}_j)=\underline{s}_j$ and $\mathcal{T}_{k}(\overline{t}_j)=\overline{s}_j$ where $\underline{t}_j:=\min I_j$, $\overline{t}_j=\sup I_j$, $\underline{s}_j:=\min \hat{I}_j$, $\overline{s}_j=\sup \hat{I}_j$. Similarly, for every $(s,j)\in\text{dom}(\hat{z})$ such that $(s,j+1)\in\text{dom}(\hat{z})$, we have that $\hat{z}(s,j+1)\in G(\hat{z}(s,j))$, and therefore $z(t,j+1)\in G(z(t,j))$. Thus $(z,u)$ is a maximal solution to \eqref{mainHDSmodel}.

\vspace{0.1cm}\noindent 
(b) Let $(z,u)$ be a maximal solution pair of \eqref{mainHDSmodel} from $z_0$. Using again the chain rule, and the definition of $\hat{z}$, we obtain that for each $j$ such that the interior of $I_j:=\{t\geq0:(t,j)\in\text{dom}(z)\}$ is nonempty, the signal $\hat{z}$ satisfies:
\begin{equation*}
\frac{\text{d}}{\text{d} s}\hat{z}(s,j)=\frac{\text{d} z}{\text{d} \mathcal{T}_k^{-1}}\frac{\text{d} \mathcal{T}_k^{-1}}{\text{d} s}= \frac{\dot{z}(t,j)}{\mu_k(t,j)}\in \frac{F_T(z(t,j),u(t,j))}{\mu_k(t,j)},
\end{equation*}
where we used \eqref{flowshybrid1a} and \eqref{ODEtau01}. Since by construction $\hat{z}(s,j)=z(t,j)$, $\hat{\mu}_k(s,j)=\mu_k(s,j)$, and $\hat{u}(s,j)=u(s,j)$ via the time contraction $t=\mathcal{T}_k^{-1}(s)$, substituting in the above expression we obtain that $\hat{z}$ satisfies $\dot{\hat{z}}\in \frac{1}{\hat{\mu}_k}F_T(\hat{z},\hat{u}_k)$ for almost all $s\in \hat{I}_j=\{s\geq0:(s,j)\in\text{dom}(\hat{z})\}$. Moreover, note that $\mathcal{T}_k^{-1}(\underline{s}_j)=\underline{t}_j$ and $\mathcal{T}_k^{-1}(\overline{s}_j)=\overline{t}_j$ where $\underline{t}_j:=\min I_j$, $\overline{t}_j=\sup I_j$, $\underline{s}_j:=\min \hat{I}_j$, $\overline{s}_j=\sup \hat{I}_j$. Since for every $(t,j)\in\text{dom}(z)$ such that $(t,j+1)\in\text{dom}(z)$, we have that $z(t,j+1)\in G(z(t,j))$, and therefore $\hat{z}(s,j+1)\in G(\hat{z}(s,j))$, it follows that $(\hat{z},\hat{u})$ is a maximal solution pair to \eqref{rescaledHDS}.  \hfill $\blacksquare$

\vspace{0.1cm}
\begin{remark}
Proposition \ref{dilationProposition} establishes a relationship between the solutions of the HDS \eqref{mainHDSmodel} in the $(t,j)$ time scale, and the solutions of \eqref{rescaledHDS} in the $(s,j)$ time scale via the family of $k$-parameterized dilations $s=\mathcal{T}_k(t)$ and contractions $\mathcal{T}_k^{-1}(s)$. In particular, the function $\mathcal{T}_k:[0,\Upsilon_{T,k})\to\mathbb{R}_{\geq0}$ will define a \emph{diffeomorphism} that preserves the structure of the HTD of the hybrid arcs of \eqref{rescaledHDS}. This observation is at the core of our analysis, as it enables us to carry out the analysis of the HDS \eqref{mainHDSmodel} based on the qualitative behavior of the solutions of system \eqref{rescaledHDS}, which has a flow map normalized by $\hat{\mu}_k$, leading to a system with no finite escape times in $\hat{\mu}_k$ via Proposition \ref{completemulemma}. This normalized HDS can be seen as a ``target'' system that can be studied and designed using the rich set of tools that already exist in the literature of hybrid control \cite{bookHDS,sanfelice2021hybrid}. \QEDB
\end{remark}

\vspace{0.1cm}
\begin{remark}\label{remark:equivalentTransform}
Using the definition of $\mu_k$ in \eqref{tpmu} with $k>1$, equation \eqref{inverseformula2} can also be written as 
\begin{equation}\label{remark:equivalentTransformEquation}
\mathcal{T}_k(t)=\frac{T\mu_0^{\frac{k-1}{k}}}{k-1}\left(\frac{T^{k-1}}{\left(T-t\mu_0^{\frac{1}{k}}\right)^{k-1}}-1\right),
\end{equation}
which is defined for all $t\in[0,\Upsilon_{T,k})$, and %
which recovers the common dilation found in the literature of ODEs when $\mu_0=1$, see \cite{song2017time}. Other types of transformations are presented in \cite{TRAN2020104605} for the study of finite-time control of ODEs. Proposition \ref{dilationProposition} provides an extension of these results to hybrid systems.
\QEDB
\end{remark}

\vspace{0.1cm}
\begin{remark}
Analyses of HDS based on the time scaling of the continuous-time dynamics are common in the study of singular perturbations \cite{averaging_singularHDS,WangTeelNesic} and averaging theory \cite{TeelNesicAveraging,poveda2017framework}. However, in contrast to \eqref{remark:equivalentTransformEquation}, the time scaling in those scenarios is usually \emph{linear}, i.e., from $\mathbb{R}_{\geq0}$ to $\mathbb{R}_{\geq0}$. \QEDB
\end{remark}
\subsection{PT-S via Flows in HDS}
Since the solutions of the HDS \eqref{mainHDSmodel} can only flow for at most $\Upsilon_{T,k}$ total amount of time, in this paper we are interested in steering the state $z$ to a desired set $\mathcal{A}$ as $t\to \Upsilon_{T,k}$, or before $\Upsilon_{T,k}$, where
\begin{equation}\label{setstable}
\mathcal{A}=\mathcal{A}_1\times\mathbb{R}_{\geq1}.
\end{equation}
For systems with inputs, the following definition aims to capture this property, which makes use of the transformation $\mathcal{T}_k$ defined in \eqref{inverseformula2}, and which extends \cite[Defs. 1]{song2017time} from ODEs to HDS. %

\vspace{0.1cm}
\begin{definition}\label{definitionptflows}
Let $\mathcal{A}$ be given by \eqref{setstable}, where $\mathcal{A}_1\subset\mathbb{R}^{n}$ is compact. Then, for the HDS \eqref{HDS0} with data given by \eqref{mainHDSmodel}, the set $\mathcal{A}$ is said to be \emph{Prescribed-Time Input-to-State Stable via Flows} (PT-ISS$_\text{F}$) if there exists $\beta\in\mathcal{KLL}$ and $\gamma\in\mathcal{K}$ such that for every $z(0,0)\in\mathbb{R}^{n+1}$, all solutions $z$ satisfy the bound
\begin{equation}\label{eq:PTflows1}
|z(t,j)|_{\mathcal{A}} \leq \beta\big(|z(0,0)|_{\mathcal{A}},\mathcal{T}_k(t),j \big)+\gamma\left(\vert u \vert_{(t,j)}\right),
\end{equation}
for all $(t,j)\in\text{dom}(z)$. If \eqref{eq:PTflows1} holds with $u\equiv0$, the set $\mathcal{A}$ is said to be \emph{Prescribed-Time Stable via Flows} (PT-S$_{\text{F}}$). \QEDB
\end{definition}

\vspace{0.1cm}
In some cases, it might be possible to completely suppress the residual effect of the input $u$ in the bound \eqref{eq:PTflows1} via PT feedback. This property, termed \emph{PT-ISS with Convergence} in \cite[Defs. 1]{song2017time}, can also be obtained in hybrid systems:

\vspace{0.1cm}
\begin{definition}\label{definitionptflowsconv}
Let $\mathcal{A}$ be given by \eqref{setstable}, where $\mathcal{A}_1\subset\mathbb{R}^{n}$ is compact. Then, for the HDS \eqref{HDS0} with data given by \eqref{mainHDSmodel}, the set $\mathcal{A}$ is said to be \emph{Prescribed-Time Input-to-State Stable and Convergent via Flows} (PT-ISS-C$_\text{F}$) if there exists $\beta\in\mathcal{K}\mathcal{L}\mathcal{L}$, $\gamma\in\mathcal{K}$, and $\beta_c\in\mathcal{K}\mathcal{L}$ such that for every $z(0,0)\in\mathbb{R}^{n+1}$, all solutions $z$ satisfy the bound
\begin{equation}\label{eq:PTflows2}
|z(t,j)|_{\mathcal{A}} \leq \beta_c\left(\beta\left(|z(0,0)|_{\mathcal{A}}, \mathcal{T}_k(t),j \right) +  \gamma\left(|u|_{(t,j)}\right),\mathcal{T}_k(t)\right),
\end{equation}
for all $(t,j)\in\text{dom}(z)$. \QEDB
\end{definition}

\vspace{0.1cm}
\begin{remark}
Definitions \ref{definitionptflows} and \ref{definitionptflowsconv} extend existing PT-S notions in the literature of ODEs \cite[Defs. 1]{song2017time} to hybrid systems. The use of $\mathcal{K}\mathcal{L}\mathcal{L}$ functions is common in the analysis of hybrid systems, enabling us to differentiate convergence behaviors in the continuous-time domain from those in the discrete-time domain. Additionally, since $|\mu_{k}(t)|_{\mathbb{R}_{\ge 1}}=0$ for all $t\in[0, \Upsilon_{T,k})$, we can equivalently express the bounds \eqref{eq:PTflows1} and \eqref{eq:PTflows2} with $z$ replaced by $\psi$ and $\mathcal{A}$ replaced by $\mathcal{A}_1$.\QEDB
\end{remark}

\vspace{0.1cm}
\begin{remark}[On the Uniformity with Respect to $\mu_0$]
By using equations \eqref{omegaformula}-\eqref{inverseformula2}, the bounds \eqref{eq:PTflows1} and \eqref{eq:PTflows2} can also be written in terms of $\omega(\mu_k,\mu_0)$ to stress their dependence on $\mu_k(0)$, thus providing some ``uniformity'' in the bounds. The example below, which follows as a particular case of the main results presented in the next section, illustrates this observation. \QEDB 
\end{remark}

\vspace{0.1cm}
\begin{example}\label{example1positive}
Consider the HDS \eqref{mainHDSmodel} with $k=1$, $T=1$, $\psi=(x,\tau)$, $F_{\Psi}=\{-x+u\}\times\{1\}$,~ $G_{\Psi}=\{\frac{1}{2}x\}\times\{0\}$, $\Psi_C=\mathbb{R}^n\times[0,1]$, and $\Psi_D=\mathbb{R}^n\times\{1\}$. Then, every solution $z=(x,\tau,\mu_1)$ satisfies the bound (see proof of Theorem \ref{theorem1}):
\begin{equation*}
|z(t,j)|_{\mathcal{A}}\leq k_1 e^{-k_2\mathcal{T}_1(t)}\left(e^{-k_3(\mathcal{T}_1(t)+j)}|z(0,0)|_{\mathcal{A}}+k_4|u|_{(t,j)}\right),
\end{equation*}
where $k_i>0$ and $\mathcal{A}_1=\{0\}\times[0,1]$, for all $(t,j)\in\text{dom}(z)=\text{dom}(u)$ and any continuous and bounded input $u$. Moreover, since $|z|_{\mathcal{A}}=|x|$ for all $z=(\psi,\mu_1)\in C\cup D$, and using \eqref{lnk1}, the above bound can be written as:
\begin{equation*}
\!\!\!\!\left|x(t,j)\right|\leq \frac{\mu_1(0,0)^{\alpha_1}}{\mu_1(t,j)^{\alpha_2}}\!\left(\frac{e^{-qj}}{\mu_1(t,j)^{\alpha_3}}|x(0,0)|+\alpha_4\cdot|u|_{(t,j)}\right)
\end{equation*}
where $\alpha_i>0$, $\mu_1(0,0)=\mu_0\ge 1$, and for all $(t,j)\in\text{dom}(z)$. It follows that $\lim_{(t,j)\in\text{dom}(x),t\to \Upsilon_{1,1}}x(t,j)=0$.
\QEDB 
\end{example}

\vspace{0.15cm}
Unlike ODEs, in HDS the existence of bounds of the form \eqref{eq:PTflows1}-\eqref{eq:PTflows2} do not necessarily guarantee that the solutions $z$ will converge to the set $\mathcal{A}$ as $t\to \Upsilon_{T,k}$, for any $\Upsilon_{T,k}>0$, even if $u\equiv0$ and $z$ is complete. The following scalar example illustrates this scenario.
\vspace{0.1cm}

\begin{example}
Consider the HDS \eqref{mainHDSmodel} with $F_{\Psi}=\{-\psi\}$,~ $G_{\Psi}=\frac{1}{2}\psi$, $\Psi_C=(-\infty,-1]\cup[1,\infty)$, and $\Psi_D=[-1,1]$. For this system, we study stability with respect to the set $\mathcal{A}=\mathcal{A}_1\times\mathbb{R}_{\geq1}$, with $\mathcal{A}_1=\{0\}$. For any initial condition $z(0,0)=(\psi(0,0),\mu(0,0))$ satisfying $|\psi(0,0)|>1$ and $\mu_0=1$, the unique solutions to the HDS satisfy $x(t,0)=x(0,0)\left(\frac{T-t}{T}\right)^T$, for all $(t,j)\in[0,t']\times\{0\}$, where $t'=T(1-\psi(0,0)^{-\frac{1}{T}})$, and $x(t,j)=\left(\frac{1}{2}\right)^jx(t',0)$, for all $(t,j)\in \bigcup_{j\in\mathbb{Z}_{\geq1}}\{t'\}\times\{j\}$. It follows that $x(t,j)\to \mathcal{A}_1$ only as $j\to\infty$. Yet, every solution of the HDS satisfies \eqref{eq:PTflows1} with $u=0$. This follows by a direct application of item (a) of Proposition \ref{dilationProposition}, the result of \cite[Thm. 1]{LyapunovExponentialHDS}, and item (b) of Proposition \ref{dilationProposition}, in that order. \QEDB
\end{example}

\vspace{0.1cm}
The previous example shows that bounds of the form \eqref{eq:PTflows1} or \eqref{eq:PTflows2} only guarantee PT-S-like behaviors via the \emph{flows} of the HDS. Therefore, to emulate the existing PT-S bounds obtained for ODEs \cite{song2017time,Orlov2022}, the HDS \eqref{mainHDSmodel} with $\mu_k \equiv 1$ must generate maximal solutions with hybrid time domains $E$ satisfying $\text{sup}_t E=\infty$, such as those in Example \ref{example1positive}. In general, this is not possible whenever $C=\emptyset$, or whenever system \eqref{mainHDSmodel} has eventually discrete, or purely discrete solutions. However, as shown in the next section, for R-Switching systems, discrete solutions can be ruled out by the design of an appropriate switching signal that also exploits the ``blow-up'' nature of the functions $\mu_k$.  
\section{PT-ISS IN R-SWITCHING SYSTEMS}\label{sectionISSTheorems}
\subsection{Modeling Framework and Main Assumptions}
We consider time-varying R-Switching systems of the form \eqref{vectorfield1}-\eqref{resetmap1}, where:
\begin{align}\label{actualflow1}
\tilde{f}_{\sigma(t)}(x,u,t)=\mu_k(t)\cdot f_{\sigma(t)}(x,u,\tau(t)),~~~t\notin \mathcal{W}(\sigma),
\end{align}
characterizes the continuous-time evolution of the state $x$. In \eqref{actualflow1}, $\tau:\mathbb{R}_{\geq0}\to\mathbb{R}$ is a locally absolutely continuous function that is uniformly bounded, and that is generated by the following system:
\begin{subequations}\label{blowupautomaton}
\begin{align}
\dot{\tau}&\in\mu_k(t)\cdot \left[0,\frac{1}{\tau_d}\right],~~~t\notin \mathcal{W}(\sigma),\label{setvaluedflows}\\
\tau^+&=\tau-1,~~~~~~~~~~~~~t\in \mathcal{W}(\sigma),
\end{align}
\end{subequations}
where $\mu_k$ is given by \eqref{tpmu} and $\tau_d>0$. To contextualize this model, some remarks are in order.
\begin{remark}
When $\mu_k\equiv 1$ and $f_{\sigma}$ does not depend on $\tau$ and $u$, the model \eqref{actualflow1} coincides with the conventional nonlinear switching systems examined in \cite{AverageDwellTime,Average_Dwell_time}. On the other hand, when $f_{\sigma}$ depends on $u$, \eqref{actualflow1} captures nonlinear switching systems with inputs, similar to those studied  \cite{GuosongLiberzon,Liu_Tanwani_Liberzon2022}. 
\QEDB
\end{remark}
\begin{remark}
When $\mu_k\equiv 1$ and $f_{\sigma}$ depends on $\tau$, system \eqref{actualflow1} describes a class of $\tau$-parameterized nonlinear switching systems. In this class, $\tau$ is not necessarily constant throughout time, and the function $t\mapsto \tau(t)$ may not be differentiable or even continuous. Such models emerge in a class of time-triggered reset systems \cite{zero_order_poveda_Lina,ochoa2021momentum} suitable for optimization and learning problems see also Section \ref{resetsexample}. %
\QEDB
\end{remark}
\begin{remark}
In many applications, the original system of interest might not have the exact form of \eqref{actualflow1}. This holds true, for instance, in the context of PT-regulation of affine dynamical systems with non-zero drift, where multiplying the entire vector fields by the gain $\mu_k$ is not feasible. Nevertheless, as shown later in Section \ref{sectionapplications}, through an appropriate feedback design or variable transformation, the systems of interest can be reformulated into the form \eqref{actualflow1}. This can be achieved either directly in the dynamics of the system or in the analysis of the system using Lyapunov-based techniques.   \QEDB 
\end{remark}

\vspace{0.1cm}
\begin{assumption}\label{assumptionLipschitz}
For each $q\in \mathcal{Q}$, the function $f_q:\mathbb{R}^n\times\mathbb{R}^m\times\mathbb{R}_{\geq0}\to\mathbb{R}^n$ is locally Lipschitz, $R_q:\mathbb{R}^n\to\mathbb{R}^n$ is continuous, and the input $u:\mathbb{R}_{\geq0}\to\mathbb{R}^l$ is continuous and bounded.  \hfill \QEDB
\end{assumption}
\begin{remark}
It is possible to relax the Lipschitz and continuity assumptions on $f_q$ and $R_q$ by considering set-valued mappings with suitable regularity properties (outer-semicontinuity, local boundedness, etc), see \cite{Liu_Tanwani_Liberzon2022}. However, since the main focus of this paper is on obtaining Prescribed-Time stability properties, we work with Lipschitz and continuous functions to simplify our presentation. \QEDB
\end{remark}

\vspace{0.1cm}
Next, we introduce our main stability assumptions on the ``target'' HDS \eqref{rescaledHDS}. In this assumption, we use the function $\Delta(\mu_k)$ to characterize the effect of $\mu_k$ on the input $u$ in \eqref{actualflow1}. Specifically, we are interested in the cases $\Delta(\mu_k)=0$, $\Delta(\mu_k)=1$, and $\Delta(\mu_k)=\mu_k^{-\ell}$ with $\ell>0$. Particular applications of each case will be discussed in Section V.
\vspace{0.1cm}
\begin{assumption}\label{assumptionunstablemodes}
There exist $\tau_d\in\mathbb{R}_{>0}$, $N_0\in\mathbb{R}_{\geq1}$, smooth functions $V_q:\mathbb{R}^n\times\mathbb{R}_{\geq0}\to\mathbb{R}_{\geq0}$, where $q\in\mathcal{Q}$, and constants $c_{q,i}>0$,~$i\in\{1,2,3,4,5\}$, $p>0$, such that:
\begin{enumerate}[(a)]
\item For all $(x,\tau,q)\in\mathbb{R}^n\times[0,N_0]\times\mathcal{Q}$: 
\begin{equation}\label{quadraticbounds0}
c_{q,1}|x|^p\leq V_q(x,\tau)\leq c_{q,2}|x|^p.
\end{equation}
\item For all $(x,\tau,q,\mu,\eta)\in\mathbb{R}^n\times[0,N_0]\times\mathcal{Q}_s\times\mathbb{R}_{\geq1}\times[0,\tau_d^{-1}]$: 
\begin{align}\label{stablestabilitycondition}
\!\!\!\!\left\langle \nabla V_q(x,\tau), \left(\begin{array}{c}
f_q(x,u,\tau)\\
\eta
\end{array}\right)
\right\rangle &\leq - c_{q,3} V_q(x,\tau)\notag\\
&+c_{q,4}\Delta(\mu_k)|u|^{p}.
\end{align}
\item  For all $(x,\tau,q,\mu,\eta)\in\mathbb{R}^n\times[0,N_0]\times\mathcal{Q}_u\times\mathbb{R}_{\geq1}\times[0,\tau_d^{-1}]$: 
\begin{align}\label{unstablestabilitycondition}
\!\!\!\!\left\langle \nabla V_q(x,\tau), \left(\begin{array}{c}
f_q(x,u,\tau)\\
\eta
\end{array}\right)
\right\rangle &\leq c_{q,5} V_q(x,\tau)\notag\\
&+c_{q,4}\Delta(\mu_k)|u|^{p}.
\end{align}
\item  For all $(x,\tau)\in\mathbb{R}^n\times[1,N_0]$ and $o,q\in \mathcal{Q}$ such that $q\neq o$:
\begin{align}\label{inequalityresets}
V_q(R_o(x,u),\tau-1)\leq \chi V_o(x,\tau),
\end{align}
where $\chi\in(0,1]$.\QEDB
\end{enumerate}
\end{assumption}

\vspace{0.1cm}
\begin{remark}
Inequalities \eqref{quadraticbounds0}-\eqref{stablestabilitycondition} are common in the context of exponential stability in continuous-time and hybrid systems. For the case when the vector field $f_q$ in \eqref{actualflow1} does not depend on $\tau$, the function $V_q$ can also be taken to be independent of $\tau$. This is the most common situation in switching systems and systems with resets. An example where $f_q$ does depend on $\tau$ will be studied in Section \ref{resetsexample}. \QEDB
\end{remark}
\begin{figure*}[t!]
    \centering
    \includegraphics[width=.48\linewidth]{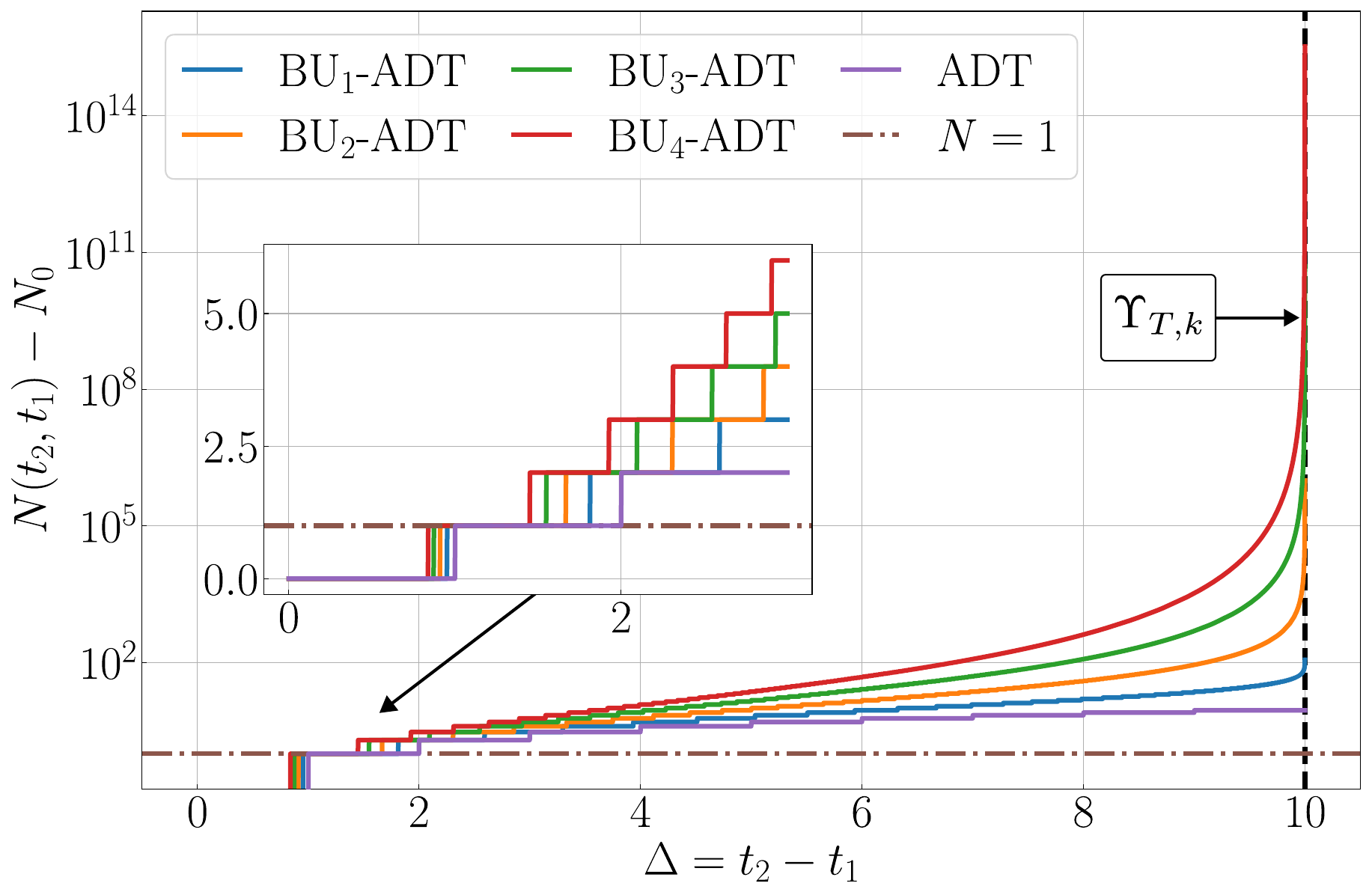}
     \includegraphics[width=.48\linewidth]{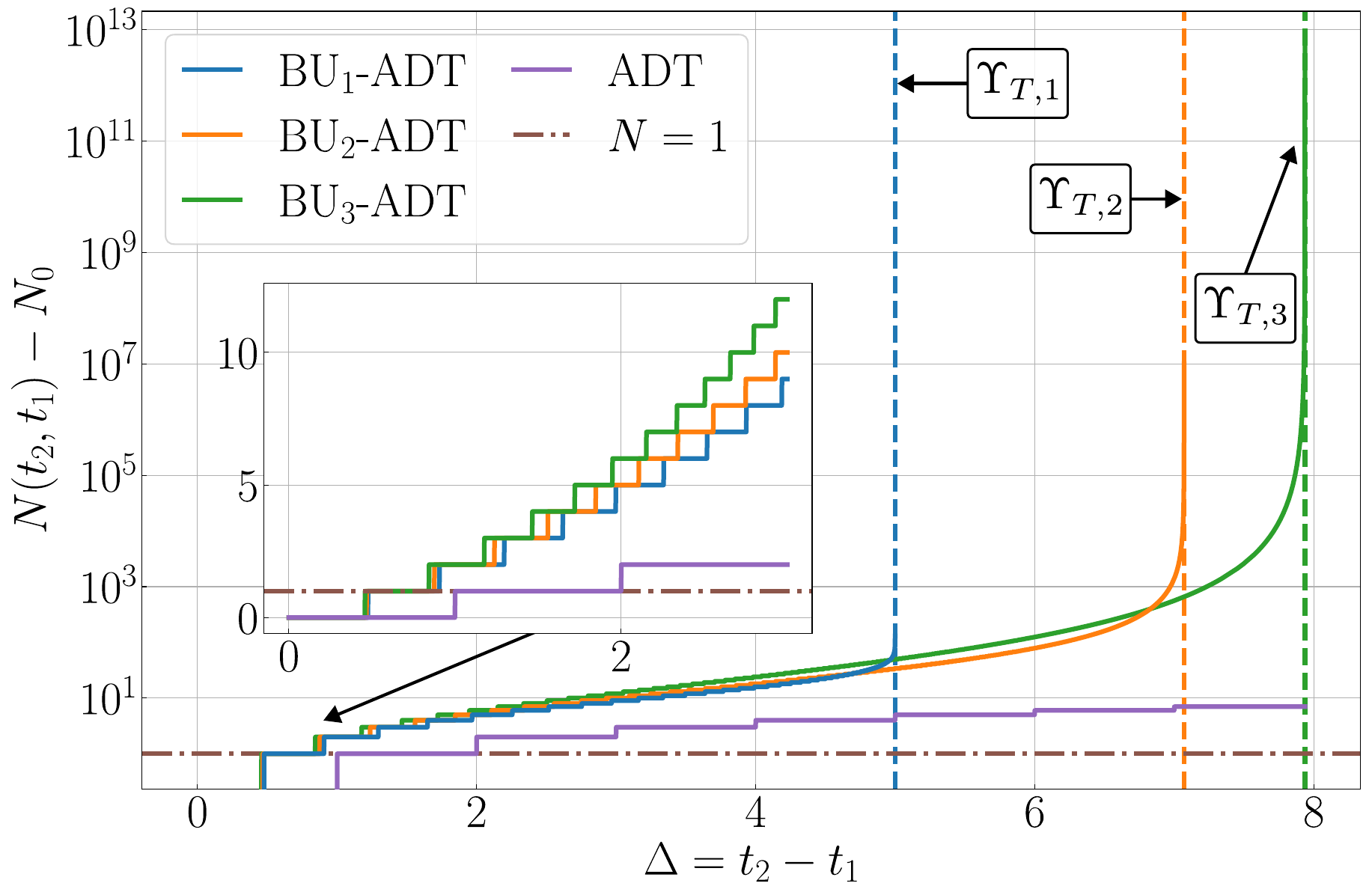}
    \caption{BU$_k$-ADT condition \eqref{BUADTk} for $k\in\{1,2,3,4\}$. Left: When $\mu_0=1$, $T=10$, and $t_1=0$, there exists a single common terminal time $T=\Upsilon_{T,k}$ for all $k$. Right: When $\mu_0=2$, $T=10$, and $t_1=0$, the dependence of $\Upsilon_{T,k}$ on $\mu_0$ (see Proposition \ref{lemmaodek}) leads to the emergence of three distinct terminal times.}\vspace{-.2cm}
    \label{fig:ADTmu01}
\end{figure*}
\vspace{0.1cm}
\begin{remark}
Inequality \eqref{stablestabilitycondition} in item (b) gives a standard decrease condition on the Lyapunov functions $V_q$, for each stable mode $q\in\mathcal{Q}_s$, and up to a neighborhood of the origin, whose size is parameterized by $\Delta(\mu_k)|u|^p$. When $\Delta(\mu_k)=0$, and by \cite[Thm. 1]{LyapunovExponentialHDS}, conditions \eqref{quadraticbounds0}-\eqref{stablestabilitycondition} imply that each mode $q\in\mathcal{Q}_s$ renders the origin exponentially stable in the dilated time scale $s=\mathcal{T}_k(t)$ (see Proposition \ref{dilationProposition}). When $\Delta(\mu_k)=1$, and by \cite[Prop. 1]{Cai2009}, conditions \eqref{quadraticbounds0}-\eqref{stablestabilitycondition} imply that each mode $q\in\mathcal{Q}_s$ renders the origin ISS with exponential decay in the dilated time scale. On the other hand, the case $\Delta(\mu_k)=\mu_k^{-\ell}$, with $\ell>0$, will emerge in the context of PT-regulation where convergence bounds of the form \eqref{eq:PTflows2} are sought-after. An example is presented in Section \ref{sectionapplications}. \QEDB
\end{remark}

\vspace{0.1cm}
\begin{remark}
Inequality \eqref{unstablestabilitycondition} in item (c) is used to rule out finite escape times during the unstable modes $q\in\mathcal{Q}_u$. Similar assumptions are consider in the context of asymptotic and exponential stability in switching systems \cite{GuosongLiberzon,poveda2017framework}. When $\mathcal{Q}_u= \emptyset$ (i.e., there are no unstable modes), item (c) holds vacuously. \QEDB
\end{remark}

\vspace{0.1cm}
\begin{remark}
Inequality \eqref{inequalityresets} in item (d) considers the effect of the resets on the Lyapunov functions related to each of the modes. Usually (e.g., in standard switching systems) $R_q\!=\!\text{id}(\cdot)$ and $V_q$ is independent of $\tau$, and in this case, inequality \eqref{inequalityresets} holds trivially with $\chi\!=\!1$.  When $V_q$ is independent of $\tau$ but $R_q\!\neq\!\text{id}(\cdot)$, item (d) recovers the main assumptions of \cite{Liu_Tanwani_Liberzon2022}. \QEDB
\end{remark}
\subsection{A Class of Blow-Up Average Dwell-Time Conditions}
To achieve asymptotic or exponential stability in switching systems it is common to assume that for all times $t_2\geq t_1\geq0$ the switching signal $\sigma$ satisfies an average dwell-time (ADT) condition of the form:
\begin{equation}\label{ADT00}
N(t_2,t_1)\leq \frac{1}{\tau_d}(t_2-t_1)+N_0,
\end{equation}
where $N(t_1,t_2)$ is the number of switches of $\sigma$ in the interval $(t_1,t_2]$, $\tau_d$ is called the dwell-time, and $N_0$ is the chatter bound, see \cite{AverageDwellTime,Average_Dwell_time} and \cite[Ch. 2.4]{bookHDS}. Since, unlike asymptotic convergence results, PT-S properties are defined with respect to finite intervals of the form $[0,\Upsilon_{T,k})$, we introduce a family of switching signals defined on similar bounded time domains, but which are allowed to switch arbitrarily fast as $t\to \Upsilon_{T,k}$. 
\vspace{0.1cm}
\begin{definition}\label{definitionEADT}
Let $\mu_k$ be given by \eqref{tpmu}. A switching signal $\sigma:[0,\Upsilon_{T,k})\to\mathcal{Q}$ is said to satisfy the \emph{blow-up average dwell-time condition of order $k$} (BU$_k$-ADT) if there exist $N_0\geq1$ and $\tau_d>0$ such that for all $t_2,t_1\in\text{dom}(\sigma)$ the number of switches in the interval $(t_1,t_2]$, denoted $N(t_1,t_2)$, satisfies:
\begin{equation}\label{BUADTk}
N(t_2,t_1)\leq \frac{1}{\tau_d}\omega_k\left(\mu_k(t_2),\mu_k(t_1)\right)+N_0.
\end{equation}
We use $\Sigma_{\text{BU$_k$-ADT}}(\tau_d,N_0,T,\mu_0)$ to denote the family of such signals. \QEDB
\end{definition}

\vspace{0.1cm}
\begin{lemma}\label{equivalentBUADTk}
The BU$_{k}$-ADT condition satisfies:%
\begin{enumerate}[(a)]
\item If $k=1$, then \eqref{BUADTk} is equivalent to:
\begin{equation}\label{BUKADT111}
N(t_2,t_1)\leq \frac{T}{\tau_d}\ln\left(\frac{\Upsilon_{T,1}-t_1}{\Upsilon_{T,1}-t_2}\right)+N_0.
\end{equation}
\item If $k\in\mathbb{Z}_{>1}$, then \eqref{BUADTk} is equivalent to:
\begin{equation*}
N(t_2,t_1)\leq \gamma_k(t_1,t_2)\left[(t_2{-}t_1)\!+\!\sum_{\ell=2}^{k-1}\tilde{c}_{\ell,k}\left(t_2^\ell-t_1^\ell\right)\right]+N_0,
\end{equation*}
where $\tilde{c}_{\ell,k}:=(-1)^{\ell+1}\frac{b_{k,l}}{k-1}\Upsilon_{T,k}^{1-\ell}$, $b_{k,l}=\frac{(k-1)!}{\ell!(k-\ell-1)!}$ and
\begin{align*}
\gamma_k(t_1,t_2):=\frac{\mu_0^{\frac{2-k}{k}}}{\tau_d}\left(\frac{T^2}{\left(\Upsilon_{T,k}-t_2\right)\left(\Upsilon_{T,k}-t_1\right)}\right)^{k-1}.
\end{align*}
\item For all $k\in\mathbb{Z}_{\geq1}$, \eqref{BUADTk} satisfies:
\begin{equation*}
\lim_{T\to\infty}\frac{1}{\tau_d}\omega_k\left(\mu_k(t_2),\mu_k(t_1)\right)+N_0=\frac{\mu_0}{\tau_d}(t_2-t_1)+N_0,
\end{equation*}
thus recovering the ADT condition \eqref{ADT00} when $\mu_0=1$.\QEDB
\end{enumerate}
\end{lemma}

\vspace{0.1cm}
\textbf{Proof:} All properties follow by direct substitution and computations using \eqref{tpmu}. For completeness, the proof is presented in the Supplemental Material. \hfill $\blacksquare$

\vspace{0.1cm}
Figure \ref{fig:ADTmu01} compares the different bounds obtained in \eqref{BUADTk} (in logarithmic scale), as a function of $\Delta=t_2-t_1$, with $t_1=0$, and for different values of $k\in\mathbb{Z}_{\geq1}$, with $\mu_0=1$ (left plot) and $\mu_0=2$ (right plot). The standard ADT bound is also shown, using purple color. Unlike the ADT bound, the BU$_k$-ADT bound grows to infinity as $\Delta\to \Upsilon_{T,k}$, allowing for a Zeno-like behavior where the number of switches grows to infinity as $t\to \Upsilon_{T,k}$. However,  for any compact sub-interval of $[0,\Upsilon_{T,k})$ the allowable number of switches is always bounded.
\subsection{PT-ISS in R-Switching Systems with Stable Modes}
We first consider the case when all the modes $f_q$ in \eqref{actualflow1} are stable, i.e., $\mathcal{Q}_u=\emptyset $ and $\mathcal{Q}=\mathcal{Q}_s$. In this case, the R-Switching system \eqref{vectorfield1}-\eqref{resetmap1} can be analyzed by considering the HDS \eqref{mainHDSmodel}, with $\psi=(x,\tau,q)\in\mathbb{R}^{n+2}$, set-valued mappings:
\begin{subequations}\label{mapsstable01}
\begin{align}
&F_{\Psi}(\psi,u,\mu_k):=\{f_q(x,u,\tau)\}\times
\left[0,\frac{1}{\tau_d}\right]\times\{0\},\\
&G_{\Psi}(\psi,u,\mu_k):=
\{R_q(x)\}\times\{\tau-1\}\times \mathcal{Q}_s\backslash\{q\},
\end{align}
and sets:
\begin{equation}\label{sets01}
\Psi_C=\mathbb{R}^n\times[0,N_0]\times\mathcal{Q}_s,~~\Psi_D=\mathbb{R}^n\times[1,N_0]\times\mathcal{Q}_s.
\end{equation}
\end{subequations}
There is a close connection between the HTDs of the solutions of system \eqref{mainHDSmodel} with data \eqref{mapsstable01}, and the switching signals $\sigma$ that satisfy the BU$_k$-ADT condition \eqref{BUADTk}: %

\vspace{0.1cm}
\begin{lemma}\label{lemmaequivalencedomains}
Let $(F_{\Psi},G_{\Psi},\Psi_C,\Psi_D)$ be given by \eqref{mapsstable01}-\eqref{sets01}, and consider the HDS \eqref{mainHDSmodel} under Assumption \ref{assumptionLipschitz}-\ref{assumptionunstablemodes}. Then, Assumption \ref{wellposedassumption} holds, and:

\begin{enumerate}[(a)]
\item For every maximal solution $z$ and for any pair $(t_1,j_1),(t_2,j_2)\in \text{dom}(z)$, with $t_2>t_1$, inequality \eqref{BUADTk} holds with  
$N(t_2,t_1)=j_2-j_1$.
\item  For every hybrid time domain satisfying property (a), there exists a solution $z$ of \eqref{mainHDSmodel} having the said hybrid time domain.  \hfill\strut \QEDB
\end{enumerate}
\end{lemma}

\vspace{0.1cm}
\textbf{Proof:} Using \eqref{mainHDSmodel}, we obtain that the overall HDS has state $z=(x,\tau,q,\mu_k)\in\mathbb{R}^{n+3}$, %
and the following data:
\begin{subequations}\label{hybridBlowUpAutomatonStable}
\begin{align}\label{flowshybrid1}
&z\in C:=\mathbb{R}^n\times [0, N_0]\times\mathcal{Q}_s\times\mathbb{R}_{\geq1},\quad
\dot{z}=\left(\begin{array}{c}
\vspace{0.3cm}
\dot{x}\\
\vspace{0.4cm}
\dot{\tau}\\
\vspace{0.4cm}
\dot{q}\\ \vspace{0.2cm} 
\dot{\mu}_k
\end{array}\right)\in F_T(z,u):=\left(\begin{array}{c}
\vspace{0.2cm}
\mu_k \cdot f_q(x,u,\tau)\\\vspace{0cm}
\left[0,\dfrac{\mu_k}{\tau_d}\right]\vspace{0.2cm}\\
0 \vspace{0.2cm}\\
 \dfrac{k}{T}\mu_k^{1+\frac{1}{k}}
\end{array}\right),\\
&z\in D:=\mathbb{R}^n\times [1, N_0]\times \mathcal{Q}_s\times\mathbb{R}_{\geq1},\quad
z^+=\left(\begin{array}{c}
\vspace{0.2cm}
x^+\\
\tau^+\vspace{0.1cm}\\
q^+ \vspace{0.1cm}\\
 \mu_k^+
\end{array}\right) \in G(z):=\left(
\begin{array}{c}
\vspace{0.2cm}
R_q(x)\\\vspace{0.1cm}
\tau-1\vspace{0.1cm}\\
\mathcal{Q}_s\backslash\{q\}\vspace{0.2cm}\\
\mu_k
\end{array}\right).
\end{align}
\end{subequations}
Since the function $\mu_k$ generated by \eqref{hybridBlowUpAutomatonStable} is precisely \eqref{tpmu}, which has a finite escape time at $t=\Upsilon_{T,k}$, any solution $z:\text{dom}(z)\to\mathbb{R}^{n+3}$ to \eqref{hybridBlowUpAutomatonStable} will necessarily satisfy $\text{length}_t(\text{dom}(z))\leq \Upsilon_{T,k}$. By Proposition \ref{dilationProposition}, the corresponding HDS \eqref{rescaledHDS} in the $(s,j)$-time scale is given by:
\begin{subequations}\label{hybridautomaton}
\begin{align}\label{flowshybrid11}
&\hat{z}\in C,~ \dot{\hat{z}}_s =\left(\begin{array}{c}
\vspace{0.2cm}
\dot{\hat{x}}_s\\\vspace{0.1cm}
\dot{\hat{\tau}}_s\vspace{0.2cm}\\
\dot{\hat{q}}_s\vspace{0.2cm}\\
\!\!\dot{\hat{\mu}}_{k_s}\!\!
\end{array}\right)\in \hat{F}_T(\hat{z},\hat{u}):= \left(\begin{array}{c}
\vspace{0.2cm}
\!\!f_{\hat{q}}(\hat{x},\hat{u},\hat{\tau})\!\!\\\vspace{0.1cm}
\left[0,\dfrac{1}{\tau_d}\right]\vspace{0.1cm}\\
0 \vspace{0.1cm}\\
 \dfrac{k}{T}\hat{\mu}_k^{\frac{1}{k}}
\end{array}\right)\\
&\hat{z}\in D,~~~~~~~~~\hat{z}^+\in G(\hat{z}),
\end{align}
\end{subequations}
where $C$, $D$, and $G$ were defined in \eqref{hybridBlowUpAutomatonStable}. %
 Since the dynamics of ($\hat{x},\hat{\tau},\hat{q}$) are decoupled from $\hat{\mu}_k$, and since $\hat{\mu}_k$ remains constant during jumps, we directly obtain $\hat{\mu}_k(s,j)$ for any $(s,j)\in\text{dom}(\hat{z})$ using Proposition \ref{completemulemma}:
\begin{equation*}%
\hat{\mu}_k(s,j)=\left(\frac{(k-1)}{T}t+\hat{\mu}(\underline{s}_j,j)^{\frac{k-1}{k}}\right)^{\frac{k}{k-1}},
\end{equation*}
for $k>1$, and $\hat{\mu}_k(s,j)=\hat{\mu}_k(\underline{s}_j,j)e^{\frac{s}{T}},$
for $k=1$, where $\underline{s}_j :=\min\{s\geq0:(s,j)\in\text{dom}(\hat{z})\}$.  By Assumption \ref{assumptionunstablemodes}, the state $\hat{x}$ in \eqref{hybridautomaton} has no finite escape times, and by \cite[Ex. 2.15]{bookHDS} every solution $\hat{z}$ of \eqref{hybridautomaton} has a hybrid time domain that satisfies the ADT bound 
\begin{equation}\label{ADT}
j_2-j_1\leq \frac{1}{\tau_d}(s_2-s_1)+N_0.
\end{equation}
for all $(s_1,j_1),(s_2,j_2) \in\text{dom}(\hat{z})$, with $s_2>s_1\geq0$. Additionally, by \cite[Ex. 2.15]{bookHDS}, for every hybrid time domain satisfying \eqref{ADT}, there exists a solution of the HDS \eqref{hybridautomaton} having said hybrid time domain. Thus, it remains to show that \eqref{ADT} is equivalent to \eqref{BUADTk} in the original $(t,j)$-time scale.
Using the time scaling function $\mathcal{T}_k$ given by \eqref{inverseformula2}, for any solution $z$ of \eqref{hybridBlowUpAutomatonStable} and all  $(t_1,j_1),(t_2,j_2)\in \text{dom}(z)$ with $0\leq t_1<t_2$, we have that $(s_1,j_1),(s_2,j_2)\in \text{dom}(\hat{z})$, where $s_1=\mathcal{T}_k(t_1)$, $s_2=\mathcal{T}_k(t_2)$, and $0\leq s_1<s_2$. Substituting in \eqref{ADT}:
\begin{align*}
j_2-j_1&\leq \frac{1}{\tau_d}(\mathcal{T}_{k}(t_2)-\mathcal{T}_{k}(t_1))+N_0\\
&=\frac{1}{\tau_d}\omega_k(\mu_k(t_2),\mu_k(t_1))+N_0,
\end{align*}
where the equality follows by using (P2) in Proposition \ref{transformationk}. \hfill $\blacksquare$

\vspace{0.1cm}
One of the main consequences of Lemma \ref{lemmaequivalencedomains} is that studying the stability properties of the R-Switching system \eqref{vectorfield1}-\eqref{resetmap1}, under the family of switching signals $\Sigma_{\text{BU-ADT}}(\tau_d,N_0,T,\mu_0)$, is equivalent to studying the stability properties of the HDS \eqref{mainHDSmodel} with $(F_{\Psi},G_{\Psi},\Psi_C,\Psi_{D})$ given by \eqref{mapsstable01}-\eqref{sets01}. For this system, we study stability with respect to the set $\mathcal{A}$ given by \eqref{setstable}, where $\mathcal{A}_1$ is the following compact set
\begin{equation}\label{stableset1}
\mathcal{A}_1=\{0\}\times [0,N_0]\times \mathcal{Q}_s.
\end{equation} 
The following Theorem is the first main result of this paper. The proof leverages the dilation/contraction property of HTDs, established in Proposition \ref{dilationProposition}, as well as the ratio
\begin{equation}\label{omegaconstant}
r=\frac{\max_{q\in\mathcal{Q}} c_{q,2}}{\min_{q\in \mathcal{Q}}{c_{q,1}}},
\end{equation}
which satisfies $r\geq1$ by Assumption \ref{assumptionunstablemodes}, and which is common in the analysis of switched systems.

\vspace{0.1cm}
\begin{thm}\label{theorem1}
Let $N_0\geq1$, $\mathcal{Q}_u=\emptyset$, and consider the HDS \eqref{mainHDSmodel} with $(F_{\Psi},G_{\Psi},\Psi_C,\Psi_D)$ given by \eqref{mapsstable01}-\eqref{sets01}. Suppose that Assumptions \ref{assumptionLipschitz}-\ref{assumptionunstablemodes} hold, and
\begin{equation*}%
\tau_d>\frac{\ln(r)}{\min_{q\in \mathcal{Q}}c_{q,3}}.
\end{equation*}
For each $(T,k)\in\mathbb{R}_{\geq0}\times\mathbb{R}_{\geq1}$, the following holds:
\begin{enumerate}[~~(a)]
\item If $\Delta(\mu_k) = 0$, then  $\mathcal{A}$ is  PT-S$_{\text{F}}$.
\item If $\Delta(\mu_k)=1$, then  $\mathcal{A}$ is PT-ISS$_{\text{F}}$. 
\item If $\Delta(\mu_k)=\mu_k^{-\ell}$, then $\mathcal{A}$ is PT-ISS-C$_{\text{F}}$.\QEDB
\end{enumerate}
\end{thm}
\begin{figure*}[t!]
    \centering
    \includegraphics[width=.45\linewidth]{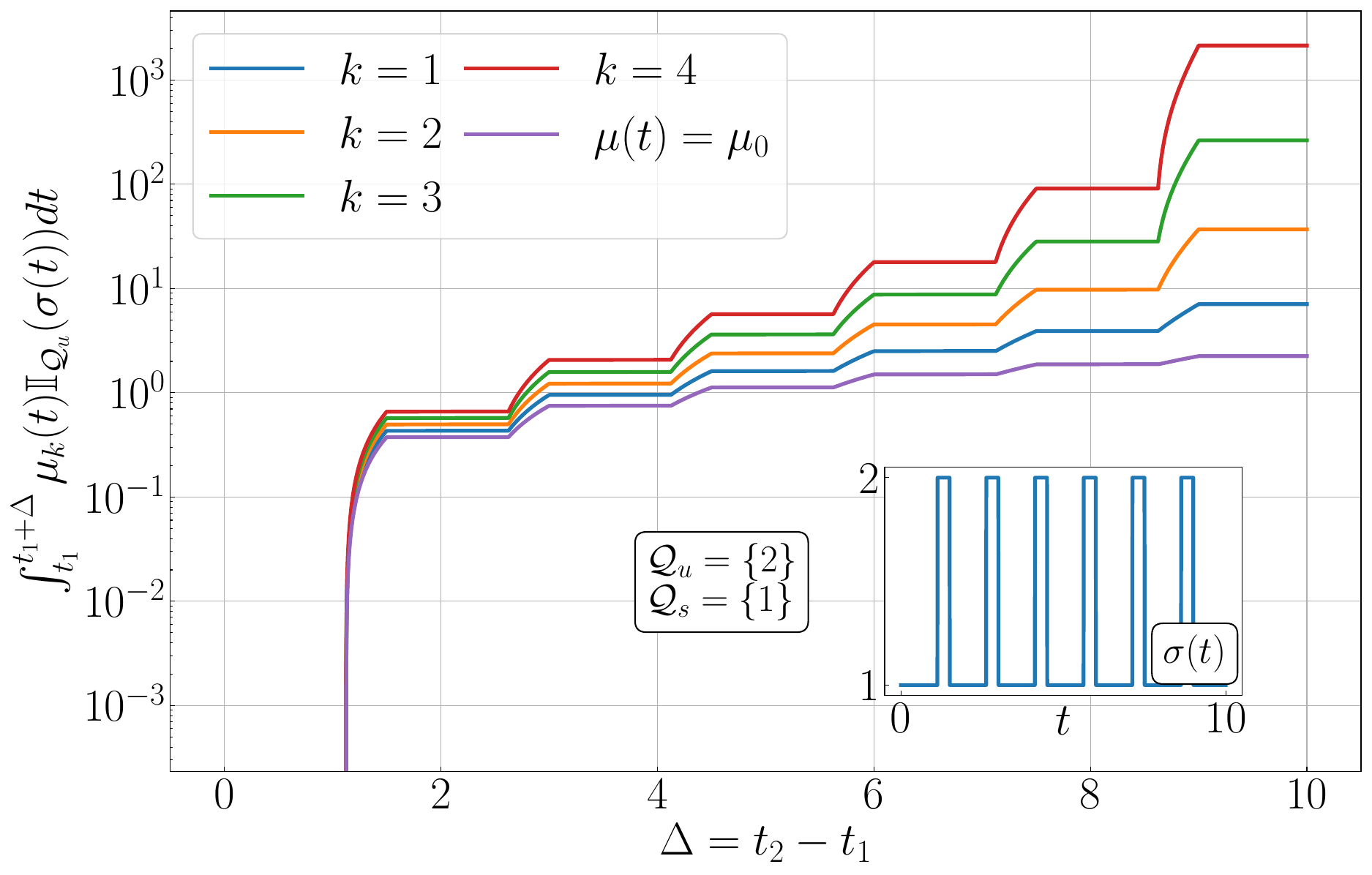} \includegraphics[width=.45\linewidth]{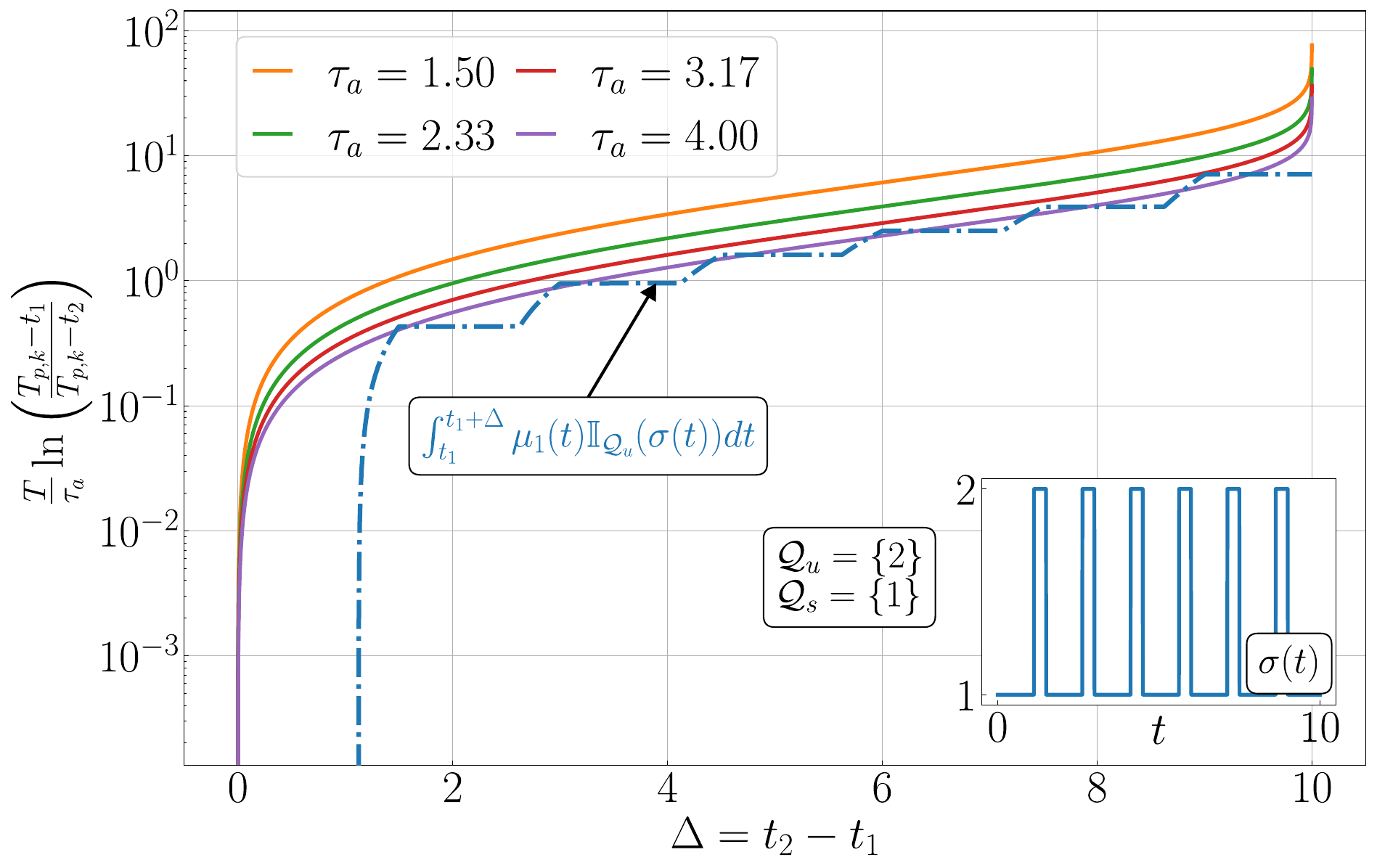}
    \caption{Functions appearing in the BU$_k$-AAT condition \eqref{BUkAAT} using the switching signal $\sigma(\cdot)$ (see inset), $T=10$, and $\mu_0=1$.}\vspace{-0.4cm}
    \label{fig:AATmu01}
\end{figure*}

\vspace{0.1cm}
\noindent
\textbf{Proof}: The proof has three main steps. %

\vspace{0.05cm}
\noindent 
\textsl{Step 1: Stability of the HDS in the $(s,j)$-Hybrid Time Scale:} The overall HDS is given by \eqref{hybridBlowUpAutomatonStable}, which in the $(s,j)$-time scale is given by \eqref{hybridautomaton}. We proceed to study the stability properties of this system with respect to the set $\mathcal{A}$.

By construction of $C$ and $D$, we have that $|\hat{z}|_{\mathcal{A}}=|\hat{x}|$ for all $\hat{z}\in C\cup D$. Thus, it suffices to study the stability properties of $\hat{x}$ with respect to the origin. To do this, we consider the Lyapunov function $W(\hat{z}):=V_{\hat{q}}(\hat{x},\hat{\tau})e^{\ln(r) \hat{\tau}}$. Using Assumption \ref{assumptionunstablemodes}, this function satisfies
\begin{equation*}
\underline{c}|\hat{z}|^p_{\mathcal{A}}\leq W(\hat{x})\leq \overline{c}|\hat{z}|^p_{\mathcal{A}},~~\forall~\hat{z}\in C\cup D,
\end{equation*}
with $\underline{c}\!\coloneqq\!\min_{p\in \mathcal{Q}}c_{1,p}$, $\overline{c}\!\coloneqq\! e^{\ln(r) N_0}\overline{c}_2$, and $\overline{c}_2\!\coloneqq\! \max_{p\in \mathcal{Q}}c_{2,p}$. When $\hat{z}\in C$, for all $\eta\in[0,1/\tau_d]$,  we have:
\begin{align*}
\langle \nabla W(\hat{z}), \hat{F}_T(\hat{z},\hat{u})\rangle&=\left\langle \nabla V_{\hat{q}}(\hat{x},\hat{\tau}), \left(\begin{array}{c}
f_{\hat{q}}(\hat{x},\hat{u},\hat{\tau})\\
\eta
\end{array}\right)
\right\rangle e^{\ln(r)\hat{\tau}}\\
&~~~~+\langle \ln(r) V_{\hat{q}}(\hat{x},\hat{\tau})e^{\ln(r)\hat{\tau} } ,\dot{\hat{\tau}}_s\rangle\notag\\
\le -\underline{c}_3\bigg(1&-\frac{\ln(r)}{\underline{c}_3\tau_d}\bigg) W(\hat{z})+\overline{c}_4e^{\ln(r) N_0}\Delta(\hat{\mu}_k)|\hat{u}|^{p}
\end{align*}
where $\underline{c}_3:=\min_{p\in\mathcal{Q}}c_{p,3}$ and $\overline{c}_4:=\max_{p\in \mathcal{Q}}c_{p,4}$. On the other hand, when $\hat{z}\in D$
\begin{align*}
W(\hat{z}^+)&=V_{\hat{q}^+}(\hat{x}^+,\hat{\tau}^+)e^{\ln(r) \hat{\tau}^+}\leq \chi \max_{p\in\mathcal{Q}}V_{p}(\hat{x},\hat{\tau})e^{\ln(r) (\hat{\tau}-1)}\\
&\leq \overline{c}_2\chi|\hat{x}|^pe^{\ln(r) (\hat{\tau}-1)}\leq \chi\frac{\overline{c}_2}{\underline{c}}e^{-\ln(r)}V_{\hat{q}}(\hat{x},\hat{\tau})e^{\ln(r) \hat{\tau}}.
\end{align*}
Thus, during jumps:
\begin{equation*}
W(\hat{z}^+)-W(\hat{z})\le  -\left(1 {-} \chi\frac{\overline{c}_2}{\underline{c}\cdot r}\right)W(\hat{z}) = -\left(1{-}\chi\right)W(\hat{z})\le 0.
\end{equation*}
Using Lemma 7 in the Supplemental Material, we conclude that every solution $\hat{z}$ of system \eqref{hybridautomaton} satisfies:
\begin{equation}\label{exponential_bound0B}
|\hat{z}(s,j)|_{\mathcal{A}}\leq \kappa_1 e^{-\kappa_2 (s+j)}|\hat{z}(0,0)|_{\mathcal{A}}+\kappa_3  \cdot \sup_{0\leq \zeta\leq s}|\hat{\Delta}(\zeta)|, %
\end{equation}
for all $(s,j)\in\text{dom}(\hat{z})$, where $\kappa_1=\left(\overline{c}/\underline{c}\right)^{1/p}e^{\frac{\lambda}{2p} \frac{\tau_d}{1+\tau_d}N_0}$ $\kappa_{2}=\lambda\tau_d/(2p(1+\tau_d))$,  $\kappa_3 = \left(2\overline{c}_4r^{N_0}/[\lambda \underline{c}]\right)^{1/p}$, $\lambda = \underline{c}_3-\ln(r)/\tau_d$, and $\hat{\Delta}(s):=\Delta(\hat{\mu}_k(s))\hat{u}(s)$. 
Moreover, when $\Delta(\hat{\mu}_k)=\frac{1}{\hat{\mu}_k^{\ell}}$, via Lemma 8 in the Supplemental Matierial, there exists  $\beta_k(r,s)\!=\!r\min\{e^{-\kappa_2s},\mu_k^{-\ell}(s)\}\in\mathcal{KL}$ such that every solution $\hat{z}$ of system \eqref{hybridautomaton} satisfies:
\begin{equation}\label{ISSPCbound}
|\hat{z}(s,j)|_{\mathcal{A}}\!\leq\!\beta_k\Big(\overline{\kappa}_1|\hat{z}(0,\!0)|_{\mathcal{A}}e^{\!-\overline{\kappa}_2(s+j)}
                        \!+ \overline{\kappa}_{3}|\hat{u}|_{(s,j)},\! s\Big),
\end{equation}
for all $(s,j)\in\text{dom}(\hat{z})$, with  $\overline{\kappa}_1 \coloneqq \kappa_1^2$, $\kappa_2\coloneqq \frac{\kappa_2}{2}$, $\overline{\kappa}_{3} \coloneqq 2\kappa_1\kappa_3$.

\vspace{0.2cm}
\noindent 
\textsl{Step 2: PT-ISS$_{\text{F}}$ of the HDS in the $(t,j)$ - Time Scale:} We now use the properties of the solutions $\hat{z}$ of system \eqref{hybridautomaton} to establish properties for the solutions $z$ of system \eqref{hybridBlowUpAutomatonStable}, under a given $\mu_k(0,0)=\mu_0\in\mathbb{R}_{\geq1}$.
Using Proposition \ref{dilationProposition}, $s=\mathcal{T}_k(t)$,  and \eqref{inverseformula2}, we obtain $e^{-\kappa_2 (s+j)}=e^{-\kappa_2(\mathcal{T}_k(t)+j)}$, and since $\mathcal{T}_k(0)=\mathcal{T}_k^{-1}(0)=0$, $|\hat{z}(\mathcal{T}_k(t),j)|_{\mathcal{A}}=|z\left(\mathcal{T}_k^{-1}(\mathcal{T}_k(t)),j\right)|_{\mathcal{A}}=|z(t,j)|_{\mathcal{A}}$, and substituting in \eqref{exponential_bound0B}, it follows that when $\Delta=0$ or $\Delta=1$, every solution $z$ of the HDS \eqref{hybridBlowUpAutomatonStable} with $\mu_k(0,0)=\mu_0\geq1$ satisfies:
\begin{equation}\label{zbound1}
|z(t,j)|_{\mathcal{A}}\leq \kappa_1 e^{-\kappa_2 (\mathcal{T}_k(t)+j)}|z(0,0)|_{\mathcal{A}}+\kappa_3  \Delta|u|_{(t,j)}, 
\end{equation}
for all $(t,j)\in\text{dom}(z)$, which implies that $\mathcal{A}$ is PT-ISS$_F$. Similarly, when $\Delta(\hat{\mu}_k)=\hat{\mu}_k^{-\ell}$, \eqref{ISSPCbound} leads to:

\vspace{-0.2cm}
\begin{small}
\begin{equation}\label{PTISCproof2}
\!\!\!\!|z(s,j)|_{\mathcal{A}}\!\leq\!\beta_k\!\Big(\overline{\kappa}_1|z(0,\!0)|_{\mathcal{A}}e^{\!{-}\overline{\kappa}_2\left(\mathcal{T}_k(t){+}j\right)}
                        {+} \overline{\kappa}_{3}|u|_{(t,j)},\!\mathcal{T}_k(t)\Big),
\end{equation}
\end{small}

\noindent 
for all $(t,j)\in\text{dom}(z)$, which implies that $\mathcal{A}$ is PT-ISS-C$_F$.

\vspace{0.2cm}
\noindent 
\textsl{Step 3: Length of solutions in the $(t,j)$ - Time Scale:} Finally, we show that
$\sup_t (\text{dom}(z))=\Upsilon_{T,k}$ for all solutions $z$ of \eqref{hybridBlowUpAutomatonStable}. First, note that by the definition of $\mathcal{T}_k$ and Proposition \ref{dilationProposition}, we have $\sup_t (\text{dom}(z))=\Upsilon_{T,k}$ if and only if $\sup_s (\text{dom}(\hat{z}))=\infty$. Furthermore, based on the bound \eqref{ADT}, we obtain hat $j\leq \frac{1}{\tau_d}s+N_0$ for any $(s,j)\in \text{dom}(\hat{z})$. Considering that every complete solution $\hat{z}$ of \eqref{hybridautomaton} satisfies $\text{length}(\text{dom}(\hat{z}))=\infty$, and noting that $\text{length}(\text{dom}(\hat{z}))=\sup_s (\text{dom}(\hat{z}))+\sup_j(\text{dom}(\hat{z}))$, we can infer that if $j\to\infty$, then $s\to\infty$. Consequently, every complete solution of \eqref{hybridautomaton} must satisfy $\sup_s (\text{dom}(\hat{z}))=\infty$, which in turn implies that $\sup_t (\text{dom}(z))=\Upsilon_{T,k}$. \hfill $\blacksquare$

\vspace{0.2cm}
The following Corollary covers the case $k=1$ (and $\mu_0=1$) which is the most common in the literature of PT-S \cite{song2017time,orlov2021prescribed}.  %

\vspace{0.1cm}
\begin{cor}
Suppose that all the assumptions of Theorem \ref{theorem1} hold, and that $k=1$. Then, for every solution $z$ and for all $(t,j)\in\text{dom}(z)$, the trajectory $x$ satisfies the following:
\begin{enumerate}
\item If \eqref{stablestabilitycondition} holds with $\Delta(\mu_1)=0$ or $\Delta(\mu_1)=1$:
\begin{equation}\label{stability:k=1}
\!\!\!\!\!\!|x(t,j)|\!\leq\! \kappa_1\!\left(\frac{\mu_0}{\mu_1(t)}\right)^{\!\!\kappa_2T}\!\!\!\!\!e^{-\kappa_2j}|x(0,0)|+\kappa_3\Delta |u|_{(t,j)}, 
\end{equation}
for $\kappa_i>0$.
\item If \eqref{stablestabilitycondition} holds with $\Delta(\mu_1)=\mu_1^{-\ell}$:
\begin{equation}\label{stability+convergence:k=1}
\!\!\!\!\!\!\!\!\left|x(t,j)\right|\leq\frac{\alpha_1\mu_0^{\alpha_2}}{\mu_1(t)^{\alpha_3}}\!\left(\frac{e^{-\alpha_4j}}{\mu_1(t)^{\alpha_5}}|x(0,0)|+\alpha_6|u|_{(t,j)}\right),
\end{equation}
for $\alpha_i>0$. \hfill \QEDB
\end{enumerate}
\end{cor}

\vspace{0.1cm}
\textbf{Proof:} Using \eqref{lnk1} and the bounds obtained in Step 2 of the proof of Theorem \ref{theorem1}, it follows that $e^{-\kappa_2(\mathcal{T}_k(t)+j)}=e^{-\alpha\ln\left(\frac{\mu_1(t)}{\mu_0}\right)}e^{-\kappa_2 j}=\left(\frac{\mu_0}{\mu_1(t)}\right)^{\alpha}e^{-\frac{\alpha}{T}j}$
where $\alpha=\kappa_2 T$. Since $|z|_{\mathcal{A}}=|x|$ for every solution, inequality \eqref{zbound1} becomes \eqref{stability:k=1}. 
Similarly, inequality \eqref{PTISCproof2} becomes \eqref{stability+convergence:k=1} with 
$\alpha_1\!:=\!\overline{\kappa}_1$, $\alpha_2\!:=\!(\overline{\kappa}_3\!+\!\overline{\kappa}_2)T$, $\alpha_3\!:=\!\overline{\kappa}_3T$, $\alpha_4\!:=\!\overline{\kappa}_2$, $\alpha_5\!:=\!\overline{\kappa}_2T$,  and $\alpha_6\!:=\!\overline{\kappa}_4$. 
\hfill $\blacksquare$

\subsection{PT-ISS in R-Switching Systems with Unstable Modes}
We now consider the situation when some of the modes $f_q$ in \eqref{actualflow1} are unstable, i.e., $\mathcal{Q}_u\neq \emptyset$ and $\mathcal{Q}=\mathcal{Q}_s\cup \mathcal{Q}_u$. In this case, we introduce the following \emph{blow-up average activation-time} (BU$_{k}$-AAT) condition on the amount of time that the unstable modes can remain active in any given window of time contained in $[0,\Upsilon_{T,k})$:

\vspace{0.1cm}
\begin{definition}
A switching signal $\sigma:[0,\Upsilon_{T,k})\to\mathcal{Q}$ is said to satisfy the \emph{blow-up average activation-time condition of order $k$} (BU$_k$-AAT) if there exist $T_0\ge 0$ and $\tau_a>1$ such that for each pair of times $t_2,t_1\in\text{dom}(\sigma)$:
\begin{equation}\label{BUkAAT}
\int_{t_1}^{t_2}\mu_k(t)\cdot \mathbb{I}_{\mathcal{Q}_u}(\sigma(t))dt \leq \frac{1}{\tau_a}\omega_k\left(\mu_k(t_2),\mu_k(t_1)\right)+T_0,
\end{equation}
where $\mu_k$ is given by \eqref{tpmu}. We denote the family of such signals as $\Sigma_{\text{BU$_k$-AAT}}(\mathcal{Q}_u,\tau_a,T_0,T,\mu_0)$.\QEDB
\end{definition}

\vspace{0.1cm}
\begin{remark}
For asymptotic and exponential stability results in switching systems with stable and unstable modes \cite{GuosongLiberzon,poveda2017framework,Liu_Tanwani_Liberzon2022}, it is common to restrict the family of admissible switching signals to those that satisfy the ADT condition \eqref{ADT} \emph{and} the following average activation-time (AAT) condition
\begin{equation}\label{eq:AAT}
        \int_{t_1}^{t_2} \mathbb{I}_{\mathcal{Q}_u}(\sigma(t))dt \leq \frac{1}{\tau_a}(t_2-t_1)+T_0,
\end{equation}
where $\tau_a>1$, and $T_0\geq0$, which is recovered by taking the limit as $T\to\infty$ in both sides of \eqref{BUkAAT}, and using $\mu_0=1$. \hfill \QEDB
\end{remark}

\vspace{0.1cm}
For $k=1$, the BU$_1$-AAT condition reduces to:
\begin{equation*}%
        \int_{t_1}^{t_2}\frac{\mathbb{I}_{\mathcal{Q}_u}(\sigma(t))}{T-t\mu_0}dt \leq \frac{1}{\tau_a\mu_0} \ln\left(\frac{T-t_1\mu_0}{T-t_2\mu_0}\right)+T_0.
\end{equation*}
Similar expressions can be obtained for $k\in\mathbb{Z}_{\geq2}$ using equations \eqref{tpmu} and \eqref{omegaformula}. Figure \ref{fig:AATmu01} compares the BU$_k$-AAT bounds obtained in \eqref{BUkAAT} and the traditional AAT bound \eqref{eq:AAT}. The left plot shows the left-hand side of \eqref{BUkAAT} for different values of $k$, under a particular switching signal $\sigma$ that switches between one stable mode and one unstable mode (see inset). The classic AAT bound is shown in purple color. Similarly, the right plot shows \eqref{BUkAAT} for $k=1$ and for different values of $\tau_a$. It can be seen that for small values of $\tau_a$, inequality \eqref{BUkAAT} is satisfied. \QEDB

\vspace{0.1cm}
To study the PT-S properties of the R-Switching system \eqref{vectorfield1}-\eqref{resetmap1}  when $\mathcal{Q}$ contains unstable modes, we now consider the HDS \eqref{mainHDSmodel} with state $\psi=(x,\tau,\rho,q)\in\mathbb{R}^{n+3}$, set-valued mappings:
\begin{subequations}\label{hdsUnstable}
\begin{align}\label{mapsunstable01}
&F_{\Psi}\!:=\!\{f_q(x,u,\tau)\}{\times}\!
\left[0,\frac{1}{\tau_d}\right]\!{\times}\!
\left[0,\frac{1}{\tau_a}\right]\!\!{-}\mathbb{I}_{\mathcal{Q}_u}(q){\times}\{0\},\\\
&G_{\Psi}\!:=\!\{R_q(x)\}\times\{\tau-1\}\times\{\rho\}\times\mathcal{Q}\backslash\{q\},
\end{align}
and sets:
\begin{align}\label{sets02}
\Psi_C&=\mathbb{R}^n\times[0,N_0]\times[0,T_0]\times\mathcal{Q},\quad
\Psi_D=\mathbb{R}^n\times[1,N_0]\times[0,T_0]\times\mathcal{Q}.
\end{align}
\end{subequations}
There is a close connection between the hybrid time domains of the solutions generated by the HDS \eqref{mainHDSmodel} with data \eqref{hdsUnstable} and the switching signals that simultaneously satisfy \eqref{BUADTk} and \eqref{BUkAAT}.

\vspace{0.1cm}
\begin{lemma}\label{lemmaequivalencedomains2}
Let $(F_{\Psi},G_{\Psi},\Psi_C,\Psi_D)$ be given by \eqref{mapsunstable01}-\eqref{sets02}, and consider the HDS \eqref{mainHDSmodel} under Assumption \ref{assumptionLipschitz}-\ref{assumptionunstablemodes}. Then, Assumption \ref{wellposedassumption} holds, and:
\begin{enumerate}[(a)]
\item For every maximal solution $z$ and for any pair $(t_1,j_1),(t_2,j_2)\in \text{dom}(z)$, with $t_2>t_1$, inequality \eqref{BUADTk} holds with  
$N(t_2,t_1)=j_2-j_1$, and inequality \eqref{BUkAAT} holds with $\sigma(t)=q(t,\underline{j}(t))$, where $\underline{j}(t):=\min\{j\in\mathbb{Z}_{\geq0}:(t,j)\in\text{dom}(z)\}$.
\item  For every hybrid time domain satisfying property (a), there exists a solution $z$ of \eqref{mainHDSmodel} having the said hybrid time domain.  \QEDB
\end{enumerate}
\end{lemma}
\noindent\textbf{Proof:} The overall HDS has state $z\!=\!(x,\tau,\rho,q,\mu_k)\!\in\!\mathbb{R}^{n+4}$ and data:
\begin{subequations}\label{hybridBlowUpAutomatonStable20}
\begin{align}\label{flowshybrid2}
&z\in C:=\mathbb{R}^n\times [0, N_0]\times[0,T_0]\times\mathcal{Q}\times\mathbb{R}_{\geq1},\quad
\dot{z}=\left(
\begin{array}{c}
\vspace{0.45cm}
\dot{x}\\\vspace{0.45cm}
\dot{\tau}\\\vspace{0.45cm}
\dot{\rho}\\\vspace{0.45cm}
\dot{q}\\\vspace{0.45cm}
\dot{\mu}_k
\end{array}\right)\in F(z,u):=\left(\begin{array}{c}
    \vspace{0.2cm}
    \mu_k \cdot f_q(x,u,\tau)\\\vspace{0.1cm}
    \left[0,\dfrac{\mu_k}{\tau_d}\right]\vspace{0.1cm}\\
    \!\!\!\!\left[0,\dfrac{\mu_k}{\tau_a}\right]-\mu_k\mathbb{I}_{\mathcal{Q}_u}(q)\!\!\!\!\vspace{0.1cm}\\
    0\vspace{0.3cm}\\
    \dfrac{k}{T}\mu_k^{1+\frac{1}{k}}\vspace{0.3cm}
    \end{array}\right),
\end{align}
\begin{align}
&z\in D:=\mathbb{R}^n\times [1, N_0]\times[0,T_0]\times\mathcal{Q}\times\mathbb{R}_{\geq1},\quad
z^+=\left(\begin{array}{c}\
    \vspace{0.1cm}
    x^+\\\vspace{0.1cm}
    \tau^+\\\vspace{0.1cm}
    \rho^+\\\vspace{0.1cm}
    q^+\\\vspace{0.1cm}
    \mu_k^+\vspace{0.1cm}
    \end{array}\right)\in G(z,u):=\left(
    \begin{array}{c}
    \vspace{0.1cm}
    R_q(x)\vspace{0.1cm}\\
    \tau-1\vspace{0.1cm}\\
    \rho\vspace{0.1cm}\\
    \mathcal{Q}\vspace{0.1cm}\\
    \mu_k
    \end{array}\right).\label{jumpshybrid2}
\end{align}
\end{subequations}
This system has a finite escape time induced by the state $\mu_k$ which diverges to infinity during flows as $t\to \Upsilon_{T,k}$. Note that, by construction, the states $(\tau,
\rho,q)$ are confined to the compact sets $[0,N_0]$, $[0,T_0]$, and $\mathcal{Q}$ 
respectively.  Using again the time variable $s=\mathcal{T}_k(t)$ defined in \eqref{inverseformula2}, and Proposition \ref{dilationProposition}, we obtain the following HDS in the $(s,j)$-time scale:
\begin{subequations}\label{hybridBlowUpAutomatonUnstable}
\begin{align}
&\hat{z}\!\in\!C,~\dot{\hat{z}}_s \!=\! \begin{pmatrix}
\vspace{0.45cm}
\dot{\hat{x}}_s\\\vspace{0.45cm}
\dot{\hat{\tau}}_s\\\vspace{0.45cm}
\dot{\hat{\rho}}_s\\\vspace{0.45cm}
\dot{\hat{q}}_s\\\vspace{0.45cm}
\dot{\hat{\mu}}_{k_s}
\end{pmatrix}\!\!\in\!
\hat{F}(\hat{z},\hat{u}):=\!\left(\begin{array}{c}
\vspace{0.2cm}
f_q(\hat{x},\hat{u},\hat{\tau})\\\vspace{0.1cm}
\left[0,\dfrac{1}{\tau_d}\right]\vspace{0.1cm}\\
\!\!\!\!\left[0,\dfrac{1}{\tau_a}\right]\!\!-\!\mathbb{I}_{\mathcal{Q}_u}(\hat{q})\!\!\!\!\vspace{0.1cm}\\
0 \vspace{0.3cm}\\
 \dfrac{k}{T}\hat{\mu}_k^{\frac{1}{k}}\vspace{0.3cm}
\end{array}\right),\label{flowsHDS2}\\
&\hat{z}\in D,~~~~~~\hat{z}^+\in G(\hat{z}),
\end{align}
\end{subequations}
where the subscript $s$ in \eqref{flowsHDS2} indicates that the time derivative is taken with respect to $s$.  Since \eqref{hybridBlowUpAutomatonUnstable} incorporates an ADT automaton $\hat{\tau}$ and a time-ratio monitor $\hat{\rho}$, by \cite[Lemma 7]{poveda2017framework} every solution $\hat{z}$ of \eqref{hybridBlowUpAutomatonUnstable} has a hybrid time domain such that for any pair $(s_1,j_1),(s_2,j_2)\in\text{dom}(\hat{z})$ the bound \eqref{ADT} is satisfied, along with the bound:
\begin{align}\label{proofAAT}
 \mathbb{T}(s_1,s_2){:=}\int_{s_1}^{s_2}\!\!\!\mathbb{I}_{\mathcal{Q}_u}(\hat{q}(s,\underline{\hat{\jmath}}(s)))ds \leq \frac{1}{\tau_a} (s_2-s_1)+T_0,
\end{align}
where $\underline{\hat{\jmath}}(s)\coloneqq \min\left\{j\in\mathbb{Z}_{\ge 0}~:~(s,j)\in\text{dom}(\hat{q})\right\}$.
Recall that, by Lemma \ref{equivalentBUADTk}-c), as $T\to\infty$, the BU$_k$-ADT condition \eqref{BUADTk} converges to the ADT condition \eqref{ADT} in the $(t,j)$-time scale. Similarly, using $s=\mathcal{T}_k(t)$, the left-hand side of \eqref{proofAAT} can be expressed in the $t$-variable as:
\begin{align}
\mathbb{T}(\mathcal{T}_k(t_2),\mathcal{T}_k(t_1))&{=}\int_{t_1}^{t_2}\!\!\frac{\partial \mathcal{T}_k(t)}{\partial t}\cdot \mathbb{I}_{\mathcal{Q}_u}\bigg(\hat{q}\Big(\mathcal{T}_k(t),\hat{\jmath}\big(\mathcal{T}_k(t)\big)\Big)\bigg)dt\notag\\
&=\int_{t_1}^{t_2}\mu_k(t)\cdot \mathbb{I}_{\mathcal{Q}_u}\big(q(t,\underline{j}(t))\big)dt,\label{integralAATtscale}
\end{align}
where we used Proposition \ref{transformationk}-(P3), together with the equality $\hat{q}(\mathcal{T}_k(t), \underline{\hat{\jmath}}\left(\mathcal{T}_k(t)\right))=q\left(\mathcal{T}_k^{-1}(\mathcal{T}_k(t)),~\underline{j} \left(\mathcal{T}^{-1}_k\left(\mathcal{T}_k(t)\right)\right)\right)=q(t,\underline{j}(t))$. Using \eqref{integralAATtscale} in \eqref{proofAAT}, together with Proposition \ref{transformationk}-(P2), it follows that the AAT condition in the $(t,j)$-time scale becomes
 \begin{equation*}
\int_{t_1}^{t_2}\mu_k(t)\cdot \mathbb{I}_{\mathcal{Q}_u}(q(t,\underline{j}(t)))dt \leq \frac{1}{\tau_a}\omega_k\left(\mu_k(t_2),\mu_k(t_1)\right)+T_0,
\end{equation*}
which is precisely \eqref{BUkAAT}. \hfill $\blacksquare$

\vspace{0.1cm}
Similar to Lemma \ref{lemmaequivalencedomains}, the result of Lemma \ref{lemmaequivalencedomains2} enables the study of the stability properties of the R-Switching system \eqref{vectorfield1}-\eqref{resetmap1}, under switching signals $\sigma$ satisfying \eqref{BUkAAT}, by studying the stability properties of the HDS \eqref{hybridBlowUpAutomatonStable20}. In this case, we consider the set $\mathcal{A}$ given by \eqref{setstable}, where $\mathcal{A}_1$ is now given by
\begin{equation}\label{stableset2}
\mathcal{A}_1=\{0\}\times [0,N_0]\times [0,T_0]\times \mathcal{Q}.
\end{equation}
The next theorem is the second main result of this paper.

\vspace{0.1cm}

\vspace{0.1cm}
\begin{thm}\label{theorem2}
Let $N_0\geq1$, $T_0\geq0$, $\mathcal{Q}_u\neq \emptyset$, and consider the HDS \eqref{mainHDSmodel} with $(F_{\Psi},G_{\Psi},\Psi_C,\Psi_D)$ given by \eqref{mapsunstable01}-\eqref{sets02}. Suppose that Assumptions \ref{assumptionLipschitz}-\ref{assumptionunstablemodes} hold, and that
\begin{equation}\label{conditionstability2}
1 > \frac{1}{\underline{c}_3\tau_d}\log(r) + \frac{1}{\tau_a}\left(1+\frac{\overline{c}_5}{\underline{c}_3}\right).
\end{equation}
where $r$ is given by \eqref{omegaconstant}, $\underline{c}_{3}\!=\!\min_{p\in\mathcal{Q}}\!c_{q,3}$, $\overline{c}_{5}\!=\!\max_{p\in\mathcal{Q}}\!c_{q,5}$. For each $(T,k)\in\mathbb{R}_{>0}\times\mathbb{R}_{\geq1}$ the following holds:
\begin{enumerate}[~~(a)]
\item If $\Delta(\mu_k)\triangleq0$, then $\mathcal{A}$ is  PT-S$_{\text{F}}$.
\item If $\Delta(\mu_k)\triangleq1$, then  $\mathcal{A}$ is PT-ISS$_{\text{F}}$. 
\item If $\Delta(\mu_k)\triangleq\mu_k^{-\ell}$, $\ell>0$, then  $\mathcal{A}$ is PT-ISS-C$_{\text{F}}$. 
\QEDB
\end{enumerate}
\end{thm}

\vspace{0.1cm}
\textbf{Proof:} The proof follows the same three steps as in the proof of Theorem \ref{theorem1}. Indeed,  we start by  using the time dilation $\mathcal{T}^{-1}_k$ and Proposition \ref{dilationProposition}. Hence, we consider the HDS \eqref{hybridBlowUpAutomatonUnstable} in the $(s,j)$-time scale, with state $\hat{z}=(\hat{x},\hat{\tau},\hat{\rho},\hat{q},\hat{\mu}_k)$. To study the stability properties of this system, let 
\begin{equation*}
    \hat{\xi}:=\log(r)\hat{\tau}+(\underline{c}_3+\overline{c}_5)\hat{\rho},
\end{equation*}
and consider the Lyapunov function $W_2(\hat{z})=V_{\hat{q}}(\hat{x},\hat{\tau})e^{\hat{\xi}}$, which, by Assumption \ref{assumptionunstablemodes}-(a), satisfies
\begin{equation*}%
\underline{\varphi}|\hat{z}|^2_{\mathcal{A}}\leq W_2(\hat{z})\leq \overline{\varphi}|\hat{z}|^2_{\mathcal{A}},
\end{equation*}
with $\underline{\varphi}=\min_{p\in\mathcal{Q}}c_{p,1},~\overline{\varphi}=\max_{p\in\mathcal{Q}}c_{p,2}e^{\log(r)N_0+(\underline{c}_3+\overline{c}_5)T_0}$.
When $\hat{z}\!\in\! C$, the time derivative of $\hat{\xi}$ with respect to $s$ satisfies:
\begin{align*}
    \dot{\hat{\xi}}_s&=\log(r)\dot{\hat{\tau}}_s+(\underline{c}_3+\overline{c}_5)\dot{\hat{\rho}}_s
    \in[0,\delta]-(\underline{c}_3+\overline{c}_5)\mathbb{I}_{\mathcal{Q}_u}(\hat{q})
\end{align*}
where $\delta:=\frac{1}{\tau_d}\log(r)+\frac{1}{\tau_a}(\underline{c}_3+\overline{c}_5)$. Using the above expression together with Assumption \ref{assumptionunstablemodes}, we evaluate the change of $W_2$ during the flows of stable and unstable modes. In particular, when $\hat{z}\in C$ and $\hat{q}\in \mathcal{Q}_s$, we have
\begin{align}
    \langle \nabla W_2(\hat{z}),\dot{\hat{z}}_s\rangle &= e^{\hat{\xi}}\left\langle\nabla V_{\hat{q}}(\hat{x},\hat{\tau}),~ \dot{\hat{x}}_s\right\rangle + e^{\hat{\xi}} V_{\hat{q}}(\hat{x},\hat{\tau})\dot{\hat{\xi}}_s  \notag \\
    &\le -(\underline{c}_3 -\delta)W_2(\hat{z}) + \frac{\overline{c}_4}{\overline{c}_2}\overline{\varphi}\hat{\Delta}(s)|\hat{u}|^{p},\label{issFlows:stableModes}
\end{align}
where $\hat{\Delta}(s):=\Delta(\hat{\mu}_k(s))\hat{u}(s)$, $\overline{c}_2\coloneqq \max_{p\in \mathcal{Q}}c_{2,p}$ and $\overline{c}_4 = \max_{p\in\mathcal{Q}}c_{4,p}$, and where $\underline{c}_3 -\delta>0$ since \eqref{conditionstability2} is satisfied by assumption.
On the other hand, when $\hat{z}\in C$ and $\hat{q}\in \mathcal{Q}_u$:
\begin{align*}
    \langle \nabla W_2(\hat{z}),\dot{\hat{z}}_s\rangle&\leq  \left(\overline{c}_5V_{\hat{q}}(\hat{x},\hat{\tau}) + \overline{c}_4\hat{\Delta}(s)|\hat{u}|\right)e^{\hat{\xi}}+V_{\hat{q}}(\hat{x},\hat{\tau})e^{\hat{\xi}} \dot{\hat{\xi}}_s\notag\\
    &\le\left(\delta-\underline{c}_3\right)W_2(\hat{z}) + \overline{c}_4\hat{\Delta}(s)|\hat{u}|e^{\hat{\xi}}\notag\\
    &\le -\left(\underline{c}_3-\delta\right)W_2(\hat{z}) +\frac{\overline{c}_4}{\overline{c}_2}\overline{\varphi}\hat{\Delta}(s)|\hat{u}|^p
\end{align*}
which is the same bound as \eqref{issFlows:stableModes}. On the other hand, when $\hat{z}\in D$, it follows that
$\hat{\xi}^+=\log(r)\hat{\tau}^+ +(\underline{c}_3+\overline{c}_4)\hat{\rho}^+ %
=\hat{\xi}-\log(r).$
Therefore, during jumps the Lyapunov function satisfies:
\begin{align*}
W_2(\hat{z}^+)%
\le \chi\max_{\hat{q}^+\in \mathcal{Q}} V_{\hat{q}^+}(\hat{x},\hat{\tau})e^{\hat{\xi}-\log(r)}
\le W_2(\hat{z})\notag.
\end{align*}
It follows that $W_2(\hat{z}^+)-W_2(\hat{z})\leq0$ for all $\hat{z}\in D$.  Using Lemma 7 in the Appendix, we conclude that every solution $\hat{z}$ satisfies the bound
\begin{equation*}%
|\hat{z}(s,j)|_{\mathcal{A}}\leq \kappa_1|\hat{z}(0,0)|_{\mathcal{A}}e^{-\kappa_2(s+j)}+\kappa_3\hat{\Delta}(s)|\hat{u}|_{(s,j)},
\end{equation*}
for all $(s,j)\in\text{dom}(\hat{z})$, where $\kappa_1=\left(\overline{\varphi}/\underline{\varphi}\right)^{1/p}e^{\frac{\lambda}{2p} \frac{\tau_d}{1+\tau_d}N_0}$ $\kappa_{2}=\lambda\tau_d/(2p(1+\tau_d))$,  $\kappa_3 = \left( 2\overline{c}_4\overline{\varphi}/[\overline{c}_2\lambda \underline{\varphi}]\right)^{1/p}$, $\lambda = \underline{c}_3-\delta$, and $\hat{\Delta}(s):=\Delta(\hat{\mu}_k(s))\hat{u}(s)$.   The bounds \eqref{eq:PTflows1}-\eqref{eq:PTflows2} are obtained following the exact same arguments used in Steps 2 and 3 of the proof of Theorem  \ref{theorem1}. \hfill $\blacksquare$

\vspace{0.1cm}

\begin{remark}[Switching with Non-PT Unstable Modes]

It is reasonable to consider a situation where the unstable modes in \eqref{actualflow1} do not have time-varying gains, i.e., $\mu_k\equiv1$ when $q\in\mathcal{Q}_u$. In particular, consider a system switching between
\begin{equation*}
\dot{x}=\mu_kf_q(x),~~q\in \mathcal{Q}_s,~~\text{and}~~\dot{x}=f_p(x),~~p\in \mathcal{Q}_u,
\end{equation*}
where the modes in $\mathcal{Q}_s$ satisfy \eqref{stablestabilitycondition}, and the modes in $\mathcal{Q}_u$ satisfy \eqref{unstablestabilitycondition} with $u\equiv 0$. Following the same approach of Theorem \ref{theorem2}, and operating in the $s$-time scale for the flows, we now obtain the following two type of modes: 
\begin{equation*}
\dot{\hat{x}}_s=f_q(\hat{x}),~~q\in \mathcal{Q}_s,~~\text{and}~~\dot{\hat{x}}_s=\frac{1}{\hat{\mu}_k}f_p(x),~~p\in \mathcal{Q}_u.
\end{equation*}
For this system, the same Lyapunov-based analysis can be applied as in the proof of Theorem \ref{theorem2} to obtain the bound \eqref{issFlows:stableModes} for all the stable modes. On the other hand, during unstable modes, we now obtain
\begin{align*}
\langle \nabla W_2(\hat{z}),\dot{\hat{z}}_s\rangle %
&\leq -\left(\underline{c}_3-\delta\right)W_2(\hat{z})-\overline{c}_5\left(1-\frac{1}{\hat{\mu}_k}\right)W_2(\hat{z})
\end{align*}
Note that $1-\frac{1}{\hat{\mu}_k}\!\geq\! 0$ since $\hat{\mu}_k\!\geq\! 1$ by Proposition \ref{completemulemma}. This implies that $\langle \nabla W_2(\hat{z}),\dot{\hat{z}}_s\rangle \!\leq\!-\left(\underline{c}_3\!-\!\delta\right)W_2(\hat{z})$. From here, the proofs follow the same steps as in the proof of Theorem \ref{theorem2}. \QEDB
\end{remark}

\vspace{0.1cm}
We finish this section by noting that (with some moderate effort) the stability results of Theorems \ref{theorem1}-\ref{theorem2} can be extended to systems for which Lyapunov functions with monomial bounds of the form \eqref{quadraticbounds0} do not exist. However, while describing an interesting research direction, such characterizations are out of the scope of this paper and could be studied in the future in the context of integral-ISS, as described in \cite{Liu_Tanwani_Liberzon2022}. Indeed, for our applications of interest, described in the next section, as well as others not discussed here due to space constraints (e.g., concurrent learning \cite{CL1}, extremum seeking \cite{book_ESC,zero_order_poveda_Lina,todorovski2023practical}, feedback-optimization, etc), Assumption \ref{assumptionunstablemodes} is usually satisfied.  
%
\section{APPLICATIONS TO PT-CONTROL AND PT-DECISION MAKING}\label{sectionapplications}
In this section, we present three applications that illustrate our main results, and which rely on the mathematical models introduced in Sections \ref{analyticalresults}-\ref{sectionISSTheorems}. Thus, when we mention the state $q$ or the blow-up gain $\mu_k$, we assume they follow the hybrid dynamics \eqref{mainHDSmodel} with data \eqref{mapsstable01} or \eqref{hdsUnstable}.
\subsection{PT-Regulation of Switching Plants via Feedback Linearization}
\begin{figure}[t]
    \centering
    \includegraphics[width=0.5\linewidth]{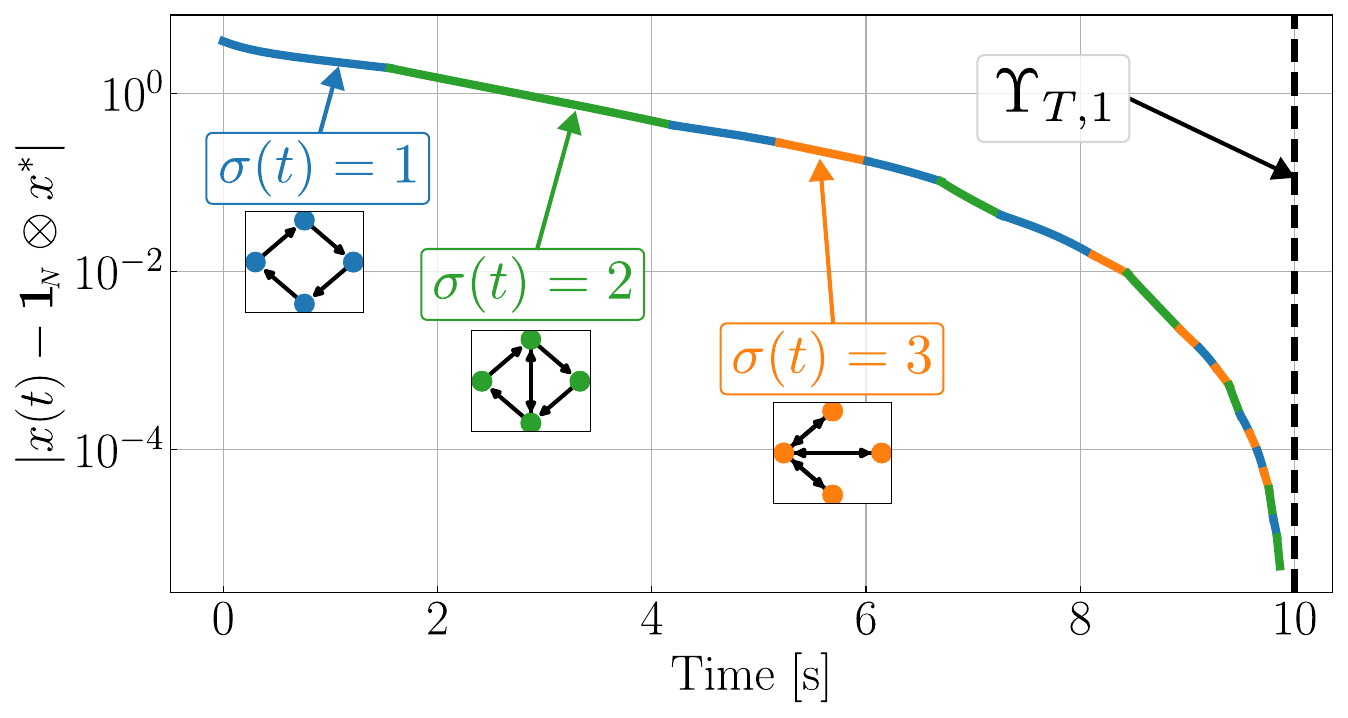}
    \caption{Norm of the error between the state of a group of $N=4$ agents and the consensus location $x^*=(-1,1)$. The error converges to $0$ as $t\to \Upsilon_{T,1}=10$.
    }\vspace{-0.3cm}
    \label{fig:switchingGraphv2Norm}
\end{figure}

Consider a switched input-affine system of the form:
\begin{equation}\label{regulation:InputAffineSwitchedSystem}
\dot{x}=d_{q}(x) + b_{q}(x)u_{q}(x,\mu_k),
\end{equation}
where  $x\in\mathbb{R}^n$, $q\in\mathcal{Q}$, $\mu_k$ is as defined in \eqref{tpmu}, $d_q(x)\in \mathbb{R}^n$, $b_q(x)\in\mathbb{R}^{n\times n}$ is invertible for all $q\in\mathcal{Q}$ and all $x\in\mathbb{R}^n$, and $u_{q}:\mathbb{R}^n\times\mathbb{R}_{\ge1} \to \mathbb{R}^n$ is the control input. Assume that $b_q(\cdot)$ and $d_q(\cdot)$ are known for all $q\in\mathcal{Q}$ (this assumption is relaxed in the next section), and consider the switched feedback law:
\begin{equation}\label{regulation:feedback}
u_q(x,\mu_k)\coloneqq b_q(x)^{-1}\left( \mu_k A_q x - d_q(x)\right).
\end{equation}
The closed-loop system has the form of the HDS \eqref{mainHDSmodel} with data \eqref{mapsstable01}, and leads to the following proposition, which extends the results of \cite[Sec. 3]{song2017time} for input-affine systems without disturbances to the case where the system switches between a finite number of stable modes. The proof follows by direct substitution and application of Theorem \ref{theorem1}, and for completeness, it is presented in the Supplemental Material.
\vspace{0.05cm}
\begin{prop}\label{Prop:Regulation_Switching_Plants}
Suppose that $A_{q}$ is Hurwitz for all $q\in\mathcal{Q}$, and let $\mathcal{A}_1$ be given by \eqref{stableset1}. Then, there exists $\tau_d>0$ such that the closed-loop system \eqref{mainHDSmodel} with data \eqref{mapsstable01} renders the set $\mathcal{A}_1\times\mathbb{R}_{\geq1}$ PT-S$_\text{F}$.\QEDB
\end{prop}
\begin{figure}[t!]
    \centering
    \includegraphics[width=0.95\linewidth]{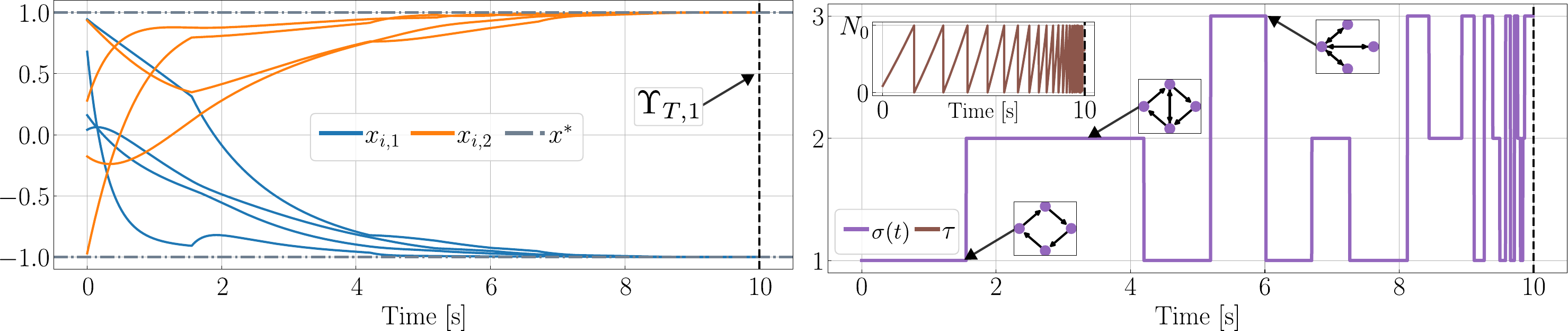}
    \caption{PT-Regulation of a group of $N=4$ agents with switching input-affine dynamics to the consensus location $x^*=(-1,1)$ using ``deficient'' information and switching communication digraphs. Top: System's state trajectory. Bottom: Trajectories of the switching signal $\sigma(t)$ (color purple), and dwell-time state $\tau$ (color brown). }
    \label{fig:switchingGraphv2}
\end{figure}
\begin{example}
   To illustrate Proposition \ref{Prop:Regulation_Switching_Plants}, we consider a group of $N=4$ agents, each with an internal state $x_i\in\mathbb{R}^n$, where $n=2$. The state of the entire group is denoted by $x\coloneqq (x_1,x_2,x_3,x_4)\in\mathbb{R}^{N\cdot n}$. The dynamics of the agents are described by \eqref{regulation:InputAffineSwitchedSystem} with:
    \begin{align*}
        d_{q}(x) = q \tanh(x), \quad \text{and} \quad b_{q}(x) = I, \quad \forall q\in\mathcal{Q},~x\in\mathbb{R}^{N\cdot n},
    \end{align*}
    where the function $\tanh(\cdot)$  operates component-wise on $x\in\mathbb{R}^{N\cdot n}$, and $\mathcal{Q}=\{1,2,3\}$.
    The agents aim to collectively converge to a common location $x^*=(-1,1)$, i.e., $x_i\to x^*$ for every $i\in\mathcal{V}=\{1,2,3,4\}$, and as $t\to\Upsilon_{T,k}$. Each agent has access to measurements of the error signal:
     \begin{align*}
         \eta_i(x_i) = B_i(x_i-x^*),
     \end{align*}
     where $B_i\in\mathbb{R}^{n}$, and $B_i\succeq 0$. We assume that this error function is ``deficient'' by letting $\text{ker}(B_i)\neq {0}$ for all $i$. To overcome this limitation, the agents communicate their current state $x_i$ through a switching digraph $\mathcal{G}_q = (\mathcal{V},\mathcal{E}_q)$, where $\mathcal{E}_q$ represents a set of switching edges. Assuming that there exists $\alpha_q>0$ such that $\mathcal{L}_q + B\succ \alpha_q I$, where $B\coloneqq \text{diag}\left(\{B_i\}_{i=1}^3\right)$, $\mathcal{L}_{q}\coloneqq\mathbb{L}_q\otimes I_n$, and  $\mathcal{L}_q$ is the Laplacian matrix of the digraph $\mathcal{G}_q$, the $i^{th}$ agent implements the control law:
     \begin{align*}
         u_{i,q}(x,\mu_1) &=\!-\mu_1\!\left(k_{\text{r}}\cdot\eta_i(x_i) + k_{\text{c}}\cdot\!\!\!\!\!\sum_{j\in\mathcal{\mathcal{N}}^{\text{out}}(i)}\!\!\!\!a_{ij,q}(x_i{-} x_j)\right)\\
         &\qquad {-} q\tanh(x_i).
     \end{align*}
     where $a_{ij,q}$ are the entries of the adjacency matrix of the digraph, and $k_{\text{r}},~k_{\text{c}}>0$ are tunable gains. The control law can be written in vector form as:
     \begin{align}\label{application1:controlLaw}
         u_{q}(x,t) &= \mu_1A_{q}(x-\mathbf{1}_N\otimes x^*) - q\tanh(x),
    \end{align}
    where $A_{q}\coloneqq -(B+L_q)$. This control law is a realization of \eqref{regulation:feedback}. Since the positive definiteness of $\mathcal{L}_q + B$ implies that $A_{q}$ is Hurwitz for all $q\in\mathcal{Q}$, all the conditions of Proposition \ref{Prop:Regulation_Switching_Plants} are satisfied. Accordingly, to numerically verify the PT-S$_{\text{F}}$ property, we simulate the closed-loop dynamics using a switching signal $\sigma\in \Sigma_{\text{BU-ADT}}(\tau_d,N_0,T,\mu_0)$ with $\tau_d=0.3129$, $T=10$, and $N_0=3$. Figure \ref{fig:switchingGraphv2} shows the trajectories of the states $q$ and $\tau$ in \eqref{mapsstable01}, associated with the switching signal $\sigma(t)$, as well as the trajectories of the agents' states $x$. As expected, Figure \ref{fig:switchingGraphv2Norm} reveals that the norm of the error with respect to the consensus location $x^*$, plotted on a logarithmic scale, rapidly approaches zero as $t\to\Upsilon_{T,1}$. \QEDB 
\end{example}
\subsection{PT-Regulation with Intermittent Feedback}
We now relax the assumptions of the previous section by considering two modifications: First, we introduce intermittent feedback, which captures scenarios where the control input is not able to affect the dynamics of the system. Second, we assume that $b_q(\cdot)$ and $d_q(\cdot)$ are unknown. Formally, we now consider systems of the form
\begin{align}
    \dot{x} = d_{q}(x) + \mathbb{I}_{\mathcal{Q}_s}\left(q\right)b_{q}(x)u_{q}(x,\mu_k),\label{application2:dynamics}
\end{align}
where $q\in\mathcal{Q}=\mathcal{Q}_s\cup\mathcal{Q}_u$ is a logic state and $\mathcal{Q}_u\neq \emptyset$. We assume that $b_q(\cdot)$ and $d_q(\cdot)$ are unknown, but locally Lipschitz, and satisfy the following ``matching'' and positive definiteness conditions:
\begin{align*}
&|d_q(x)|\leq \overline{d}_q(x), \quad ~~~~~~\forall~ q\in\mathcal{Q},~x\in\mathbb{R}^n,\\
&b_q(x) + b_q(x)^\top \succeq \epsilon I_n,~~\forall~q\in \mathcal{Q}_s,~x\in\mathbb{R}^n,
\end{align*}
where $\epsilon>0$, and $\overline{d}_q(x)>0$ is a known scalar-valued function assumed to be continuous for all $q\in\mathcal{Q}$. We also assume that $\overline{d}_q(x)$ is $\ell_q$-globally Lipschitz for all $q\in\mathcal{Q}_u$. To regulate the state $x$ to the origin in a prescribed time, we consider the following switching feedback-law:
\begin{align}
u_q(x, \mu_k) &= -\mu_k \left(\eta_q + \delta_q\overline{d}_q(x)^2\right)x,\label{feedback}
\end{align}
with $\delta_q>0$ and $\eta_q>0$. The closed-loop system has the form of the HDS \eqref{mainHDSmodel} with the data \eqref{hdsUnstable}, and leads to the following Proposition by a direct application of Theorem \ref{theorem2}. For completeness, the proof and computations are presented in the Supplemental Material.
\vspace{0.1cm}
\begin{prop}\label{Prop:Regulation_intermittenUpsilonlants}
    There exists $\tau_d>0$ and $\tau_a>0$ such that the closed-loop system renders the set $\mathcal{A}_1\times\mathbb{R}_{\ge1}$ PT-ISS-C$_{\text{F}}$, where $\mathcal{A}_1$ is as given in \eqref{stableset2}.~\QEDB
\end{prop}

\begin{figure}[t]
    \centering
    \includegraphics[width=0.65\linewidth]{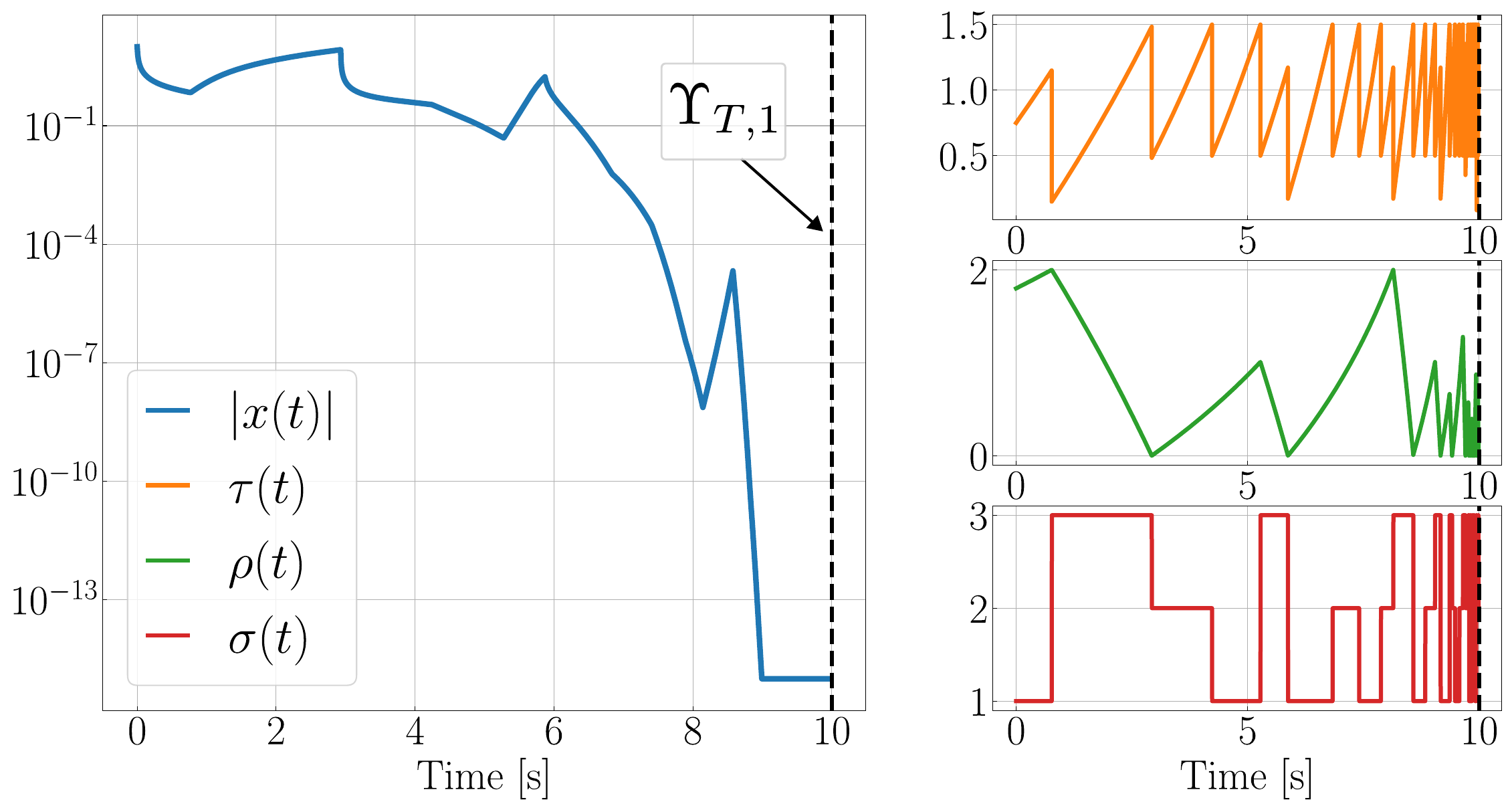}
    \caption{PT-Regulation with intermittent feedback. Left: System's state trajectory. Right: Trajectories of the dwell-time state $\tau$ (color orange), the monitor state $\rho$ (color green), and switching signal $\sigma(t)$ (color red). }
\label{application2:figureSims}\vspace{-0.4cm}
\end{figure}

\vspace{0.1cm}
\begin{example}
    To illustrate Proposition \ref{Prop:Regulation_intermittenUpsilonlants} with a numerical example, consider $\mathcal{Q}_s = \{1,2\}$, $\mathcal{Q}_u =\{3\}$, and $x\in\mathbb{R}$. Let
     $ d_{q}(x) = q\tanh(x),~ b_{q}(x)=1,~\forall q\in\mathcal{Q}$. After choosing the control-law $u_q(x,t)=-\mu_1(t) (1 + q|x|^2)x$, all the conditions to apply Proposition \ref{Prop:Regulation_intermittenUpsilonlants} are satisfied. We numerically verify the PT-ISS-C${\text{F}}$ property by using a switching signal $\sigma\in\Sigma_{\text{BU-ADT}}(\tau_d,N_0,T,\mu_0)\cap \Sigma_{\text{BU-AAT}}(\mathcal{Q}_u,\tau_a,T_0,T,\mu_0)$ with $\tau_a=2$, $\tau_d=1$, $T=10$, $T_0=2$, and $N_0=1.5$. Figure \ref{application2:figureSims} displays the trajectories of the state $x$, the switching signal $\sigma$, and the associated average dwell-time and average activation time states $\tau$ and $\rho$. As shown in the figure, the state $x$ rapidly approaches zero as $t\to\Upsilon_{T,1}$. The overshoots occur when the system is in one of the unstable modes $(q\in\mathcal{Q}_u)$. %
     \QEDB 
\end{example}
%
%
%
%
\subsection{PT-Decision-Making in Switching Games}
\label{resetsexample}
We consider a non-cooperative game with $n\in\mathbb{Z}_{\geq2}$ players \cite{poveda2022fixed}, where the cost functions defining the game switch in time. Formally, for each $i\in\mathcal{V}=\{1,2,\cdots, n\}$, we consider that the $i^{th}$ player has an associated mode-dependent and continuously differentiable cost function $\phi_{q}^i:\mathbb{R}^n\to\mathbb{R}$, where $q\in\mathcal{Q}$. We refer to the $q^{th}$ game as the game with the set of cost functions $\left\{\phi_{q}^i\right\}_{i\in\mathcal{V}}$. The action of the $i^{th}$ player is denoted by $x_1^{i}\in\mathbb{R}$, and the action profile of the game is given by the vector $x_1 \coloneqq \left(x_1^1, x_1^2, \ldots, x_1^n\right) \in \mathbb{R}^n$. The goal of the players is to converge to the unique common Nash equilibrium (NE) of the games, defined as the vector $\tilde{x}_{1}\in\mathbb{R}^n$ that satisfies:
\begin{align*}
\phi_{q}^i\left(\tilde{x}_{1}^{i},~\tilde{x}_{1}^{-i}\right) = \inf_{x^{i}\in\mathbb{R}}\phi^i_q\left(x_{1}^{i},~\tilde{x}_{1}^{-i}\right),~~\forall i\in\mathcal{V},
\end{align*}
for each $q\in\mathcal{Q}$, where $x_{1}^{-i}\in\mathbb{R}^{n-1}$ denotes the vector that contains all actions except those of player $i$. To study this problem, let $\mathcal{G}_q:\mathbb{R}^n\to\mathbb{R}^n$ denote the pseudo-gradient of the $q^{th}$ game, which is given by: 
\begin{equation*}
    \mathcal{G}_{q}(x_1) \coloneqq \left(\frac{\partial\phi_{q}^1}{\partial x_1^1},~\frac{\partial\phi_{q}^2}{\partial x_1^2},\ldots,\frac{\partial\phi_{q}^n}{\partial  x_1^n},\right).
\end{equation*}
 For all $q\in\mathcal{Q}$, we assume that there exists $\kappa_q>0$ and $\ell_q>0$ such that $\mathcal{G}_q$ is a $\kappa_q$-strongly monotone and $\ell_q$-globally Lipschitz mapping. These properties guarantee the existence and uniqueness of the NE $\tilde{x}_1$ \cite{poveda2022fixed}. %
To efficiently achieve convergence to the NE in a prescribed time, we introduce \emph{PT high-order NE-seeking dynamics with momentum and resets} (PT-NESmr), given by the system \eqref{mainHDSmodel} with data \eqref{mapsstable01} and maps $f_q$ and $R_q$ defined as follows:
\begin{align}\label{pT:NESwMomentum}
     f_q(x,\tau) = \begin{pmatrix}
        \dfrac{2}{\eta(\tau)}\left(x_2 {-} x_1\right)\vspace{0.1cm}\\
        -2\eta(\tau)\mathcal{G}_{q}(x_1)
    \end{pmatrix},~~R_q(x) = \begin{pmatrix}
        x_1\vspace{0.1cm}\\x_1
    \end{pmatrix},
\end{align}
where $x\coloneqq(x_1,x_2)\in\mathbb{R}^{2n}$, and $x_2\coloneqq (x_2^1, x_2^2, \ldots, x_2^n)\in\mathbb{R}^n$, and where $\eta:[0,N_0]\to[\underline{\eta},\overline{\eta}]$ is an affine bounded mapping defined as:
\begin{equation}
\eta(\tau) \coloneqq \tau\frac{\left(\overline{\eta} - \underline{\eta}\right)}{N_0}+ \underline{\eta}
\end{equation}
with $\overline{\eta}>\underline{\eta}>0$ being tunable parameters. In the context of asymptotic convergence, mappings of the form \eqref{pT:NESwMomentum}, which incorporate momentum (via the state $x_2$) and resets (via the update $x_2^+=x_1$), have been recently shown to improve the transient performance of NE-seeking dynamics in strongly monotone games  \cite{ochoa2021momentum}. 

To study convergence to the NE in prescribed time, we let\footnote{Alternatively, one could use a change of variables to shift the NE to the origin so that $\mathcal{A}_x=\{0\}\in\mathbb{R}^{2n}$ as in \eqref{stableset1}.}:
\begin{align}\label{application3:set}
\mathcal{A}_x \coloneqq \{\tilde{x}_1\}\times\{\tilde{x}_1\}\in\mathbb{R}^{n}\times\mathbb{R}^n,
\end{align}
and we make use of the following assumptions:

\vspace{0.1cm}
\begin{assumption}[Game-Jacobian Regularity Condition]\label{PTNES:gameRegularity}
      For every $q\in\mathcal{Q}$, there exists $\sigma_q\in \mathbb{R}_{>0}$ such that $\sigma_{\max}\left(I-\partial \mathcal{G}_q(x)\right) \le \sigma_q$, for all $x\in\mathbb{R}^n$, where $\sigma_{\max}(\cdot)$ is the maximum singular value of its argument. \QEDB
\end{assumption}

\vspace{0.1cm}
\begin{assumption}[Tuning Guidelines]\label{PTNES:tuningGuidelines}
   Let $\zeta_{q}\coloneqq \kappa_q/\ell_q^2$, and consider the choice of parameter  
   \begin{equation*}
   \overline{\eta}^2\le  \delta_{\eta}\frac{\min_{q\in\mathcal{Q}}\zeta_q}{\left(\max_{q\in\mathcal{Q}}\sigma_q\right)^2},~~~~\frac{1}{\tau_d} \le \delta_d\frac{N_0}{\overline{\eta}-\underline{\eta}} \min_{q\in\mathcal{Q}}\zeta_q
   \end{equation*}
   where $\delta_\eta\in(0,1)$, $\delta_d>0$, and $\delta_\eta + \delta_d\eqqcolon \delta \in(0,1)$. Moreover, suppose that $\overline{\gamma} \in (0,1],$ where
    $
            \overline{\gamma} \coloneqq \frac{\overline{\ell}^2}{\underline{\kappa}^2}\frac{\eta(N_0-1)^2}{\eta(1)^2} {+} \frac{1}{2\underline{\kappa}^2 \eta(1)^2}.$
    \QEDB
\end{assumption}

\vspace{0.15cm}
The following Proposition follows by a direct application of Theorem \ref{theorem1}, and a standard Lyapunov-based argument that shows that, under Assumptions \ref{PTNES:gameRegularity}-\ref{PTNES:tuningGuidelines}, the target HDS satisfies the conditions of Assumption \ref{assumptionunstablemodes}. For the purpose of completeness, the step-by-step proof and computations are presented in the Supplemental Material.

\vspace{0.1cm}
\begin{prop}\label{Prop:Nash_equilibrium} 
    Suppose that Assumptions \ref{PTNES:gameRegularity}-\ref{PTNES:tuningGuidelines} are satisfied. Then, the PT-NESmr dynamics render the set $\mathcal{A}_{x}\times [0,N_0]\times \mathcal{Q}\times \mathbb{R}_{\ge 1}$ PT-S$_{\text{F}}$, whenever
    \begin{equation*}
        \tau_d>\frac{\max\left\{3, 2\left(\frac{1}{\underline{\kappa}^2} + \overline{\eta}^2\right)\right\}}{4\underline{\eta}\nu}\ln\left(\frac{\max\{3,~2+2\overline{\ell}^2\overline{\eta}^2\}}{\min\{1, 2\underline{\eta}\underline{\kappa}^2\}}\right),
    \end{equation*}
    where $\nu= \frac{(1-\delta_d - \delta_\eta)\overline{\sigma}^2}{\delta_{\eta}(1-\delta_d)\underline{\zeta} + \overline{\sigma}^2}$, where $\underline{\zeta}\coloneqq \min_{q\in\mathcal{Q}}\zeta_q$\QEDB
\end{prop}

\vspace{0.2cm}
\noindent The following numerical example illustrates Proposition \ref{Prop:Nash_equilibrium}.
\begin{figure}[t]
    \centering
    \includegraphics[width=0.75\linewidth]{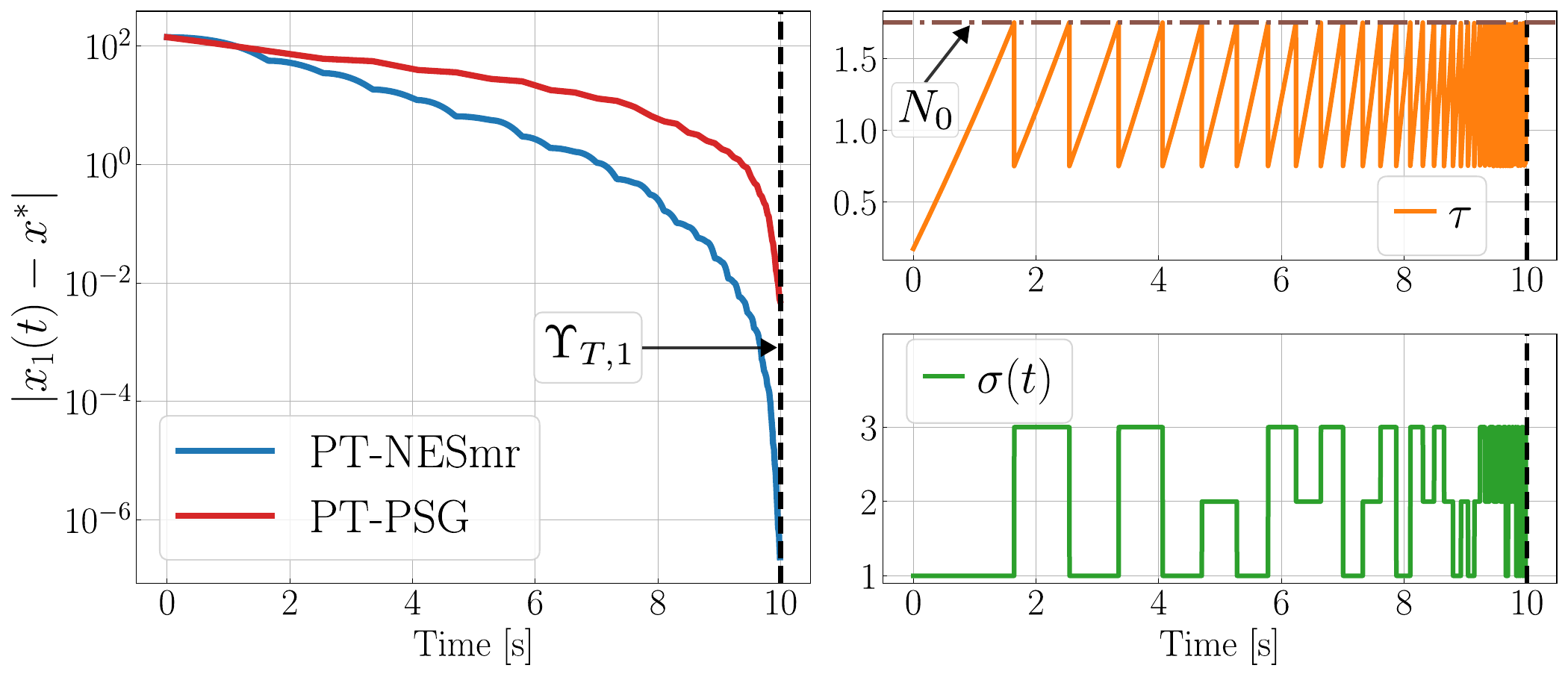}
    \caption{PT Nash-Equilibrium Seeking in a Strongly-Monotone Switching Game. Left: Trajectory of the errors to the NE generated by the PT-NESmr and PT-PSG dynamics. Right: Trajectories of the dwell-time state $\tau$ (color orange) and the switching signal $\sigma(t)$ (color green).}
    \label{application3:figureSims}\vspace{-0.5cm}
\end{figure}
\vspace{0.2cm}
\begin{example}
    Let $\mathcal{Q} = \{1,2,3\}$ and $\mathcal{G}_q(x_1) = \vartheta A_q(x_1-\tilde{x}_1)$, with $\tilde{x}_1=(1,1)$, $A_1 = [6,-1.5; -1.5, 6]$, $A_2 = [8, -2; 2, 8]$, $A_3=[4,0;0;8]$, and  $\vartheta=5\times 10^{-2}$.
    Note that the pseudo-gradients $\mathcal{G}_q(\cdot)$ are $\kappa_q$-strongly monotone and $\ell_q$-globally Lipschitz. Using $k=1$, $\tau_d = 1.14$, $N_0=1.75$, and $\Upsilon_{T,k}=10$ all the conditions to apply Proposition \ref{Prop:Nash_equilibrium} are satisfied. We simulate the system using a switching signal $\sigma\!\in\!\Sigma_{\text{BU-ADT}}(\tau_d,N_0,T,\mu_0)$ with $T=10$. We compare our results with the continuous-time prescribed-time pseudo-gradient-flows (PT-PSG), recently introduced in \cite{ZhaoTaoprescribedTimeNash2023}, and given by $\dot{x}_1=\mu_1(t)\mathcal{G}_{\sigma(t)}(x_1)$. The resulting trajectories are shown in Figure \ref{application3:figureSims}. As observed the synergistic incorporation of momentum, resets, and PT techniques leads to an improvement compared to the continuous-time PT algorithm, under the same switching signal. \QEDB  %
\end{example}
%

%

\section{CONCLUSIONS}
\label{secconclusions}
The property of prescribed-time stability was studied for a class of hybrid dynamical systems incorporating switching nonlinear vector fields with time-varying increasing gains, exogenous inputs, and resets. Switching conditions that preserve the stability of the system were derived using tools from hybrid dynamical systems theory and under a suitable contraction/dilation of the hybrid time domains. The switching conditions permit the incorporation of unstable modes. The results were illustrated in three different applications in the context of control and decision-making. Future applications will include prescribed-time concurrent learning and prescribed-time switching extremum seeking where the assumptions of this paper are satisfied. Additional potential extensions might also include switching ODEs interconnected with linear hyperbolic PDEs.
%
\bibliographystyle{ieeetr}
\bibliography{REFERENCES}

\vspace{-0.05cm}
\clearpage
\section{APPENDIX}\label{sec:appendix}
\subsection{Proof of Proposition \ref{transformationk}:}

\vspace{0.05cm}
\noindent (P1) Follows by the monotonicity of $\omega_k(\cdot, \cdot)$ in its first argument, combined with the limit $\lim_{t\to\Upsilon_{T,k}}\mu_k(t)=\infty$.

\vspace{0.1cm}\noindent (P2) For $k>1$, the result follows by direct computation. For $k=1$, the result is obtained by the properties of the logarithm.

\vspace{0.1cm}\noindent (P3) By definition, the equality  $\mathcal{T}_k(0)=0$ holds for all $k\in\mathbb{R}_{\ge 1}$. For $k=1$, by direct computation, we have: $\frac{d\mathcal{T}_1(t)}{dt} = \frac{T}{\mu_{1}(t)}\dot{\mu}_1(t) = \mu_1(t).$ For $k>1$, by the chain rule, we obtain:
$$
    \frac{d\mathcal{T}_k(t)}{dt} = \frac{\partial\omega_k(b,\mu_k(0))}{\partial b}\big|_{b=\mu_k(t)}\dot{\mu}_k
                                 = \mu_k(t).$$     
\vspace{0.1cm}\noindent (P4) For $k=1$, we have that  $\mu_1(t) = \frac{\mu_0 T}{T-\mu_0t}$. It then follows that $s=\left(\mathcal{T}_1\circ\mathcal{T}_1^{-1}\right)(s) = T\ln\left(\frac{\mu_{1}\left(\mathcal{T}^{-1}(s)\right)}{\mu_0}\right)$. Solving for $\mathcal{T}_1^{-1}(s)$ leads to $\mathcal{T}_1^{-1}(s) = \Upsilon_{T,1}\left(1 - e^{-\frac{s}{T}}\right)$.
For $k>1$, let $y_k\coloneqq \mathcal{T}_k^{-1}$. By using \eqref{tpmu}, and the inverse function theorem, we obtain that$\frac{dy_k}{ds}=\frac{\left(\Upsilon_{T,k}-y\right)^k}{T^k}$. 
Then, by direct integration and using the fact that $y_k(0)=0$, we obtain the following equality $\Upsilon_{T,k}-y_k(s)=\left(\frac{(k-1)s}{T^k}+\Upsilon_{T,k}^{1-k}\right)^{\frac{1}{1-k}}.$
Solving for $\mathcal{T}_k^{-1}(s)$, we obtain that
$\mathcal{T}^{-1}_k(s)=\Upsilon_{T,k}-\Upsilon_{T,k}\left(1+\frac{(k-1)s}{\Upsilon_{T,k}\mu_0}\right)^{\frac{1}{1-k}}.$

\vspace{0.1cm
}\noindent (P5) Follows directly by the inverse function theorem.

\vspace{0.1cm}
\noindent (P6) For $k=1$, using the equality $\ln(1-x)=\sum_{l=1}^\infty \frac{-1}{l}x^l$, $|x|<1$, we obtain that $\mathcal{T}_1(t) =\mu_0t + \sum_{l=2}^\infty\frac{1}{l}\mu_0^lt^lT^{1-l}$, for all $t\mu_{0}< T$.
Letting $T\to\infty$, the second term in this expression vanishes, and we obtain that the equality $\lim_{T\to\infty}\mathcal{T}_{1}(t) = \mu_0t$ holds for all $(t,\mu_0)\in\mathbb{R}_{\ge 0}\times \mathbb{R}_{\ge 1}$. For $k>1$, from Remark \ref{remark:equivalentTransform} it follows that
\begin{align}
    \mathcal{T}_k(t) &= \frac{T\mu_0^{\frac{k-1}{k}}}{k-1}\left(\left(1-\frac{t\mu_0^{\frac{1}{k}}}{T}\right)^{1-k}-1\right).\label{proof:P6:a}
\end{align}
Now, using the binomial theorem we have that
\begin{align*}
   \left(1-\frac{t\mu_0^{\frac{1}{k}}}{T}\right)^{1-k}-1 = \frac{(k-1)t\mu_{0}^{\frac{1}{k}}}{T} + \sum_{l=2}^{\infty}g_{k,l}\left(\frac{t\mu_0^{\frac{1}{k}}}{T}\right)^{l},
\end{align*}
for all $t\mu_{0}^{\frac{1}{k}}< T$, and where $g_{k,l}=\frac{(k-1)k(k + 1)\cdots (k + l - 2)}{l!}$. Thus, for all $t\mu_{0}^{\frac{1}{k}}< T$, equality \eqref{proof:P6:a} can be written as
$
     \mathcal{T}_k(t) = \mu_0^{\frac{k-1}{k^2}}t + \sum_{l=2}^\infty \frac{g_{k,l}}{k-1}t^{l}\mu_0^{\frac{(k-1)l}{k^2}} T^{1-l}.
$
Letting $T\to\infty$, the second term in this expression vanishes. Thus, it follows that the limit $\lim_{T\to\infty}\mathcal{T}_{k}(t) = \mu_0^{\frac{k-1}{k^2}}t$ holds for all $(t,\mu_0)\in\mathbb{R}_{\ge 0}\times \mathbb{R}_{\ge 1}$.
\hfill $\blacksquare$

%
\subsection{Supplemental Material}
We present detailed proofs of all the auxiliary lemmas and propositions used in the paper. These results follow directly by computations and/or straightforward extensions or specializations of existing results in the literature. 

\vspace{0.2cm}
\noindent \textbf{Proof of Lemma \ref{lemmaodek}}:

\noindent First, we have that:
\begin{align*}
\int_{\mu_0}^{\mu_{k}(t)}\frac{d\mu_k}{\mu_k^{1+\frac{1}{k}}}=\int_{0}^t\frac{k}{T}dt\implies
&-k\mu_k^{\frac{-1}{k}}\Big|_{\mu_0}^{\mu_k(t)}=\frac{k}{T}t.
\end{align*}
Thus, it follows that $k\left(-\mu_k(t)^{\frac{-1}{k}}+\mu_k(0)^{\frac{-1}{k}}\right)=\frac{k}{T}t$, and:
\begin{align*}
\frac{1}{\mu_k(t)^\frac{1}{k}}=\frac{1}{\mu_k(0)^{\frac{1}{k}}}-\frac{t}{T}=\frac{T-t\mu_k(0)^\frac{1}{k}}{T\mu_k(0)^{\frac{1}{k}}},
\end{align*}
from which we obtain the result. \hfill $\blacksquare$

\vspace{0.2cm}
\noindent\textbf{Proof of Lemma \ref{completemulemma}:} 

\noindent We have that:
\begin{align*}
\int_{\mu_0}^{\hat{\mu}_k(t)}\frac{d\hat{\mu}_k}{\hat{\mu}_k^{\frac{1}{k}}}=\int_{0}^t\frac{k}{T}dt\implies &\frac{1}{1-\frac{1}{k}}\hat{\mu}_k^{1-\frac{1}.{k}}\Big|_{\mu_0}^{\hat{\mu}_k(t)}=\frac{k}{T}t.
\end{align*}
Therefore, we obtain $\frac{k}{k-1}\left(\hat{\mu}_{k}^{1-\frac{1}{k}}(t) - \mu_0^{1-\frac{1}{k}} \right)=\frac{k}{T}t$, and:
\begin{align*}
\hat{\mu}_{k}(t) = \left(\frac{k-1}{T}t + \mu_0^{\frac{k-1}{k}}\right)^{\frac{k}{k-1}}.
\end{align*}
This obtains the result. \hfill $\blacksquare$

\vspace{0.2cm}

\vspace{0.1cm}
\noindent\textbf{Proof of Lemma \ref{equivalentBUADTk}:}

\noindent The case $k=1$ follows directly by the definition of $\mathcal{T}_1$ and Remark \ref{casek1}. For $k>1$, consider expanding the right-hand side of \eqref{BUADTk}:

\begin{small}
\begin{align*}
&N(t_2,t_1)\leq \frac{T}{\tau_d}\left(\frac{\mu_k(t_2)^\frac{k-1}{k}}{k-1}-\frac{\mu_k(t_1)^\frac{k-1}{k}}{k-1}\right)+N_0\\
&\qquad=\frac{T^k}{\tau_d(k-1)}\left(\frac{\left(\Upsilon_{T,k}-t_1\right)^{k-1}-\left(\Upsilon_{T,k}-t_2\right)^{k-1}}{\left(\left(\Upsilon_{T,k}-t_2\right)\left(\Upsilon_{T,k}-t_1\right)\right)^{k-1}}\right)+N_0.
\end{align*}
\end{small}
Taking the limit as $k\to1$, one obtains \eqref{BUKADT111}, see also Remark \ref{casek1}. On the other hand, when $k\in\mathbb{R}_{>1}$, the Binomial theorem can be use to write $(\Upsilon_{T,k}-t_i)^{k-1}=\sum_{\ell=0}^{k-1}b_{k,l}\Upsilon_{T,k}^{k-1-\ell}(-t_i)^\ell$, for $i\in\{1,2\}$, where $b_{k,l}:=\frac{(k-1)!}{\ell!(k-\ell-1)!}$ are the so-called Binomial coefficients. Let
\begin{align*}
S&\coloneqq\sum_{\ell=0}^{k-1}b_{k,l}\Upsilon_{T,k}^{k-1-\ell}(-t_1)^\ell-\sum_{\ell=0}^{k-1}b_{k,l}\Upsilon_{T,k}^{k-1-\ell}(-t_2)^\ell\\
&=b_{k,1}\Upsilon_{T,k}^{k-2}(t_2-t_1)+\sum_{\ell=2}^{k-1}(-1)^{\ell+1}b_{k,l}\Upsilon_{T,k}^{k-1-\ell}\left(t_2^\ell-t_1^\ell\right).
\end{align*}
Therefore, the BU$_k$-ADT bound can be written as
\begin{align*}
N(t_2,t_1)&\leq \frac{T^k}{\tau_d(k-1)}\left(\frac{S}{\left((\Upsilon_{T,k}-t_2)(\Upsilon_{T,k}-t_1)\right)^{k-1}}\right)+N_0\\
&=\gamma_k(t_1,t_2)\left[(t_2-t_1)+\sum_{\ell=2}^k\tilde{c}_{\ell,k}\left(t_2^\ell-t_1^\ell\right)\right]+N_0,
\end{align*}
where
\begin{align*}
\gamma_k(t_1,t_2)%
&=\frac{b_{k,1}\mu_0^{\frac{2-k}{k}}}{\tau_d(k-1)}\left(\frac{T^2}{\left(\Upsilon_{T,k}-t_2\right)\left(\Upsilon_{T,k}-t_1\right)}\right)^{k-1},%
\end{align*}
and
\begin{align*}
\tilde{c}_{\ell,k}&=(-1)^{\ell+1}b_{k,l}\Upsilon_{T,k}^{k-1-\ell} \left(b_{k,1}\Upsilon_{T,k}^{k-2}\right)^{-1}=(-1)^{\ell+1}\frac{b_{k,l}}{b_{k,1}}\Upsilon_{T,k}^{1-\ell}.
\end{align*}
The result follows by using $b_{k,1}=k-1$.  \hfill $\blacksquare$
%

\vspace{0.4cm}
\begin{lemma}
Suppose that Assumption \ref{assumptionunstablemodes} holds. Then, there exists $\omega\geq1$ such that:
\begin{equation}\label{boundjumpsLyapunov}
V_{o}(x) \leq \omega V_q (x),
\end{equation}
for all $(x,o,q) \in (C\cup D)\times \mathcal{Q}\times \mathcal{Q}$. \QEDB
\end{lemma}

\vspace{0.1cm}
\noindent\textbf{Proof:} Using item (a) in Assumption \ref{assumptionunstablemodes}, we have:
\begin{subequations}
\begin{align}
 &V_q (x)\leq \max_{q\in\mathcal{Q}} c_{q,2}|x|^p,~\forall~x\in\mathbb{R}^n,q\in\mathcal{Q}\\
 &\min_{q\in \mathcal{Q}}c_{q,1}|x|^p\leq V_o (x),~\forall~x\in\mathbb{R}^n,o\in\mathcal{Q}
\end{align}
\end{subequations}
for all $x\in\mathbb{R}^n,q\in\mathcal{Q}$. It follows that
\begin{equation*}
V_q (x)\leq \frac{\max_{q\in\mathcal{Q}} c_{q,2}}{\min_{q\in \mathcal{Q}}c_{q,1}}V_o(x), ~~\forall~x\in\mathbb{R}^n,~(q,o)\in \mathcal{Q}\times\mathcal{Q},
\end{equation*}
which implies that \eqref{boundjumpsLyapunov} holds with $\omega=\frac{\max_{q\in\mathcal{Q}} c_{q,2}}{\min_{q\in \mathcal{Q}}c_{q,1}}\geq1$. \hfill $\blacksquare$

\vspace{0.2cm}
The following lemma is a specialization of \cite[Prop. 2.7]{Cai2009} for the case when the system is exponentially ISS. We present the proof for the purpose of completeness.

\vspace{0.2cm}
\begin{lemma}\label{auxlemma01}
Consider the HDS \eqref{HDS0}, and suppose there exist constants $\underline{\alpha},\overline{\alpha},\rho,p>0$, $\lambda\in(0,1)$, and a smooth function $V:C\cup D\to\mathbb{R}_{\geq0}$, such that the following inequalities hold:
\begin{align*}
\underline{\alpha}|z|^p_{\mathcal{A}}&\leq V(z)\leq \overline{\alpha}|z|^p_{\mathcal{A}},~~~~\forall~z\in C\cup D\cup G(D),\\
\langle\nabla V(z),F(z,u)\rangle&\leq -\lambda V(z)+\rho|u|^p,~~\forall~(z,u)\in C\times\mathbb{R}^m,\\
V(G(z))-V(z)&\leq -\lambda V(z)+\rho|u|^p,~~\forall~(z,u)\in D\times\mathbb{R}^m.
\end{align*}
Then, every solution of \eqref{HDS0} satisfies
\begin{equation}\label{lemmaISS:statement}
|z(s,j)|_{\mathcal{A}}\leq \kappa_1e^{-\kappa_2(s+j)}|z(0,0)|_{\mathcal{A}}+\kappa_3\sup_{0\leq \tau\leq s}|u(\tau)|,
\end{equation}
for all $(s,j)\in\text{dom}(z)$, and  where $\kappa_1= \left(\overline{\alpha}/\underline{\alpha}\right)^p$, $\kappa_2 = \lambda/2p$, and $\kappa_3 = \left(\frac{2\rho}{\lambda\underline{\alpha}}\right)^{1/p}$. \QEDB
\end{lemma}
\vspace{0.1cm}\noindent 
\textbf{Proof:}
We follow similar ideas as in the proof of \cite[Prop. 2.7]{Cai2009}, but considering set-valued flow and jump maps. The proof has four main steps:

\vspace{0.1cm}\noindent
\textsl{Step 1:} First, note that for all $(z,u)\in (C\cup D)\times\mathbb{R}^m$:
\begin{equation}\label{condition1}
-\lambda V(z)+ \rho|u|^p\leq -\frac{\lambda}{2}V(z),~~\text{if}~V(z)\geq \frac{2\rho}{\lambda}|u|^p.
\end{equation}
Therefore, whenever $V(z)\geq\frac{2\rho}{\lambda}|u|^p$ we have that
\begin{align*}
\langle\nabla V(z),F(z,u)\rangle&\leq -\tilde{\lambda} V(z),~\forall(z,u)\in C\times\mathbb{R}^m,\\
V(G(z))-V(z)&\leq -\tilde{\lambda} V(z),~\forall(z,u)\in D\times\mathbb{R}^m,
\end{align*}
where $\tilde{\lambda}:=\lambda/2$.

\vspace{0.1cm}\noindent
\textsl{Step 2:} For any $r\geq0$, define $\gamma_{c_4}(r,s,j)=e^{-\tilde{\lambda}(s+j)}r$. We first show that when $V(z)\geq \frac{2\rho}{\lambda}|u|^p$, the function $V$ evaluated along the solutions of \eqref{HDS0} satisfies 
\begin{equation}\label{exponentialdecayLyapunov}
V(z(s,j))\leq \gamma_{\lambda}(V(z(0,0)),s,j),~~\forall~(s,j)\in\text{dom}(z).
\end{equation}
To establish this property, note that since $V(z(\cdot,\cdot))$ is not increasing during flows and jumps, if there is $(s',j')\in\text{dom}(z)$ with $0<s'+j'<t+j$ and such that $V(z(s',j'))=0$, then we necessarily must have $V(z(\tilde{s},\tilde{j}))=0$ for all $(\tilde{s},\tilde{j})\in\text{dom}(z)$ such that $s'+j'\leq \tilde{s}+\tilde{j}\leq s+j$, and \eqref{exponentialdecayLyapunov} would hold for such times $(\tilde{s},\tilde{j})$. Suppose there is no $(s',j')\in\text{dom}(z)$ with $0<s'+j'<t+j$ such that $V(z(s',j'))=0$. For each $(s,j)\in\text{dom}(z)$, we partition the hybrid time domain of $z$ up to time $(s,j)$ as $\text{dom}(z)=\bigcup_{n=0}^j[s_n,s_{n+1}]\times\{n\}$, with $s_0=0$ and $s_{j+1}=s$. For any $n\in\{0,1,\ldots,j\}$, $V$ satisfies
\begin{equation*}
\int_{s_n}^{s_{n+1}}\frac{\dot{\overbrace{V(z(\tau,n))}}}{\tilde{\lambda} V(z(\tau,n))}d\tau\leq -\int_{s_n}^{s_{n+1}}d\tau= -(s_{n+1}-s_n).
\end{equation*}
Using the new variable $\varrho=V(z(\tau,n))$, we obtain $d\varrho=\dot{V}d\tau$ and the above integral can be written as
\begin{equation}\label{decreaseflowsproofA}
\int_{V(z(s_{n},n))}^{V(z(s_{n+1},n))}\frac{d\varrho}{\tilde{\lambda} \varrho}\leq -(s_{n+1}-s_n).
\end{equation}
Similarly, note that
\begin{align*}
\int_{V(z(s_{n+1},n))}^{V(z(s_{n+1},n+1))}\frac{d\varrho}{\tilde{\lambda} \varrho}&\leq \int_{V(z(s_{n+1},n))}^{V(z(s_{n+1},n+1))}\frac{d\varrho}{\tilde{\lambda} V(z(s_{n+1},n))}\\
&\leq -1,
\end{align*}
where the last inequality follows by the inequality $V(z(s,j+1))-V(z(s,j))\leq -\tilde{\lambda}V(z(s,j))$. Combining the above two inequalities, we obtain
\begin{align}
\int_{V(z(0,0))}^{V(z(s,j))}\frac{d\rho}{\tilde{\lambda}\varrho}&=\sum_{n=0}^j\int_{V(z(s_n,n))}^{V(z(s_{n+1},n))}\frac{d\varrho}{\tilde{\lambda}\varrho}\notag\\
&~~~+\sum_{n=1}^j\int_{V(z(s_{n+1},n))}^{V(z(s_{n+1},n+1))}\frac{d\varrho}{\tilde{\lambda} \varrho}\notag\\
&\leq -\left(\sum_{n=0}^j (s_{n+1}-s_n)+\sum_{n=1}^j1\right)\notag\\
&=-(s_{j+1}-s_0+j)=-(s+j).\label{inequaintegral}
\end{align}
Integrating the left-hand side, we obtain $\frac{1}{\tilde{\lambda}}\ln\left(\frac{V(z(s,j))}{V(z(0,0))}\right)\leq -(s+j)$, from which we directly get
\begin{equation}\label{boundStep2}
V(z(s,j))\leq V(z(0,0))e^{-\frac{\lambda}{2}(s+j)}
\end{equation}
\textsl{Step 3:} Let $(z,u)$ be a maximal solution pair of \eqref{HDS0}. Define the set
\begin{equation}
\Omega:=\left\{z\in\mathbb{R}^n:V(z)\leq \frac{2\rho}{\lambda}|u|^p_{\infty} \right\}.
\end{equation}
For each $z_0\in\mathbb{R}^n$, let
\begin{align*}
T_{z,u,z_0}&:=\sup\Big\{\tau\in\mathbb{R}_{\geq0}:z(s,j)\notin\Omega,~z(0,0)=z_0,\\
&~~~~~~~~~~~~~~~\forall~(s,j)\in\text{dom}(z),~0\leq s+j\leq \tau\Big\}.
\end{align*}
It follows that for all solutions of \eqref{HDS0} with $z(0,0)=z_0$ and $(s,j)\in\text{dom}(z)$ such that $0\leq s+j<T_{z,u,z_0}$ we have that $V(z)>\frac{2\rho}{\lambda}|u|^p_{\infty}$, which, by Step 2, implies that $V$ satisfies \eqref{boundStep2}. Using the quadratic upper and lower bounds on $V$, we obtain:
\begin{equation}\label{exponentialdecay}
|z(s,j)|_{\mathcal{A}}\leq \left(\frac{\overline{\alpha}}{\underline{\alpha}}\right)^{\frac{1}{p}}|z(0,0)|_{\mathcal{A}}e^{-\frac{\lambda}{2p}(s+j)},
\end{equation}
which holds for all $(s,j)\in\text{dom}(z)$ such that $0\leq s+j<T_{z,u,z_0}$.

\vspace{0.1cm}\noindent 
\textsl{Step 4:} The last step is to prove forward invariance of $\Omega$.
Suppose there exist $(s',j')\in\text{dom}(z)$ such that $z(s',j')\in \Omega$ and $(s',j'+1)\in\text{dom}(z)$. Since $\tilde{\lambda}<\lambda$, $V$ satisfies
\begin{align*}
V(z(s',j'+1))&\leq (1-\tilde{\lambda})V(z(s',j'))+\rho|u|^p_\infty,\\
&\leq \left(1-\frac{\lambda}{2}\right)\frac{2\rho}{\lambda}|u|^p_\infty+\rho|u|^p_\infty=\frac{2\rho}{\lambda}|u|^p_\infty.
\end{align*}
Moreover, if $(s',j'+1)\in\text{dom}(z)$, then $z$ cannot leave $\Omega$ via flows because $\dot{V}\leq 0$ if $V(z)\geq \frac{2\rho}{\lambda}|u|_{\infty}^p$. It follows that for all $(s,j)\in \text{dom}(z)$ such that $s+j\geq T_{z,u,z_0}$ the solution $z$ satisfies:
\begin{equation}
\underline{\alpha}|z(s,j)|_{\mathcal{A}}^p\leq V(z(s,j))\leq \frac{2\rho}{\lambda}|u|^p_{\infty},
\end{equation}
that is, $|z(s,j)|_{\mathcal{A}}\leq \left(\frac{2\rho}{\lambda\underline{\alpha}}\right)^\frac{1}{p}|u|_{\infty}$, for all $s+j\geq T_{z,u,z_0}$. Combining this bound with \eqref{exponentialdecay} we obtain
\begin{equation}\label{maxExponentialDecay}
|z(s,j)|_{\mathcal{A}}\leq \max\left\{\left(\frac{\overline{\alpha}}{\underline{\alpha}}\right)^{\frac{1}{p}}|z(0,0)|e^{-\frac{\lambda}{2p}(s+j)},\left(\frac{2\rho}{\lambda\underline{\alpha}}\right)^\frac{1}{p}|u|_{\infty}\right\},
\end{equation}
for all $(s,j)\in\text{dom}(z)$. Since $\max\{a,b\}\leq a+b$, we obtain  
\begin{equation}\label{ISSbound}
|z(s,j)|_{\mathcal{A}}\leq \kappa_1|z(0,0)|e^{-\kappa_2(s+j)}+\kappa_3|u|_\infty,
\end{equation}
with $\kappa_1=\left(\frac{\overline{\alpha}}{\underline{\alpha}}\right)^{\frac{1}{p}}$, $\kappa_2=\frac{\lambda}{2p}$ and $\kappa_3=\left(\frac{2\rho}{\lambda\underline{\alpha}}\right)^{\frac{1}{p}}$. The result follows from the above inequality by time-invariance and causality.$\blacksquare$

\vspace{0.2cm}\noindent
The following result relaxes the third condition in Lemma \ref{auxlemma01} under a standard average dwell-time condition on the jumps. The proof follows similar ideas, but we present the additional steps.

\vspace{0.1cm}
\begin{lemma}\label{auxlemma012}
Consider the HDS \eqref{HDS0}, and suppose that: (a) every solution satisfies the ADT constraint \eqref{ADT00}; (b) there exist constants $\underline{\alpha},\overline{\alpha},\rho,p>0$, $\lambda\in(0,1)$, and a smooth function $V:C\cup D\to\mathbb{R}_{\geq0}$, such that the following inequalities hold:
\begin{align*}
\underline{\alpha}|z|^p_{\mathcal{A}}&\leq V(z)\leq \overline{\alpha}|z|^p_{\mathcal{A}},~~~~\forall~z\in C\cup D\cup G(D),\\
\langle\nabla V(z),F(z,u)\rangle&\leq -\lambda V(z)+\rho|u|^p,~~\forall~(z,u)\in C\times\mathbb{R}^m,\\
V(G(z))-V(z)&\leq 0,~~~~~~~~~~~~~~~~~~~~\forall~z\in D.
\end{align*}
Then, every solution of \eqref{HDS0} satisfies
\begin{equation}
|z(s,j)|_{\mathcal{A}}\leq \kappa_1e^{-\kappa_2(s+j)}|z(0,0)|_{\mathcal{A}}+\kappa_3\sup_{0\leq \tau\leq s}|u(\tau)|,
\end{equation}
for all $(s,j)\in\text{dom}(z)$, where $\kappa_i>0$, for $i\in\{1,2,3\}$. \QEDB
\end{lemma}
\vspace{0.1cm}\noindent
\textbf{Proof:} The proof follows similar steps as the proof of Lemma \ref{auxlemma01}. In particular, inequality \eqref{decreaseflowsproofA} still holds. On the other hand, during jumps, we now have
\begin{equation}
V(z(s_{n+1},n+1))-V(z(s_{n+1},n))\leq 0
\end{equation}
Dividing both sides by $\tilde{\lambda}V(z(s_{n+1},n))$, we obtain
\begin{align*}
0&\geq \frac{V(z(s_{n+1},n+1))-V(z(s_{n+1},n))}{\tilde{\lambda}V(z(s_{n+1},n))}\\
&=\int_{V(z(s_{n+1},n))}^{V(z(s_{n+1},n+1))}\frac{d\varrho}{\tilde{\lambda} V(z(s_{n+1},n))}.
\end{align*}
It follows that inequality \eqref{inequaintegral} now becomes $\int_{V(z(0,0))}^{V(z(s,j))}\frac{d\rho}{\tilde{\lambda}\varrho}\leq -s$, from which we obtain after integration:
\begin{equation}\label{boundStep222}
V(z(s,j))\leq V(z(0,0))e^{-\frac{\lambda}{2}s}
\end{equation}
Finally, the ADT condition \eqref{ADT} guarantees that $j\leq \frac{1}{\tau_d}s+N_0$ for any $(s,j)\in\text{dom}(\hat{z})$, which implies that $s+j\leq (\frac{1}{\tau_d}+1)s+N_0$. In turn, this inequality can be written as $s\geq\frac{\tau_d}{1+\tau_d}(s+j)-\frac{\tau_d}{1+\tau_d}N_0$, so that \eqref{boundStep222} can be upper-bounded as follows:
\begin{equation}\label{hybridexponentialdecrease2}
V(z(s,j))\leq \kappa_7e^{-\kappa_8 (s+j)}V(z(0,0)),
\end{equation}
where $\kappa_7:=e^{\frac{\lambda}{2} \frac{\tau_d}{1+\tau_d}N_0} $ and $\kappa_8:=\frac{\lambda}{2} \frac{\tau_d}{1+\tau_d}$. From here the proof follows the same Steps 3-4 from the proof of Lemma \ref{auxlemma01}.  
In particular,  the inequality \eqref{ISSbound} now becomes
\begin{equation*}
    |z(s,j)|_{\mathcal{A}}\leq \tilde{\kappa}_1|z(0,0)|e^{-\tilde{\kappa}_2(s+j)}+\tilde{\kappa}_3|u|_\infty,
\end{equation*}
with
$
    \tilde{\kappa}_1 \coloneqq \left(\frac{\overline{\alpha}}{\underline{\alpha}}\right)^{\frac{1}{p}}e^{\frac{\lambda}{2p} \frac{\tau_d}{1+\tau_d}N_0}$, $\tilde{\kappa}_2 \coloneqq \frac{\lambda}{2p}\frac{\tau_d}{1+\tau_d}$, and $\kappa_3 = \left(\frac{2\rho}{\lambda\underline{\alpha}}\right)^{\frac{1}{p}}.$
\hfill $\blacksquare$
%

\vspace{0.2cm}
\begin{lemma}\label{lemmaiss2}
Suppose that every solution pair $(\hat{z},\hat{u})$ of the HDS \eqref{hybridautomaton} satisfies the bound \eqref{exponential_bound0B} for all $(s,j)\in\text{dom}(\hat{z})$. Assume that $\Delta(\hat{\mu}_k)=\hat{\mu}_k^{-\ell}$, where $\ell>0$. Then, there exists $\beta_k(r,s)=r\min\{e^{-\kappa_2s},\mu_k^{-\ell}(s)\}\in\mathcal{KL}$, such that $(\hat{z},\hat{u})$ satisfies
\begin{equation*}
|\hat{z}(s,j)|_{\mathcal{A}}\leq \beta_k\Big(\bar{\kappa}_1|\hat{z}(0,0)|_{\mathcal{A}}e^{-\bar{\kappa}_2(s+j)}
                        +\quad \bar{\kappa}_{3,\mu_0}|\hat{u}|_{(s,j)}, s\Big),
\end{equation*}
for all $(s,j)\in\text{dom}(\hat{z})$, and where $\bar{\kappa}_1 \coloneqq \kappa_1^2$, $\kappa_2\coloneqq \frac{\kappa_2}{2}$, $ \bar{\kappa}_{3,\mu_0} \coloneqq \kappa_1\kappa_3(\mu_0 + 1)$.\QEDB
\end{lemma}

\vspace{0.1cm}
\noindent\textbf{Proof:} Consider a complete solution pair $(\hat{z},\hat{u})$ of the HDS \eqref{hybridautomaton} satisfying the bound \eqref{exponential_bound0B}. Then, we have that
\begin{align}\label{evaluatingthisbound}
|\hat{z}(s,j)|_{\mathcal{A}}\leq \kappa_1 e^{-\kappa_2 (s+j)}|\hat{z}(0,0)|_{\mathcal{A}}+\kappa_3  \cdot \sup_{0\leq \zeta\leq s}|\hat{\Delta}(\zeta)|,
\end{align}
for all $(s,j)\in\text{dom}(\hat{z})$, and where $\hat{\Delta}(s) \coloneqq \Delta(\mu_k^{-\ell}(s))\hat{u}(s)$.
Next, pick an arbitrary time $(\bar{s},\bar{j})\in\text{dom}(\hat{z})$, and let $\hat{y}(r,k):=\hat{z}(r+\bar{s},k+\bar{j})$,  and $v(r,k):=\mu_{k}^{-\ell}(\bar{s}+r)$. Since $\hat{y}$ is also a hybrid arc that is a solution to \eqref{hybridautomaton}, using the above bound and by time-invariance, it satisfies:
\begin{align}\label{expoboundA}
|\hat{y}(r,k)|_{\mathcal{A}}&\leq \kappa_1|y(0,0)|e^{-\kappa_2(r+k)} +  \kappa_3|\hat{u}|_{(r,k)}|v|_{(r,k)}\notag\\
                            &\!\!\!= \kappa_1|z(\bar{s},\bar{j})|e^{-\kappa_2(r+k)} +  \kappa_3|\hat{u}|_{(r,k)}\sup_{0\le \tau \le r} \hat{\mu}^{-\ell}(\bar{s} + \tau)\notag\\
                            &\le \kappa_1|z(\bar{s},\bar{j})|e^{-\kappa_2(r+k)} +  \kappa_3|\hat{u}|_{(r,k)}\hat{\mu}_k^{-\ell}(\bar{s}).
\end{align}
Now, using \eqref{evaluatingthisbound} with $s=\bar{s}$ and $j=\bar{j}$, we obtain:
\begin{equation}\label{expoboundB}
|\hat{z}(\bar{s},\bar{j})|_{\mathcal{A}}\leq \kappa_1|\hat{z}(0,0)|_{\mathcal{A}}e^{-\kappa_2(\bar{s}+\bar{j})} + \kappa_3|\hat{u}|_{(\bar{s},\bar{j})}\sup_{0\le \tau \le r}\hat{\mu}_k^{-\ell}(\tau) .
\end{equation}
Combining \eqref{expoboundA} and \eqref{expoboundB}, and using Remark 2, we have
\begin{align*}
|\hat{y}(r,k)|_{\mathcal{A}}&\leq \kappa_1\Big( \kappa_1|\hat{z}(0,0)|_{\mathcal{A}}e^{-\kappa_2(\bar{s}+\bar{j})}\\
                &\qquad\qquad + \kappa_3\sup_{0\le \tau \le r}|\hat{u}(\tau)|\sup_{0\le \tau \le r}\hat{\mu}_k^{-\ell}(\tau) \Big )e^{-\kappa_2(r+k)}\\
                &\qquad\qquad\quad  + \kappa_3\sup_{0\le \tau \le r}|\hat{u}(\tau)|\hat{\mu}_k^{-\ell}(\bar{s}).
\end{align*}
Evaluating the above bound at $r=\bar{s}$ and $\tilde{j}\in\mathbb{Z}_{\geq0}$ such that $(\bar{s},\tilde{j})\in\text{dom}(y)$, we obtain:
\begin{align*}
|\hat{y}(\bar{s},\tilde{j})|_{\mathcal{A}}&\leq \kappa_1\Big( \kappa_1|\hat{z}(0,0)|_{\mathcal{A}}e^{-\kappa_2(\bar{s}+\bar{j})}\\
                &\qquad + \kappa_3\sup_{0\le \tau \le \bar{s}}|\hat{u}(\tau)|\sup_{0\le \tau \le \bar{s}}\hat{\mu}_k^{-\ell}(\tau) \Big )e^{-\kappa_2(\bar{s}+\tilde{j})}\\
                &\qquad\quad  + \kappa_3\sup_{0\le \tau\le \bar{s}}|\hat{u}(\tau)|\hat{\mu}_k^{-\ell}(\bar{s}).
\end{align*}    
Using the definition of $\hat{y}$, and letting $\lambda\coloneqq 2\bar{s}$, $i\coloneqq\tilde{j}+\bar{j}$:
 \begin{align*}
|\hat{z}(\lambda, i)|_{\mathcal{A}}&\leq \kappa_1\Big( \kappa_1|\hat{z}(0,0)|_{\mathcal{A}}e^{-\kappa_2(\frac{\lambda}{2} + \bar{j})}\\
                &\quad + \kappa_3\sup_{0\le \tau \le \lambda/2}|\hat{u}(\tau)|\sup_{0\le \tau \le \lambda/2}\hat{\mu}_k^{-\ell}(\tau) \Big )e^{-\kappa_2\left(\frac{\lambda}{2}+\tilde{j}\right)}\\
                &\quad\quad  + \kappa_3\sup_{0\le \tau \le \lambda/2}|\hat{u}(\tau)|\hat{\mu}_k^{-\ell}\left(\lambda/2\right).
\end{align*}
Let $\eta_k(\lambda)\coloneqq\min\left\{e^{-\kappa_2\lambda},\mu_{k}^{-\ell}(\lambda)\right\}$, which is continuous and satisfies $\eta(\lambda)\to 0 $ as $\lambda\to\infty$. Then, for all $\tilde{j}\in\mathbb{Z}_{\ge0}$ we obtain
\begin{align}\label{lemmaISS+C:preConclusion}
    |\hat{y}(\lambda, i)&|_{\mathcal{A}}\leq \kappa_1\Big( \kappa_1|\hat{z}(0,0)|_{\mathcal{A}}e^{-\kappa_2(\frac{\lambda}{2} + i)}\notag\\
                &\quad + \kappa_3\sup_{0\le \tau \le \lambda/2}|\hat{u}(\tau)|\sup_{0\le \tau \le \lambda/2}\hat{\mu}_k^{-\ell}(\tau)e^{-\kappa_2\tilde{j}}\notag\\
                &\quad\quad  + \kappa_3\sup_{0\le \tau \le \lambda/2}|\hat{u}(\tau)|\Big )\eta_k(\lambda/2)\notag\\
                &\le \Big( \kappa_1^2|\hat{z}(0,0)|_{\mathcal{A}}e^{-\kappa_2(\frac{\lambda}{2} + i)} \notag\\
                &\qquad\qquad + 2\kappa_1\kappa_3\sup_{0\le \tau \le \lambda/2}|\hat{u}(\tau)|\Big )\eta_k(\lambda/2),
\end{align}
where we used the fact that $\sup_{0\le \tau \le \lambda/2}\hat{\mu}_k^{-\ell}(\tau)\le \mu_0^{-\ell} \le 1$ since $\mu_0\ge 1$ and $\ell>0$.
Since the choice of $(\bar{s},\bar{j})\in\text{dom}(z)$ was arbitrary, $z$ is complete, and the previous inequality holds for all $\tilde{j}\in\mathbb{Z}_{\ge 0}$, in particular we can use $s=2\bar{s}$, $j=\bar{j}$, and $\tilde{j}=0$ such that $(s,j)\in\text{dom}(z)$. Thus, from \eqref{lemmaISS+C:preConclusion} we obtain that there exists $\beta_k(r,s) \coloneqq r\cdot\eta_k(s)\in\mathcal{KL}$ such that
\begin{equation*}
|\hat{z}(s,j)|_{\mathcal{A}}\leq \beta_k\Big(\bar{\kappa}_1|\hat{z}(0,0)|_{\mathcal{A}}e^{-\bar{\kappa}_2(s+j)}
                        + \bar{\kappa}_{3}|\hat{u}|_{(s,j)}, s\Big),
\end{equation*}
with $\bar{\kappa}_1 \coloneqq \kappa_1^2$, $\kappa_2\coloneqq \frac{\kappa_2}{2}$, $ \bar{\kappa}_{3} \coloneqq 2\kappa_1\kappa_3$, and where we used Remark 2. \hfill $\blacksquare$
\subsection{Proofs of Section \ref{sectionapplications}}
\vspace{0.1cm}
\noindent \textbf{Proof of Propostion \ref{Prop:Regulation_Switching_Plants}:} The closed-loop system is given by $\dot{x} = \mu_k(t)A_{\sigma(t)}x$, which has the form \eqref{actualflow1}. In particular, we have that $f_q(x, u, \tau(t)) = A_qx$ which satisfies Assumption \ref{assumptionLipschitz}. Moreover, since $A_q$ is Hurwitz, for every $H\in\mathbb{R}^{n\times n},~H\succ 0$, there exists $P_q\succ 0$ satisfying $P_qA_q + A_q^\top P_q = -H$. From there, fix $H\succeq 0$, and for every $q$ consider the corresponding $P_q$. Then, for every $q\in\mathcal{Q}$ and all $\tau_d>0$ there exists $V_q(x) = x^\top P_q x$, $c_{q,1} = \lambda_{\min}(P_q),~c_{q,2}=\lambda_{\max}(P_q)$, $c_{q,3}=\lambda_{\min}\left(P_q\right)/\lambda_{\max}(P_q)$, and $c_{q,4}=0$ satifying the items of Assumption 4. By using $\tau_d>\ln(\max_{q\in\mathcal{Q}}c_{q,2}/\min_{q\in\mathcal{Q}}c_{q,1})/\min_{q\in\mathcal{Q}}c_{q,3}$, we obtain the result via Theorem \ref{theorem1}-a).\hfill$\blacksquare$

\vspace{0.1cm}
\noindent\textbf{Proof of Proposition \ref{Prop:Regulation_intermittenUpsilonlants}:} The interconnection of \eqref{application2:dynamics} and \eqref{feedback} results in a switching system of the form
\begin{align}
    \dot{x} &= \mu_k(t)\cdot\phi_{\sigma(t)}(x,1/\mu_k(t)),\label{application2:closedLoop}
\end{align}
where, for every $q\in\mathcal{Q}$, $\phi_{q}:\mathbb{R}^n\times\mathbb{R}\to \mathbb{R}_{\ge 0}^n$ is given by
\begin{align*}
    \phi_{q}(x,\upsilon) \coloneqq -\mathbb{I}_{\mathcal{Q}_s}\left(q\right) \left(\eta_q + \delta_q\psi_q(x)^2\right)b_q(x) x + \upsilon d_q(x).
\end{align*}
We show that, under Assumptions \ref{PTNES:gameRegularity}-\ref{PTNES:tuningGuidelines}, a suitable Lyapunov function can be used to show that Assumption \ref{assumptionunstablemodes} is satisfied. First, we analyze the R-Switching system $\dot{x} = \mu_k(t)\cdot \phi_{\sigma(t)}(x,\upsilon)$, with $\upsilon\in\mathbb{R}_{\ge0}$. Due to the assumptions on $b_q(\cdot)$, there exists $\underline{\sigma}_{q}>0$ such that $2\underline{\sigma}_{q} \le \lambda_{\min}\left(b_{q}(x) + b_{q}(x)^\top\right)$ holds for all $x\in\mathbb{R}^n$. Now, let $V_{q}(x)=\frac{1}{2\underline{\sigma}_{q}}|x|^2$ for every $q\in\mathcal{Q}_s$. By employing Young's inequality, we get:
\begin{align}
    \left\langle\nabla V_{q}(x) ,~\phi_{q}(x,\upsilon)\right\rangle 
                 \le   -&2\underline{\sigma}_q\eta_{q}V_{q}(x) + \upsilon^2\frac{1}{4\underline{\sigma}_{q}^2 \delta_{q}},
                \label{lieDerivativeRegulationUncertainty:stable}
\end{align}
for all $q\in\mathcal{Q}_s$ and $\upsilon\in\mathbb{R}_{\ge 0}$.
Similarly, for all $q\in\mathcal{Q}_u$ let $V_{q}(x)=\frac{|x|^2}{2}$. This function satisfies
\begin{align}
    \left\langle\nabla V_{q}(x) ,~\phi_{q}(x,\upsilon)\right\rangle %
                                                                          &\le V_{q}(x) + \upsilon^2\frac{\bar{d}_q^2}{2},\label{lieDerivativeRegulationUncertainty:unstable}
\end{align}
for all $q\in\mathcal{Q}_u,~\upsilon\in\mathbb{R}_{\ge 0}.$  Also, inequalities \eqref{lieDerivativeRegulationUncertainty:stable} and \eqref{lieDerivativeRegulationUncertainty:unstable} hold when $\upsilon = 1/\mu_k$. In this case, we have that
\begin{align*}
    &\left\langle\nabla V_{q}(x) ,~\phi_{q}(x,1/\mu_k)\right\rangle \le {-}2\underline{\sigma}_q\eta_{q}V_{q}(x) {+} \frac{1}{\mu_k^2}\frac{1}{4\underline{\sigma}_{q}^2 \delta_{q}},~q\in\mathcal{Q}_s,\\
    & \left\langle\nabla V_{q}(x) ,~\phi_{q}(x,1/\mu_k)\right\rangle \le  V_{q}(x) {+} \frac{1}{\mu_k^2}\frac{\bar{d}_q^2}{2},~q\in\mathcal{Q}_u.
\end{align*}
Using $c_{q,1}=c_{q,2}= 1/2\underline{\sigma}_q$, $c_{q,3} = 2\underline{\sigma}_{q}\eta_{q}$, $c_{q,4}=1/4\underline{\sigma}_{q}^2\delta_{q}$, when $q\in \mathcal{Q}_s$, and 
$c_{q,1}=c_{q,2}=1/2$, $c_{q,5} = 1$, $c_{q,4}= \bar{d}_q^2/2$ when $q\in\mathcal{Q}_u$, together with the set of smooth functions $\{V_{q}\}_{q\in\mathcal{Q}}$, Assumption \ref{assumptionunstablemodes} is satisfied. 
Thus, we can always pick $\tau_a>1$ and $\tau_d>0$ large enough to satisfy the stability condition \eqref{conditionstability2}.
Moreover, Assumption 3 is satisfied by the Lipschitz assumption on $d_{q}(\cdot)$ and $b_{q}(\cdot)$. Using these facts, we obtain the result via Theorem 2-c). \strut\hfill$\blacksquare$

\vspace{0.1cm}
%
%
The following Lemma is instrumental to study the stability properties of the HDS with data \eqref{pT:NESwMomentum}.

\begin{lemma}\label{lem:Mc0q}
Consider the matrix
\begin{align}\label{matrixM}
    M_{\zeta_q}(x_1,\tau) \coloneqq \begin{pmatrix}
        \frac{1}{\eta(\tau)^2}I & I -\partial\mathcal{G}_q(x_1)^\top\\
        I - \partial\mathcal{G}_q(x_1) & (\zeta_q-\rho\eta'(\tau))I
    \end{pmatrix},
\end{align}
where $q\in\mathcal{Q}$, $\tau\in[0,N_0]$, $\eta(\tau)\in [\underline{\eta}, \overline{\eta}]$, $\rho\in[0, 1/\tau_d]$, and $\eta'(\tau)\coloneqq\frac{d\eta}{d\tau}(\tau)$, $\mathcal{G}_q(\cdot)$, and $\zeta_q$ are as introduced in Section \ref{resetsexample}. Suppose that Assumptions \ref{PTNES:gameRegularity} and \ref{PTNES:tuningGuidelines} are satisfied. Then, 
\begin{equation}\label{eq:Mc0_bound}
    M_{\zeta_q}(x_1,\tau) \succeq \nu_M I,~~~\forall~\tau\in[0,N_0],~x_1\in\mathbb{R}^n
\end{equation}
where 
 $   \nu_M{\coloneqq}\frac{(1-\delta_d - \delta_\eta)\overline{\sigma}^2}{\delta_{\eta}(1-\delta_d)\underline{\zeta} + \overline{\sigma}^2}$, with $\underline{\zeta}\coloneqq \min_{q\in\mathcal{Q}}\zeta_q$ and $\overline{\sigma}\coloneqq \max_{q\in\mathcal{Q}}\sigma_{q}$. \QEDB
\end{lemma}

\vspace{0.2cm}
\noindent\textbf{Proof:} 
The proof follows the same ideas of \cite[Lemma B.3]{decentralizedMomentumResets23}. First we show that matrix-valued function $M_{\zeta_q}(\cdot,\cdot)$ is positive-definite uniformly over $\rho\in[0,\tau_d^{-1}]$, $x_1\in\mathbb{R}^n$, and $\tau\in[0,N_0]$. To this end, we decompose the matrix $M_{\zeta_q}(x_1,\tau)$ as follows:
\begin{subequations}    
\begin{align}\label{eq:Mc0_decompose}
&M_{\zeta_q}(x_1,\tau)= U_q(x_1,\tau) W_q(\tau,x_1)U_q(x_1,\tau)^\top,\\
&W_q(\tau,x_1)\!\coloneqq\!\begin{pmatrix} \frac{I}{\eta(\tau)^{2}} &  0  \\ 0 &\!\!\!\!\! \varrho_q(\tau)I \!-\! \eta^2(\tau)\Sigma_{q}(x_1)\Sigma_q(x_1)^{\!\top}\!\!\end{pmatrix}\!,\\
&\varrho_q(\tau)\coloneqq \zeta_q\!-\!\rho \eta'(\tau),\quad \Sigma_q(x_1)\coloneqq I-\partial \mathcal{G}_q(x_1),\\
&U_q(x_1,\tau)\coloneqq \begin{pmatrix}
    I & 0\\
    \eta^2(\tau)\Sigma_q(x_1) & I
\end{pmatrix}.
\end{align}
\end{subequations}
By the fact that $\eta(\tau)\in[\underline{\eta},\overline{\eta}]$ for all $\tau\in[0,N_0]$ it follows that 
\begin{equation}\label{lemma:PositiveDefinite:Diag1}
    \frac{1}{\eta(\tau)^2}I\succeq \frac{1}{\overline{\eta}^2}I.
\end{equation}
Also,  by Assumptions \ref{PTNES:gameRegularity} and \ref{PTNES:tuningGuidelines}, we have that
\begin{align}\label{lemma:PositiveDefinite:Diag2}
    \varrho_q(\tau)I-\eta(\tau)^2\Sigma_q(x_1)\Sigma_q(x_1) &\succeq \left(\varrho_q(\tau)-\overline{\eta}^2\sigma_q^2\right)I\notag\\
    &\succeq \left(\zeta_q-\frac{\overline{\eta}-\underline{\eta}}{\tau_d N_0} - \overline{\eta}^2\sigma_q^2\right)I\notag\\
    &\succeq \tilde{\delta} I,
\end{align}
where $\tilde{\delta}\coloneqq  (1-\delta)\underline{\zeta}$, with $\underline{\zeta}\coloneqq \min_{q\in\mathcal{Q}} \zeta_q$. Therefore, via \cite[Theorem 7.7.7]{horn2012matrix},  the matrix  $M_{\zeta_q}(x_1,\tau)$ is positive definite for all $x_1\in\mathbb{R}^n$ and $\tau\in[0,N_0]$.
Now, we establish the matrix inequality \eqref{eq:Mc0_bound}. To do so, we use \eqref{lemma:PositiveDefinite:Diag1} and \eqref{lemma:PositiveDefinite:Diag2} in \eqref{eq:Mc0_decompose} to obtain that
\begin{align}
		M_{\zeta_q}(x_1,\tau)&\succeq U_q(x_1,\tau) 
  \begin{pmatrix}
 \!\frac{1}{\overline{\eta}^2}I&0\\
  0 &\tilde{\delta}I
  \end{pmatrix}U_q^\top(x_1,\tau)\notag\\
  &\succeq Z_q(x_1,\tau)Z_q(x_1,\tau)^\top,\label{lemma:PositiveDefinite:Vmatrix}
\end{align}
where $Z_q(x_1,\tau)^\top$ is the upper block triangular matrix
\begin{align*}
Z_q(x_1,\tau)^\top\coloneqq \begin{pmatrix} \frac{1}{\overline{\eta}}I & \frac{\eta(\tau)^2}{\overline{\eta}}\Sigma_q(x_1,\tau)^\top \\  0 & \sqrt{\tilde{\delta}}I \end{pmatrix}.
\end{align*}
By  applying \cite[Lemma B.2]{decentralizedMomentumResets23}, and using \eqref{lemma:PositiveDefinite:Vmatrix} together with the fact that $Z_{q}(x_1,\tau)$ has full column rank for all $x_1\in\mathbb{R}^n$ and $\tau\in[0,N_0]$ and thus that $\sigma_{\min}(Z_{q}(x_1,\tau)Z_{q}(x_1,\tau)^\top)\ge\sigma_{\min}(Z_{q}(x_1,\tau))\sigma_{\min}(Z_{q}(x_1,\tau)^\top) = \sigma^2_{\min}(Z_{q}(x_1,\tau)^\top)$, we obtain
\begin{align*}
M_{\zeta_q}(x_1,\tau)& \succeq \frac{1}{\overline{\eta}^2\left(1 + \frac{\overline{\eta}^2}{\tilde{\delta}}\lVert \Sigma_q(x_1)^2\rVert\right) + \frac{1}{\tilde{\delta}}}I\\
    &\succeq \frac{(1-\delta)\underline{\zeta}}{\overline{\eta}^2((1-\delta)\underline{\zeta} + \delta_{\eta}\underline{\zeta})+1}I\\
    &\succeq \frac{(1-\delta_d - \delta_\eta)\overline{\sigma}^2}{\delta_{\eta}(1-\delta_d)\underline{\zeta} + \overline{\sigma}^2}I,
\end{align*}	
where in  the last two steps we used Assumption \ref{PTNES:tuningGuidelines}. This completes the proof.\hfill$\blacksquare$

\vspace{0.1cm}
\noindent \textbf{Proof of Proposition \ref{Prop:Nash_equilibrium}:} 
For every $q\in\mathcal{Q}$ consider the Lyapunov function
\begin{align*}
    V_q(x,\tau) &= \frac{1}{4}|x_2 - x^*|^2 + \frac{1}{4}|x_2-x_1|^2 + \frac{\eta(\tau)^2}{2}\left|\mathcal{G}_q(x_1)\right|^2,
\end{align*}
which in the flow set and jump set satisfies:
\begin{align}\label{application3:quadraticBounds}
 v_{q,1}|x|_{\mathcal{A}}^2 \le V_q&(x,\tau) \le v_{q,2}|x|_{\mathcal{A}}^2,
\end{align}
where $v_{q,1},v_{q,2}>0$ are given by $v_{q,1}\coloneqq 0.25\min\left\{1, 2\kappa_q^2\underline{\eta}^2\right\}$, and $v_{q,2}\coloneqq 0.25\max\left\{3,~2 + 2\ell_q^2\overline{\eta}^2\right\}$. Let
\begin{align*} \mathcal{L}_{\left(f_q,\rho\right)}V_q(x,\tau)\coloneqq \left\langle \nabla V_q(x,\tau), \left(\begin{array}{c}
f_q(x,\tau)\\
\rho
\end{array}\right)\right\rangle,    
\end{align*}
where $\rho\in [0,\tau_d^{-1}]$. Since  $\mathcal{G}_q(\cdot)$ is  $\kappa_q$-strongly-monotone and $\ell_q-$Lipschitz,  we have that $\left\langle x_1-\tilde{x}_1,~ \mathcal{G}_q(x_1)\right\rangle \ge \zeta_q\left|\mathcal{G}_q(x_1)\right|^2$, with $\zeta_q= \kappa_q^2/\ell_q$. Thus during flows, for all $(x,\tau, \rho) \in \mathbb{R}^{2n}\times [0, N_0]\times [0,\tau_d^{-1}]$, we have:
\begin{align}
  \mathcal{L}_{(f_q,\rho)}V_q(x,\tau)  &= -\frac{1}{\eta(\tau)}|x_2-x_1|^2\notag\\
   &~ - 2\eta(\tau)\left\langle\mathcal{G}_q(x_1),~\left[I - \partial\mathcal{G}_q(x_1)\right](x_2-x_1)\right\rangle\notag \\
   &~ -\eta(\tau)\left[\left\langle x_1 - x^*,~\mathcal{G}_q(x_1)\right\rangle  -  \rho\eta'(\tau)|\mathcal{G}_q(x_1)|^2\right]\notag\\
    &\le -\eta(\tau)\left\langle\chi_q,~M_{\zeta_q}(x_1,\tau)\chi_q\right\rangle,\label{PTNES:LieDerivative}
\end{align}
where $\chi_q\coloneqq \left(x_2-x_1, \mathcal{G}_q(x_1)\right)\in\mathbb{R}^{2n}$, and $M_{\zeta_q}$ is given by \eqref{matrixM}. Using Lemma \ref{lem:Mc0q}, we conclude: 
\begin{align}\label{application3:LieDerivative}
    \mathcal{L}_{(f_q,\rho)}V_q(x,\tau) &\le -\underline{\eta}\nu_{M}|\chi_q|^2
\end{align}
$\forall x_1\in\mathbb{R}^n,~\tau\in [0,N_0]$. Since $V_q(x,\tau)\leq \frac{1}{4}\max\left\{3, 2\left(\frac{1}{\kappa_q^2} + \overline{\eta}^2\right)\right\}|\chi_q|^2$,
using \eqref{application3:LieDerivative}, we obtain
\begin{align}\label{application3:LieDerivativefinal}
     \mathcal{L}_{(f_q,\rho)}V_q(x,\tau) \le -\frac{4\underline{\eta}\nu_{M}}{\max\left\{3, 2\left(\frac{1}{\kappa_q^2} + \overline{\eta}^2\right)\right\}}V_q(x,\tau).
\end{align}
Now, for all $p,q\in\mathcal{Q}$, let
\begin{align*}
    \Delta V_{p}^{q}(x,\tau)\coloneqq V_{q}\left(R_{p}(x),\tau - 1\right) - V_{p}(x,\tau),\quad \tau\in[1,N_0].
\end{align*}
Then, during jumps we have:
\begin{align}
     &\Delta V_{p}^q (x,\tau)  = V_q\left((x_1,x_1),\tau-1\right) - V_p(x,\tau)\\
                            &\le -\frac{1}{4}|x_1-x^*|^2- \frac{1}{4}|x_1-x_2|^2 + \frac{1}{4}|x_1 - x^*|^2  \notag\\
                             &\quad + \frac{1}{2}\bigg(\eta(\tau-1)^2\frac{\ell_{q}^{2} }{\kappa_{p}^2} - \eta(\tau)^2\bigg)\left|\mathcal{G}_p(x_1)\right|^2\notag\\
                             &\le  -\frac{1}{4}|x_1-x^*|^2- \frac{1}{4}|x_1-x_2|^2 + \frac{1}{4\kappa_p^2}|\mathcal{G}_p(x_1)|^2\notag\\
                             &\quad + \frac{1}{2}\left(\eta(N_0-1)^2\frac{\ell_{q}^{2} }{\kappa_{p}^2} - \eta(1)^2\right)\left|\mathcal{G}_p(x_1)\right|^2\notag\\
                            &\le-\left(1-\gamma_{p}^q\right)V_{p}(x,\tau),\notag
\end{align}
where $\gamma_{p}^q\coloneqq \frac{2\eta(N_0-1)^2\ell_{q}^{2}  + 1}{2\kappa_p^2\eta(1)^2}$. The above inequality implies that
\begin{align}\label{application3:preDecreaseJumps}
    V_{q}\left(R_p(x),~\tau - 1\right) \leq \gamma_{p}^q V_p(x,\tau).
\end{align}
where
\begin{align*}
    \gamma_{p}^q &\le  \frac{2\eta(N_0-1)^2\ell_{q}^{2}  + 1}{2\kappa_p^{2}\eta(1)^2}
                \le \frac{\overline{\ell}^2}{\underline{\kappa}^2}\frac{\eta(N_0-1)^2}{\eta(1)^2} + \frac{1}{2\underline{\kappa}^2 \eta(1)^2}\eqqcolon \overline{\gamma},
\end{align*}
where $\underline{\ell}\coloneqq \min_{q\in\mathcal{Q}}\ell_q$, $\overline{\kappa}\coloneqq \max_{q\in\mathcal{Q}}\kappa_q$, and $\underline{\kappa}\coloneqq \min_{q\in\mathcal{Q}}\kappa_q$,
Therefore, from \eqref{application3:preDecreaseJumps}, we obtain that
\begin{align}
    V_{q}\left(R_p(x),\tau-1\right) \leq \overline{\gamma} V_p(x,\tau),\label{application3:DecreaseJumps}
\end{align}
for all $\tau\in[1,N_0],~p,q\in\mathcal{Q}$.
By the smoothness properties of $\mathcal{G}_q(\cdot)$ and the differentiability of $\eta(\cdot)$ by design, we obtain that $f_{q}(x,\tau)$ is locally Lipschitz and, thus, that Assumption \ref{assumptionLipschitz} also holds. On the other hand, note that via a simple change of coordinates, and without loss of generality, the results of Theorem \ref{theorem1} hold for $\mathcal{A}=\tilde{\mathcal{A}}_1\times \mathbb{R}_{\ge 1}$, where the compact set $\tilde{\mathcal{A}}_1$ is the same as $\mathcal{A}_1$ in \eqref{stableset1} but with the set $\{0\}$ replaced by the set $\mathcal{A}_x$ defined in \eqref{application3:set}. Therefore, inequalities \eqref{application3:quadraticBounds}, \eqref{application3:LieDerivativefinal}, \eqref{application3:DecreaseJumps}, and the fact that $\overline{\gamma}\in(0,1]$, by assumption, imply PT- PT-S$_{\text{F}}$ stability of $\mathcal{A}_x\times [0,N_0]\times \mathcal{Q}\times \mathbb{R}_{\ge 1}$ via Theorem \ref{theorem1}-a). In particular, for any solution $z=(x,\tau,q,\mu)$ to the HDS \eqref{mainHDSmodel} with data \eqref{mapsstable01}, and maps $f_q$ and $R_q$ as defined in \eqref{pT:NESwMomentum}, the following bound holds:
\begin{align*}
|z(t,j)|_{\mathcal{A}}\leq \kappa_1 e^{-\kappa_2 (\mathcal{T}_k(t)+j)}|z(0,0)|_{\mathcal{A}},
\end{align*}
for any $(t,j)\in\text{dom}(z)$, and where
\begin{align*}
    \kappa_1 &= \frac{\left(\max_{q\in\mathcal{Q}}v_{q,2}\right)^{\frac{N_0 + 1}{2}}}{\left(\min_{q\in\mathcal{Q}}v_{q,1}\right)^{\frac{N_0}{2}}}e^{\frac{\lambda}{2} \frac{\tau_d}{1+\tau_d}N_0},~~\text{and}\\
    \kappa_2 &= \frac{\tau_d}{4(1+\tau_d)}\min\Bigg\{1-\overline{\gamma},\\
            &\quad~\frac{4\underline{\eta}\nu_M}{\max\left\{3, 2\left(\frac{1}{\underline{\kappa}^2} + \overline{\eta}^2\right)\right\}}-\frac{1}{\tau_d}\ln\left(\frac{\max\{3,~2+2\overline{\ell}^2\overline{\eta}^2\}}{\min\{1, 2\underline{\eta}\underline{\kappa}^2\}}\right)\Bigg\}.
\end{align*}
This completes the proof.
\strut\hfill$\blacksquare$

\end{document}